%% file: 0_main.tex
\documentclass[11pt]{article}
\usepackage[letterpaper,hmargin=1.0in,vmargin=1.0in]{geometry}
\usepackage{amsmath,amsthm,amssymb,setspace,sectsty,bbm,graphicx,mathtools}
\usepackage[round]{natbib}
\usepackage[usenames,dvipsnames]{xcolor}
\usepackage[linktocpage=true,pagebackref=true,colorlinks,allcolors=blue,bookmarks=true,bookmarksopen,bookmarksnumbered]{hyperref}

\usepackage{booktabs} 
\usepackage{algorithm,algorithmicx}
\usepackage[noend]{algpseudocode}
\usepackage{subcaption}
\usepackage{bbm}
\usepackage{tabularx}
\usepackage{multirow}
\usepackage{makecell}
\usepackage{xcolor}
\usepackage{bm}
\usepackage{multirow}

\newtheoremstyle{DStheorem}
{\topsep}
{\topsep}
{\itshape}
{0pt}
{\scshape}
{.}
{ }
{\thmname{#1}\thmnumber{ #2}\thmnote{ (#3)}}
\theoremstyle{DStheorem}
\newtheorem{theorem}{Theorem}[section]
\newtheorem{lemma}[theorem]{Lemma}
\newtheorem{claim}[theorem]{Claim}
\newtheorem{corollary}[theorem]{Corollary}

\newtheorem{remark}[theorem]{Remark}
\newtheorem{problem}[theorem]{Problem}
\newtheorem{proposition}[theorem]{Proposition}
\newtheorem{definition}[theorem]{Definition}

\let\oldproofname=\proofname
\renewcommand{\proofname}{\rm\sc{\oldproofname}}

\makeatletter
\newcommand{\multiline}[1]{%
  \begin{tabularx}{\dimexpr\linewidth-\ALG@thistlm}[t]{@{}X@{}}
    #1
  \end{tabularx}
}
\makeatother
\makeatletter
\def\munderbar#1{\underline{\sbox\tw@{$#1$}\dp\tw@\z@\box\tw@}}
\makeatother

\newcommand{\lrp}[1]{\left(#1\right)}

\newcommand{\lrc}[1]{\left[#1\right]}
\newcommand{\lrl}[1]{\left\{#1\right\}}
\newcommand{\lra}[1]{\left|#1\right|}

\newcommand{\R}{\mathbb{R}}

\newcommand{\E}{\mathbb{E}}
\newcommand{\PP}{\mathbb{P}}
\newcommand{\indicator}{\mathbbm{1}}

\newcommand{\OPT}{\text{OPT}}

\newcommand{\ind}[1]{\mathbbm{1}_{\lrl{#1}}}

\newcommand{\cB}{\mathcal{B}}

\newcommand{\cM}{\mathcal{M}}

\newcommand{\cR}{\mathcal{R}}

\newcommand{\asize}{K} 
\newcommand{\potentials}{\mathcal{P}}

\newcommand{\backlog}{\mathcal{B}}

\newcommand{\matches}{\cM}

\newcommand{\batch}{S}

\newcommand{\plike}{\phi}

\newcommand{\cIJ}{I\cup J}

\newcommand{\clikesin}[1]{\mathcal{L}_{#1, \leftarrow}}
\newcommand{\cdislikesin}[1]{\mathcal{D}_{#1, \leftarrow}}
\newcommand{\clikesout}[1]{\mathcal{L}_{#1, \rightarrow}}
\newcommand{\cdislikesout}[1]{\mathcal{D}_{#1, \rightarrow}}

\newcommand{\likesout}[1]{L_{#1, \rightarrow}}
\newcommand{\dislikesout}[1]{D_{#1, \rightarrow}}

\newcommand{\cumlikesin}[1]{\mathcal{L}_{#1, \leftarrow}}
\newcommand{\cumdislikesin}[1]{\mathcal{D}_{#1, \leftarrow}}
\newcommand{\cumlikesout}[1]{\mathcal{L}_{#1, \rightarrow}}
\newcommand{\cumdislikesout}[1]{\mathcal{D}_{#1, \rightarrow}}

\newcommand{\RR}{\mathbb{R}}
\newcommand{\Ma}{{\sf M}}

\newcommand{\probP}{\text{I\kern-0.15em P}}

\DeclareMathOperator*{\argmaxA}{arg\,max}
\DeclareMathOperator*{\argmax}{arg\,max}

\definecolor{myblue}{RGB}{51,51,178}
\definecolor{myred}{RGB}{178,51,51}
\definecolor{myblack}{RGB}{0,0,0}
\definecolor{mygreen}{RGB}{0,128,0}
\definecolor{mygray}{gray}{0.6}
\definecolor{mysilver}{RGB}{220,220,220}

\algrenewcommand\algorithmicrequire{\textbf{Input:}}
\algrenewcommand\algorithmicensure{\textbf{Output:}}

\sectionfont{\large} \subsectionfont{\normalsize}
\allowdisplaybreaks
\makeindex
\newcommand{\changelocaltocdepth}[1]{%
  \addtocontents{toc}{\protect\setcounter{tocdepth}{#1}}%
  \setcounter{tocdepth}{#1}%
}

\begin{document}
\begin{titlepage}

\title{The Dating Heuristic: \\ A Provably Strong Matching Algorithm for Dating Platforms}

\author{%
Ignacio Rios\thanks{ School of Management, University of Texas at Dallas. Email: {\tt ignacio.riosuribe@utdallas.edu}.}
\and Alfredo Torrico\thanks{ Daniel J. Epstein Department of Industrial and Systems Engineering, University of Southern California. Email: {\tt atorrico@usc.edu}.}}

\date{}
\maketitle

\setcounter{page}{200}
\thispagestyle{empty}

\begin{abstract}
Motivated by online dating platforms, we study the problem of selecting which subset of profiles to display to each user in each period. Users observe the profiles set by the platform, decide which of them to like, and a match occurs if and only if two users mutually like each other, potentially across different periods. The platform aims to maximize the expected number of matches produced over the entire time horizon, and users' behavior---captured by their like probabilities---may depend on their history.

  We develop a general theoretical model that captures the dynamic, two-sided nature of the problem and the influence of users’ past experiences on their future behavior. We focus on one-lookahead policies and propose the Integral Dating Heuristic (DH-int), providing formal performance guarantees: DH-int achieves a uniform $1-1/e$ guarantee across all platform designs under reasonable assumptions. Our empirical analysis, using proprietary data from a major U.S.-based dating app, confirms that DH-int consistently outperforms other benchmarks such as Greedy, Perfect Matching and DH, and approaches the theoretical upper bound across multiple platform designs and variants of the history effect. The superior performance of DH-int is driven primarily by its careful balancing of initial and follow-up interactions, which accounts for the two-sided nature of the market.

\end{abstract}

    
\end{titlepage}


\newpage
\pagestyle{plain}
\setcounter{page}{1}

\input{1_introduction}

\input{2_literature}
%
\input{3_problem}
\input{4_analysis_two_periods}

\input{5_analysis_multiple_periods}

\input{6_experiments}

\input{7_conclusions}

\addcontentsline{toc}{section}{Bibliography}
\bibliographystyle{plainnat}
\bibliography{draft} 

\changelocaltocdepth{1} 
\appendix

\input{9_online_supplement}

\input{10_electronic_companion}
\end{document}

%% file: 1_introduction.tex
\section{Introduction}\label{sec: introduction}

A common feature of many dating platforms is their \emph{curated} approach, presenting users with a limited subset of profiles rather than granting unrestricted access to all the available ones. This strategy, used in platforms such as Bumble, Hinge, Coffee Meets Bagel, and The League, aims to enhance users' experience by promoting more meaningful connections, reducing app fatigue, and fostering a sense of anticipation and excitement that can translate into user satisfaction.

A central operational decision for these platforms is determining which subset of profiles to present to each user in each period so as to maximize the total number of matches. This problem is challenging for several reasons. First, it is inherently two-sided: a match is only produced if both users see and like each other. Second, the platform must handle a large state space and complex dynamics, tracking the sets of users who have liked or disliked each profile in previous periods. Third, users’ past behavior---including the profiles they have liked or disliked, the matches they have obtained, and other aspects of their history---can influence their future choices, creating dynamic interactions across periods. For example, users who have recently experienced few successful matches may become frustrated, while users with many recent matches may become pickier, thereby affecting their future decisions in the app. Consequently, platforms must carefully manage the timing and sequence of profile displays to maximize match formation. 
These challenges are further shaped by the platforms' design choices and implementation constraints. For instance, while most platforms employ a two-directional approach that allows both sides to take the first action in the sequence of evaluations that could lead to a match, Coffee Meets Bagel uses a one-directional model---known as ``Ladies Choice’’---where men evaluate profiles first and women are then shown only the profiles of men who liked them.
Other design choices, such as explicitly indicating to a user that the profile they are viewing has previously liked them---a feature offered to paid users on most major apps---can increase the probability of reciprocation (as shown in~\cite{Lee2014}) and thus affect the overall match rate. Consequently, dating platforms must carefully account for the interplay between design features and user behavior when deciding which profiles to display.

A good policy---i.e., a decision rule specifying which profiles to display in each period---for these platforms should dynamically track the state of the system, including users’ backlogs, matches and prior interactions, and then use this information to select which profiles to display in the current period while accounting for the expected impact of these decisions on all future periods. Implementing such a policy, however, is extremely challenging. The problem is NP-hard~\citep{rios2021}, the state space grows combinatorially, and users' behavior may evolve in complex, path-dependent ways. As a result, platforms may be compelled to adopt simpler heuristic approaches or to focus only on the short-term consequences of their decisions. Yet, choosing among these alternatives requires a clear understanding of how they affect overall match generation and how their performance interacts with platform design choices and users' behavior.

In this paper, we theoretically and empirically study the robustness of different approaches to decide the subset of profiles to display to each user in each period to maximize the expected number of matches, whereby a match between two users realizes if, and only if, both users like each other, potentially in different periods.  
The goal of this paper is twofold. First, we develop a general theoretical model and provide performance guarantees for commonly used algorithms, offering insights into their effectiveness in dynamic two-sided settings. Second, we empirically examine the robustness of these methods to variations in platform design and user behavior, demonstrating how design choices and behavioral dynamics interact to influence match outcomes. By combining a thorough theoretical analysis with extensive simulations, our results provide both principled guidance and practical insights for optimizing curated dating platforms.


\subsection{Contributions}

We begin by introducing a theoretical model of a dating platform that must decide, in each period, which subset of profiles to show to each user in order to maximize the expected total number of matches. A match occurs only when both users like each other, potentially across different periods, and users’ past experiences can influence their future behavior. To capture this, we model like probabilities as a function of users’ histories, which include the profiles they have liked or disliked, the matches they have obtained, and any pending likes in their backlogs. The platform observes the state of the system—fully characterized by users’ histories and their sets of potentials, i.e., the set of users they can potentially match with based on preferences and prior evaluations—and uses this information to design a policy that determines which profiles to display to each user in each period.

Given the complexity of the problem, we focus on one-lookahead policies for selecting which profiles to display. These are policies that, at any period $t$, use the current state of the system to choose profiles by explicitly considering the impact of these decisions on both the current and the immediately following periods (i.e., $t$ and $t+1$). Despite this limited horizon, optimally implementing one-lookahead policies remains challenging due to the stochastic nature of users’ decisions and the non-linear dependence of the expected number of matches in the following period on the current period's choices. We overcome this challenge by using the current like probabilities in both periods of the horizon (i.e., $\phi_{\ell, \ell'}^{t+1} \approx \phi_{\ell, \ell'}^t$ for all pairs $\ell, \ell'$), which enable us to derive an upper bound of the problem. 
Similar to the Dating Heuristic (DH) introduced by~\cite{rios2021}, we use this upper bound as a building block for the Integral Dating Heuristic (DH-int). Unlike DH, which optimizes continuous variables and then rounds the fractional solution to determine which the profiles to display, DH-int considers integral decisions in the initial period of the one lookahead horizon, producing a solution that can be directly implemented and eliminating the need for rounding procedures that can create inefficiencies and complicate the theoretical analysis.

{As a starting point of the analysis, we study the performance of DH-int and several commonly used benchmarks under a two-period horizon. 
Our main result shows that DH-int achieves a uniform $1-1/e$ approximation guarantee. This is in sharp contrast to natural heuristics such as Greedy---which selects profiles based solely on immediate match probabilities and isolated from other users' choices---or Perfect Matching---which prioritizes clearing backlogs and, when none exist, pairs users to evaluate each other simultaneously---both of which can perform arbitrarily poorly in the worst case. Although algorithms from the submodular optimization literature provide constant-factor guarantees, we demonstrate that these are weaker and design-dependent, further highlighting the robustness and distinct advantage of DH-int.}
Furthermore, the derivation of our main performance guarantee for DH-int is of independent interest and builds on two key elements. First, we develop two upper bounds for the maximum expected number of matches that can be realized in the second period as a result of the profiles displayed in the first period. One bound is obtained by adapting DH-int, and the other by formulating a novel linear program---referred to as the \emph{distribution problem}---in which the decision variables represent the general distribution of backlogs. 
Using a duality argument, we show that the upper bound obtained by the latter is at least that obtained by the former.
Second, we demonstrate that, given the subsets of profiles displayed in the first period, the induced distribution of the backlog in the second period constitutes a feasible solution to the distribution problem. Hence, we can use the notion of \emph{correlation gap}~\citep{agrawal2010correlation}, which compares the objective values of our solution and the distribution problem, to establish our guarantee. 

In the general case of $T$ periods with $T>2$, the analysis becomes substantially more challenging. As discussed earlier, the state space grows quickly with $T$, the dynamics introduce intricate dependencies, and the impact of current decisions in far-ahead periods is difficult to capture. Despite these challenges, we establish theoretical results for the general case under some additional assumptions. On one hand, we show that when like probabilities are time-invariant and the platform restricts attention to \emph{semi-adaptive} policies---those that select profiles for the initiating side in a non-adaptive manner but adaptively select profiles from the backlog for the responding side---a variant of DH-int designed for the multi-period setting, which we call DHT, achieves the same $1-1/e$ approximation guarantee (see Table~\ref{tab: summary comparison algorithms} in Appendix~\ref{app:intro} for a summary of the differences between our methods). On the other hand, we prove that this $1-1/e$ guarantee also holds relative to \emph{any} adaptive policy in the special case where the platform design is one-directional and only sequential interactions are allowed. These results highlight that even in the more complex multi-period environment, a semi-adaptive variant of DH-int can provide strong performance guarantees under reasonable assumptions. 

To complement our theoretical analysis, we empirically evaluate the performance of DH-int and the other benchmarks using proprietary data from our industry partner, a major U.S.-based dating app (kept anonymous as part of our NDA). Since our goal is to test the robustness of these methods, our simulations consider different platform designs---e.g., two and one-directional designs---and several variants and magnitudes of the history effect. These variants include no history dependence, histories driven by the number of recent matches, boosts in like probabilities when users know that a profile comes from their backlog, and potential disengagement effects as users spend longer times using the app. Across all these settings, DH-int consistently outperforms all other benchmarks. Moreover, DH-int performs remarkably close to the theoretical upper bound, suggesting that its practical effectiveness exceeds its formal guarantees. We also find that DH-int yields a modest but consistent improvement over DH, indicating that the rounding step required by the latter introduces some inefficiencies. Our analysis also reveals that the main driver of the superior performance achieved by the DH variants (i.e., DH-int and DH)  is their ability to carefully balance the like probabilities across initial and follow-up interactions. Finally, we find that the non-adaptive variant of DH-int, DHT, performs competitively relative to DH-int. This makes DHT an attractive option for platforms where computational speed and resource efficiency are critical, offering a practical balance between performance and implementability.

Overall, our results demonstrate---both theoretically and empirically---that DH-int is a robust and reliable approach for deciding which profiles to display in curated dating platforms. Beyond its robustness, DH-int offers a simple yet powerful framework that performs near-optimally in practice across diverse platform designs and user behaviors. Together, these findings indicate that platforms can substantially enhance matching outcomes and user satisfaction by implementing DH-int, thereby strengthening retention and engagement while remaining operationally efficient.

\paragraph{Organization of the paper.} Section~\ref{sec: literature} reviews the most closely related literature. Section~\ref{sec: model} introduces our model, formally defines one-lookahead policies, and presents our main heuristic. In Section~\ref{sec: analysis two period model}, we analyze the two-period version of the problem, establish performance guarantees for DH-int, and discuss other benchmarks. Section~\ref{sec: analysis T period model} extends some results to the general $T$-period setting. Section~\ref{sec: experiments} reports our empirical evaluation, and Section~\ref{sec: conclusions} concludes.

%% file: 2_literature.tex

\section{Related Literature}\label{sec: literature}

Our paper is related to several strands of the literature. First, we contribute to the literature on
assortment optimization. Most of this literature focuses on one-sided settings,
where a retailer must choose the assortment of products to show in order to maximize the
expected revenue obtained from a sequence of customers. 
We refer to \cite{kok2015} for an extensive review of the current
state of the assortment planning literature in one-sided settings.
Within this literature, the most closely related to our work is a new strand devoted to two-sided markets.
~\citet{ashlagi19} study a model where each customer simultaneously decides whether to contact a supplier from their assortment or opt out, and each supplier may then choose one contacting customer or decline to match. The platform’s goal is to design assortments that maximize expected matches. The authors show that the problem is strongly NP-hard and provide a constant-factor approximation algorithm, which is significantly improved by \citet{torrico21}.
\citet{aouad21} introduce the online version of the model in~\citep{ashlagi19}
and show that when suppliers do not accept/reject requests immediately, then a simple greedy policy achieves a 1/2-factor approximation.
All these papers study sequential two-sided matching markets, where only one side initiates and the other responds, with agents on both sides restricted to match with a single partner. In contrast, we examine markets that allow multiple selections, enabling two-directional interactions and non-sequential matches.

Within the emerging assortment optimization literature in two-sided markets, the closest paper to ours is~\citet{rios2021}, who study the same problem of selecting profiles to maximize expected matches, propose a family of one-period lookahead heuristics (Dating Heuristics), and validated it through a field experiment. Unlike their work, which focuses on estimating like probabilities and their dependence on recent matches, we establish the first performance guarantees for this problem and generalize the model to encompass different platform designs and variants of the history effect.

The second stream of literature related to our paper is on the design of matching platforms.
Starting with the seminal work of~\citet{Rochet2003}, this literature has focused on
participation, competition, and pricing, highlighting the role of cross-side externalities in different settings, including ridesharing~\citep{besbes21},
labor markets~\citep{aouad21,fonseca23}, crowdsourcing~\citep{manshadi22}, public housing~\citep{Arnosti20}, and volunteering platforms~\citep{rodilitz22}.
In the dating context,
\citet{Halaburda2018} show that two platforms can successfully coexist charging different prices by limiting the set of options offered to their users.
~\cite{cui24} analyze subscription pricing strategies and show that platforms can achieve higher revenue and welfare by varying the term of subscriptions.
\citet{Kanoria2017} show that limiting what side of the market reaches out first or hiding quality information can considerably improve the platform's outcomes.
~\cite{Immorlica22} show that platforms can achieve near-optimal welfare by limiting users' choices when they can guide their search.~\cite{celdir24} examine popularity bias in recommendation systems, showing that recommending popular profiles can increase matches and revenue if profiles remain ``within reach.'' These models assume a stylized matching market where users exit after a single match. In contrast, we allow agents to like multiple profiles and potentially form several matches over the time horizon.


{Another stream of related work examines dynamic decision making in the presence of history-dependent effects. In the context of search and optimal stopping, prior studies show that agents’ past experiences can shape subsequent behavior through frustration~\citep{rios24}, reference dependence~\citep{long20}, and other mechanisms~\citep{kagan25}. In dynamic pricing, previous work accounts for the role of past prices in shaping consumer decisions, including models with reference prices \citep{guo25}, loss aversion \citep{Hu18, Chen20}, and regret~\citep{Ozer15}. Finally, in the assortment optimization literature, recent work explores how history-dependent preferences---arising from prospect-theoretic considerations \citep{Wang18} or satiation effects \citep{He25}---affect choice behavior.  }

The last stream of the literature related to our work, motivated by applications in kidney exchange, is the stochastic matching problem in the query-commit model, also known as stochastic probing with commitment. In this problem, the matchmaker can query the edges of a general graph (e.g., to assess the compatibility between a pair donor-patient) to form a match of maximum cardinality using the accepted edges. Starting with~\cite{Chen09}, who introduce the problem (with patience constraints) and provide the first performance guarantee, most of this literature has focused on settings where the matchmaker queries only one edge at a time, see e.g.~\citep{adamczyk2011improved,bansal2012lp,gamlath2019beating,Jeloudar21}.
Nevertheless,~\cite{Chen09} also study a case closer to ours, where the planner can query a matching in each period. The authors show that a greedy algorithm that selects the edges with the highest success rate in decreasing order provides a 1/4-approximation to the optimal online algorithm when forced to commit.~\cite{Jeloudar21} study the case when there is no such commitment (i.e., the matchmaker can choose not to use an accepted edge) and show that a similar greedy algorithm achieves a 0.316-approximation guarantee. Our problem is similar in that the edges are of uncertain reward (given that like decisions are stochastic). However, we focus on selecting (or in other words, probing) subsets of profiles and each edge realization (a match) depends on the result of two outcomes. Moreover, we show that the greedy approaches used in this strand of the literature perform arbitrarily badly when adapted to our setting.


%% file: 3_problem.tex

\section{Model}\label{sec: model}

In this section, we present our model. Section~\ref{subsec: problem formulation} details its components, Section~\ref{sec:one_lookahead_policies} introduces the notion of one-lookahead policies, and Section~\ref{sec:integral_DH} describes our main method, the Integral Dating Heuristic. All proofs are deferred to Appendix~\ref{app: proofs}.


\subsection{Problem formulation}\label{subsec: problem formulation}

We now describe the dating market, the matching process, and the platform's goal and design.

\paragraph{Dating Market.}
Consider a dating platform that faces a discrete-time problem over a finite horizon of $T$ periods, where the set of periods is denoted by $[T] = \lrl{1, \ldots, T}$. Let $I$ and $J$ be two different sets of users participating on this platform, which are known at the outset and remain fixed over time, i.e., no users enter or leave the platform. This assumption captures settings involving short-term spans, e.g., platforms operating daily. To simplify the exposition, we focus on a heterosexual market; thus, users and their interactions are captured by  a bipartite graph.

At the beginning of the horizon, each user $\ell\in I\cup J$ reports their profile information (e.g., age, height, race, etc.) and their preferences regarding each of these dimensions (e.g., preferred age and height ranges, preferred races, etc.), which remain fixed over time.
The platform uses this information to compute the initial ($t=1$) set of potential partners \(\potentials^1_\ell\)---or simply \emph{potentials}---for each user $\ell$, i.e., the set of users that $\ell$ prefers over being single and for whom $\ell$ satisfies their preferences. Since our graph is bipartite, then \(\potentials^1_i \subseteq J\) for each \(i \in I\) and \(\potentials^1_j \subseteq I\) for each \(j\in J\).

\paragraph{Displays, Likes and Dislikes.}
In each period $t\in\lrc{T}$, the platform selects a subset of profiles \(S_{\ell}^{t}\subseteq \potentials^t_\ell\) to display to each user \(\ell\in I\cup J\), where $\potentials^t_\ell$ is the set of potentials for user \(\ell\) at the beginning of period \(t\). As we later discuss, the sets of potentials are updated at the end of each period to capture users' decisions and prevent them from evaluating a profile they have already seen in the past or someone who has already disliked them. To mimic our industry partner's practice, we assume that the maximum number of profiles a user can see in a given period is fixed and equal to $K_\ell$, i.e., we enforce $|S_\ell^t|\leq K_\ell$ for all $\ell\in I\cup J$ and $t\in [T]$.

For each user \(\ell\in I\cup J\) and a profile \(\ell'\in S_\ell^t\), let $\Phi^t_{\ell,\ell'}$ be the binary random variable that indicates whether $\ell$ likes $\ell'$ in period $t$, i.e., \(\Phi^t_{\ell,\ell'} = 1\) if \(\ell\) likes \(\ell'\) and \(\Phi^t_{\ell,\ell'} = 0\) otherwise.\footnote{We formalize the distribution of like probabilities below, after introducing the concept of history.} Furthermore, let \(\clikesout{\ell}^{t} = \lrl{\ell': \ell'\in S_{\ell}^{t} \text{ and } \Phi_{\ell, \ell'}^{t} = 1} \)
and \(\cdislikesout{\ell}^{t} = \{\ell': \ell' \in S_{\ell}^{t}\text{ and }\Phi_{\ell, \ell'}^{t} = 0\}\) be the random sets of profiles that user $\ell$ liked and disliked in period $t$, and let $\cumlikesout{\ell}^{<t} = \bigcup_{\tau=1}^{t-1} \clikesout{\ell}^\tau$ and $\cumdislikesout{\ell}^{<t} = \bigcup_{\tau=1}^{t-1} \cdislikesout{\ell}^\tau$ be the cumulative sets of likes and dislikes at the beginning of period $t$, respectively. 
Similarly, let \(\clikesin{\ell}^{t} = \lrl{\ell': \ell\in S_{\ell'}^{t} \text{ and } \Phi_{\ell', \ell}^{t} = 1} \)
and \(\cdislikesin{\ell}^{t} = \{\ell': \ell \in S_{\ell'}^{t}\text{ and }\Phi_{\ell', \ell}^{t} = 0\}\) be the random sets of users that liked and disliked $\ell$ in period $t$, and let $\cumlikesin{\ell}^{<t} = \bigcup_{\tau=1}^{t-1} \clikesin{\ell}^\tau$ and $\cumdislikesin{\ell}^{<t} = \bigcup_{\tau=1}^{t-1} \cdislikesin{\ell}^\tau$.\footnote{We focus on cases where users can either like or not like a profile, ruling out the skip option that is part of our partner's platform. As discussed in~\citep{rios2021}, this is without loss of generality because 5\% of profiles are skipped. Hence, throughout this paper, dislike and non-like are equivalent and, consequently, $S_\ell^t = \clikesout{\ell}^t \cup \cdislikesout{\ell}^t$. }

\paragraph{Dynamics.}
Let \(\backlog_\ell^t\) be the random backlog of user \(\ell\in I\cup J\) at the beginning of
period \(t\), i.e., the subset of users that liked \(\ell\)'s profile before period \(t\) but have not been shown to $\ell$ yet, i.e., \(\backlog_\ell^t=\potentials^t_\ell\cap\cumlikesin{\ell}^{<t}\).
The set of potentials and the backlog of user $\ell$ in period $t$ can be computed as:
{\begin{equation}\label{eq: updating formulas}
    \potentials^t_\ell = \potentials^{t-1}_\ell \setminus (S^{t-1}_\ell \cup \cdislikesin{\ell}^{t-1}), \quad
    \backlog^t_\ell = \lrp{\backlog^{t-1}_\ell \cup \mathcal{N}_{\ell,\leftarrow}^{t-1}} \setminus S^{t-1}_\ell,
\end{equation}
where $\mathcal{N}_{\ell,\leftarrow}^{t-1} = \potentials^{t-1}_\ell\cap\clikesin{\ell}^{t-1}$ is the set of likes received by $\ell$ from profiles not shown before $t-1$.}
In words, user $\ell$'s set of potentials in period $t$ can be obtained by excluding from $\ell$'s set of potentials in the previous period (i) the set of users displayed to $\ell$ in period $t-1$ and (ii) the set of users who disliked $\ell$ in $t-1$. Similarly, the backlog of user $\ell$ in period $t$ corresponds to their backlog in period $t-1$, adding the set of users who liked $\ell$ in $t-1$ and removing those profiles displayed to $\ell$ in period $t-1$. Note that, for any user $\ell \in \cIJ$ and period $t\in [T]$, (i) $\backlog_\ell^t \subseteq \potentials_\ell^t$, i.e., the backlog is a subset of the set of potentials, and (ii)
that $\potentials_\ell^t$ can only decrease as $t$ increases since no users enter or leave the platform and, consequently, the market's composition does not vary over time.

\paragraph{Matches.}
A match between users \(\ell\) and \(\ell'\) occurs if both users see and like each other. Let \(\mu_{\ell,\ell'}^t\) be the random variable whose value is one if a match between users \(\ell\) and \(\ell'\) happens in period \(t\), and zero otherwise. Then, we know that \(\mu_{\ell,\ell'}^t=1\) if and only if one of the next two events holds:
(i) users see and like each other in different periods, i.e., \(\lrl{\Phi_{\ell, \ell'}^t = 1, \; \ell'\in \backlog_\ell^t}\) or \(\lrl{\Phi_{\ell', \ell}^t = 1, \; \ell\in \backlog_{\ell'}^t}\);
or (ii) users see and like each other in period $t$, i.e., \(\ell\in S^t_{\ell'}\), \(\ell'\in S^t_{\ell}\), and \(\Phi_{\ell, \ell'}^t = \Phi_{\ell',\ell}^t = 1\). In the former case, we say that the match happens \emph{sequentially}, while in the latter, we say it happens \emph{non-sequentially}.
Notice that these two events are disjoint since users see each other at most once and, thus, we cannot simultaneously have that \(\Phi_{\ell', \ell}^t = 1\) and \(\ell'\in \backlog_\ell^t\). Finally, let $\matches_\ell^{<t}$ be the cumulative set of users that have matched with user $\ell$ at the beginning of period $t$, i.e., $\matches_\ell^{<t} = \{\ell':\; \mu_{\ell,\ell'}^\tau=1,\; \text{for all}\; \tau\in[t-1]\}$.

\paragraph{History.} As \cite{rios2021} observe, users' behavior may depend on their past experiences in the platform. To capture this, we assume that the probability that a user likes a profile may depend on their \emph{history}. Formally, let $\sigma_\ell^t$ be the history of user $\ell\in I\cup J$ at the beginning of period $t \in [T]$. Then, for each user $\ell \in \cIJ$ and $t\in [T]$, we define the realized history before period $t$ as
\begin{equation}\label{eq:history_update}
        \sigma_\ell^{t} = (L_{\ell,\rightarrow}^{<t}, D_{\ell,\rightarrow}^{<t}, B_\ell^{t}, M_\ell^{<t}),
\end{equation}
i.e., the history of each user is composed by the cumulative sets of likes and dislikes given so far, their backlog, and their cumulative set of matches. To avoid redundancy, we do not include the set of profiles observed, as this information is contained in the sets of likes and dislikes given, i.e.,  $S_\ell^\tau = L_{\ell,\rightarrow}^\tau \cup D_{\ell,\rightarrow}^\tau$. Furthermore, note that the history does not include the set of potentials nor the set of users that have previously like or disliked them---i.e., $L_{\ell,\leftarrow}^{<t}$ and $D_{\ell,\leftarrow}^{<t}$---since this information is not available to users.
Although the backlog is not observable in the free version of our partner's app, we include it in the history since it is a paid feature widely available in many dating platforms.

As we mentioned earlier, the history affects the distribution of the random variables $\Phi_{\ell,\ell'}^t$. Formally, given a user \(\ell\in I\cup J\) and a profile \(\ell'\in S_\ell^t\), 
let $\phi^t_{\ell,\ell'}(\sigma_\ell^t) = \PP(\Phi^t_{\ell,\ell'}=1\;\mid\;\sigma_\ell^t)$  be the probability that $\ell$ likes $\ell'$ in period \(t\) given their history $\sigma_\ell^t$; in the remainder of the paper, we refer to it as the \emph{like probability}.
These probabilities are independent across users and known to the platform \emph{as a function of the history}, i.e.,
the platform knows $\phi^t_{\ell,\ell'}(\sigma)$ for all $\ell'\in\potentials^1_\ell, \; t\in [T]$ and all possible realizations of the history $\sigma$, but it does not know the actual realization of the history $\sigma_\ell^t$ since the latter realizes only after period $t-1$ ends. Furthermore, we assume that like probabilities do not depend on the subset of profiles displayed, i.e., $\phi^t_{\ell,\ell'}(\cdot)$ is independent of $S^t_\ell\setminus \lrl{\ell'}$. This assumption simplifies the analysis, allowing users to potentially like several profiles without having to model preferences over sets. This is without major practical loss, as \citet{rios2021} observe that there are almost no substitution effects.
Finally, we denote by $\phi^1_{\ell,\ell'}(\sigma_\ell^1) = \phi^1_{\ell,\ell'}$ the like probabilities in the first period---absent of history---and by \(\beta^t_{\ell,\ell'}(\sigma_\ell^t,\sigma_{\ell'}^t) = \phi^t_{\ell,\ell'}(\sigma_\ell^t)\cdot \phi^t_{\ell'\ell}(\sigma_{\ell'}^t)\) the probability of a match between users $\ell$ and $\ell'$ conditional on them seeing each other in period $t$ and
their histories $(\sigma_\ell^t,\sigma_{\ell'}^t)$.

\begin{remark}
    We have not imposed any specific structure on the like probabilities as a function of the history. This is because our framework can accommodate a wide range of behavioral patterns. In Section~\ref{sec: experiments}, we provide simulation results for several parametric forms of these probabilities, including the one considered in \cite{rios2021}. The details of these specifications are provided in Section~\ref{sec: simulation setup}.
\end{remark}

\paragraph{Platform's Goal and Design.}
The platform aims to find a dynamic policy that selects a feasible subset of profiles for each user in each period to maximize the total expected number of matches throughout the entire horizon, as formalized {below} in Problem~\ref{def: general problem with multiple periods}.

A policy \(\pi \in \Pi\) for this problem prescribes a sequence of feasible subsets \(\lrl{S_\ell^{t,\pi}}_{\ell \in I\cup J}\) for \(t\in [T]\) that depends on the initial sets of potentials, the history of profiles shown, the realized like/dislike decisions, and the platform's design choices.
Let \(\Pi\) be the set of all admissible policies satisfying the design requirements (e.g., cardinality of subsets, feasibility based on potentials, etc.).
Then, the problem faced by the platform can be formulated as:
\vspace{0.25em}
\begin{problem}\label{def: general problem with multiple periods}
Given a set of admissible policies $\Pi$, the platform's problem is the following:
\[
    \max_{\pi \in \Pi} \ \E\lrc{\sum_{t=1}^T \sum_{(i,j)\in I\times J} \mu_{i,j}^{t,\pi} },
\]
where the expectation is over users' probabilistic choices and any policy randomization.
\end{problem}
\vspace{0.25em}

\begin{remark}
    Problem~\ref{def: general problem with multiple periods} can be formulated as an exponentially-sized dynamic program (DP) where the sets of potentials and history fully characterize the state of the system.
    We use the optimal value of this DP as the benchmark to compute our performance guarantees.
    One may consider a stronger benchmark that knows the realizations of all likes/dislikes in advance and computes the maximum matching in the realized graph. However, this benchmark (known as \emph{omniscient} or \emph{offline} optimum) is too strong in our setting since no dynamic policy can achieve a meaningful guarantee relative to it (even in a two-period setting). Similar observations have been previously done in the {stochastic matching} literature, see e.g., \citep{Chen09, Jeloudar21}. 
\end{remark}

\paragraph{Notation.}
Let 
\(E^1=\{\{i,j\}: i\in I,\ j\in J,\ j\in\potentials_i^1,\ i\in\potentials_j^1\}\) 
denote the set of potential undirected edges between \(I\) and \(J\). The associated directed edges are 
\(\vec E^1=\vec E_I^1\cup \vec E_J^1\), where 
\(\vec E_I^1=\{(i,j): i\in I,\ j\in\potentials_i^1\}\) and 
\(\vec E_J^1=\{(j,i): j\in J,\ i\in\potentials_j^1\}\). 
Bold symbols denote vectors or families of sets (e.g., 
\(\bm\backlog^t=\{\backlog_\ell^t\}_{\ell\in I\cup J}\)), while standard font denotes components. Calligraphic letters denote random elements and standard font their realizations (e.g., \(\mathcal B_\ell^t\) and \(B_\ell^t\)). When clear from context, dependence on the policy \(\pi\) or the history \(\sigma\) is omitted.

\subsection{One-lookahead Policies}\label{sec:one_lookahead_policies}
As discussed in~\citep{rios2021}, one-lookahead policies---i.e., policies that in each period $t$ consider $\lrl{t,t+1}$ as a time-horizon and, consequently, optimize over the current and next period only---are commonly used and perform well in practice. 
Our goal in this section is to formalize this type of policies and briefly discuss their practical limitations. Thus, in the remainder of this section, we consider a fixed period $t\in[T-1]$
and the realization of the history up to that period, {i.e., 
\begin{equation}\label{eq:def_history}
\bm{\sigma}^t = \lrl{\sigma_\ell^t}_{\ell \in \cIJ} = \lrl{\lrp{\likesout{\ell}^{<t}, \dislikesout{\ell}^{<t}, B_\ell^t, M_\ell^{<t}}}_{\ell \in \cIJ},
\end{equation}
where 
$\bm L_{\rightarrow}^{<t} = \lrl{L_{\ell,\rightarrow}^{<t}}_{\ell\in \cIJ}$ is the (realized) cumulative set of profiles liked by each user $\ell \in \cIJ$ before period $t$;
$\bm D_{\rightarrow}^{<t} = \lrl{D_{\ell,\rightarrow}^{<t}}_{\ell\in \cIJ}$ is the cumulative set of profiles disliked;
$\mathbf{B}^t = \lrl{B_\ell^t}_{\ell \in \cIJ}$ is the set of realized backlogs; and 
$\mathbf{M}^{<t} = \lrl{M_\ell^{<t}}_{\ell \in \cIJ}$ is the  set of realized cumulative matches obtained by each user at the beginning of period $t$. Furthermore, let $\mathbf{P}^t = \lrl{P_\ell^t}_{\ell \in \cIJ}$ be the set of potentials of each user $\ell$ at the beginning of period $t$, 
$\vec{E}^t = \{(\ell,\ell')\in\vec{E}^1:\; \ell'\in P^t_{\ell}\}$ the set of available arcs, and $E^t=\{\{\ell,\ell'\}\in E^1:\; \ell'\in P^t_{\ell},\; \ell\in P^t_{\ell'}\}$ the set of available edges. All the probabilities and expectations in the remainder of this section should be understood as conditional on $\mathbf{P}^t$ and $\bm{\sigma}^t$.}
{To simplify the notation, we will remove the dependency on the history and write $\phi^t_{\ell,\ell'} = \phi_{\ell,\ell'}(\sigma^t_\ell)$ and $\beta_{e}^t =\phi^t_{\ell,\ell'}(\sigma_\ell^t)\cdot\phi^t_{\ell',\ell}(\sigma_{\ell'}^t)$ for each arc $(\ell,\ell') \in \vec{E}^t$ and edge $e=\{\ell,\ell'\} \in E^t$ when clear from the context.}


The family of subsets of profiles to display in period $t$, $\{S^t_\ell\}_{\ell\in I\cup J}$, can be fully represented by two binary vectors: (i) $\mathbf{x}^t\in\{0,1\}^{\vec{E}^t}$, which captures sequential interactions and its support is the set of available arcs $\vec{E}^t$; and (ii) $\mathbf{w}^t\in\{0,1\}^{E^t}$, which captures non-sequential interactions, i.e., that both users see each other in period $t$, and its support is the set of available edges $E^t$. Then, $\ell'\in S_\ell^t$ if and only if $x_{\ell, \ell'}^{t} = 1$ or $w_{e}^{t} = 1$ with $e=\{\ell,\ell'\}$. As we discuss later, we impose the constraint $x_{\ell, \ell'}^{t}+x_{\ell', \ell}^{t}+w_{e}^{t}\leq 1$ for all $e=\{\ell,\ell'\}\in E^t$ to make the latter two events mutually exclusive so that users see each other sequentially or non-sequentially, but not both.


Our next goal is to characterize the distribution of the random events that happen during period $t$. Specifically, at the start of period $t$, the history of user $\ell$ in period $t+1$, $\sigma_{\ell}^{t+1}=\lrp{\clikesout{\ell}^{<t+1}, \cdislikesout{\ell}^{<t+1}, \cB_\ell^{t+1}, \cM_\ell^{<t+1}}$, is a random object that depends on {the history at time $t$}, the decision vectors $\lrp{\mathbf{x}^t,\; \mathbf{w}^t}$, and the like realizations in that period.
Let
\begin{equation}\label{eq: definition of set R}
    \cR=\{(\ell,\ell')\in \vec{E}^t:\; \Phi^t_{\ell,\ell'}=1\text{ and }x^t_{\ell,\ell'}=1\}\;\cup\;\{e\in E^t:\; \Phi^t_{\ell,\ell'}=\Phi^t_{\ell',\ell}=1\text{ and }w_e^t=1\}    
\end{equation}
be the random set of arcs and edges that were both displayed and liked in period $t$.  The first set captures sequential interactions, where a user $\ell$ likes another user $\ell’$. The second set captures non-sequential interactions, where two users are shown to each other in the same period and mutually like each other.
Since like decisions are independent across users, we can characterize the distribution of $\cR$ using decision vectors $\mathbf{x}^t$ and $\mathbf{w}^t$ as
\begin{equation}\label{eq:like_distribution}
    \probP_{\mathbf{z}^t}\lrp{\cR = R} = \prod_{a\in R} z_{a}^t\;\cdot\;\prod_{a\notin R} \lrp{1-z_{a}^t},
\end{equation}
where for each $a=\{\ell,\ell'\}\in R\cap E^t$ we have $z_a^t = \beta^t_{a}\cdot w^t_a$ and for each $a=(\ell,\ell')\in R\cap\vec{E}^t$ we have  $z_a^t = \plike_{\ell,\ell'}^t\cdot x_{\ell,\ell'}^t$.
In words, an element \(a=\{\ell,\ell'\}\in R\cap E^t\) is included in \(R\) if \(w_a^t=1\) and both users like each other. Conversely, it is excluded from \(R\) if either \(w_a^t=0\) or \(w_a^t=1\) but at least one of the two users did not reciprocate the like. The same logic applies to elements \(a \in \vec{E}^t\). Furthermore, recall that we impose the constraint
\(
x_{\ell, \ell'}^{t}+x_{\ell', \ell}^{t}+w_{e}^{t}\leq 1,
\)
which ensures that at most one element in \(\{(\ell,\ell'),(\ell',\ell),\{\ell,\ell'\}\}\) can belong to the realized set \(R\). This set, together with the decision variables \(\mathbf{x}^t\) and \(\mathbf{w}^t\), determines the sets of likes, dislikes, matches, and backlogs for period \(t+1\).
Specifically, we have 
\(
\likesout{\ell}^t =\{\ell': (\ell,\ell')\in R,\; x_{\ell,\ell'}^t=1\}
\), \(
\dislikesout{\ell}^t =\{\ell': (\ell,\ell')\notin R,\; x_{\ell,\ell'}^t=1\}
\), $B^{t+1}_\ell = (B^t_\ell \cup N^t_{\ell,\leftarrow})\setminus\{\ell'\in B^t_\ell:\; x_{\ell,\ell'}^t=1\}$ with $N^t_{\ell,\leftarrow}=\{\ell'\in P^t_\ell:\; (\ell',\ell)\in R\}$ as the realized version of $\mathcal{N}^t_{\ell,\leftarrow}$ defined in \eqref{eq: updating formulas}. 
Finally, the set of cumulative matches is given by
\begin{align}\label{eq:update_cumulative_matches}
    M^{<t+1}_\ell = M^{<t}_\ell\;\cup\; \{\ell':\; (\ell,\ell')\in R,\;\ell'\in B_\ell^t\}\;\cup\; \{\ell':\; (\ell',\ell)\in R,\;\ell\in B^t_{\ell'}\} \; \cup \;\{\ell':\; \{\ell,\ell'\}\in R\},
\end{align}
which accounts, in order, for the matches accumulated in previous periods, the sequential matches obtained from likes generated from the user’s backlog or from reciprocated likes coming from other users that had $\ell$ in their backlogs, and the non-sequential matches. Thus, the realized history for each user $\ell\in I\cup J$ is given by  $\sigma_\ell^{t+1} = \lrp{\likesout{\ell}^{<t+1}, \dislikesout{\ell}^{<t+1}, B_\ell^t, M_\ell^{<t+1}}$,  and we write the like probabilities in period $t+1$ as $\phi^{t+1}_{\ell,\ell'} =\phi_{\ell,\ell'}(\sigma_\ell^{t+1})$.

Given any realization $R$ resulting from the initial decisions $\lrp{\mathbf{x}^t,\mathbf{w}^t}$, the platform’s problem in period $t+1$, considering only that period and ignoring subsequent ones, is to select a subset of profiles to display to each user so as to maximize the expected number of matches. This problem can be formalized through the function $f^{t+1}(\cdot)$ defined as: 
\begin{align}\label{eq: general next problem}
        f^{t+1}(R,\mathbf{x}^t,\mathbf{w}^t)\;:=\; \max  \quad & \quad \sum_{\ell \in \cIJ}\sum_{\ell'\in B^{t}_{\ell}\cup N_{\ell,\leftarrow}^t} \phi_{\ell, \ell'}^{t+1}\cdot x^{t+1}_{\ell, \ell'} + \sum_{e\in E^{t}} \beta_{e}^{t+1}\cdot w^{t+1}_{e} \\
    \text{s.t.} \quad  & \quad \sum_{\ell'\in B^{t}_\ell\cup N_{\ell,\leftarrow}^t} x^{t+1}_{\ell,\ell'} + \sum_{e\in E^{t}:\; \ell \in e} w^{t+1}_{e} \leq \asize_\ell, \hspace{0.7cm} \forall \ell\in \cIJ, \notag\\
    & \quad w_{e}^{t+1}\leq 1-x_{\ell,\ell'}^t-x_{\ell',\ell}^t-w_{e}^t\;, \hspace{2cm} \forall e=\{\ell,\ell'\}\in E^{t}\, \notag\\
    & \quad x_{\ell,\ell'}^{t+1}\leq 1-x_{\ell,\ell'}^t\;, \hspace{4cm} \forall \ell\in I\cup J,\; \ell'\in B^{t}_\ell \notag\\
    & \quad x_{\ell,\ell'}^{t+1}\in\{0,1\},\; \hspace{4.4cm} \forall \ell\in I\cup J,\; \ell'\in B^{t}_\ell\cup N^t_{\ell,\leftarrow} \notag\\
    & \quad  w_e^{t+1}\in\{0,1\}\;, \hspace{4.4cm} \forall e\in E^{t}.\notag
\end{align}
In this formulation, the decision variables $x_{\ell,\ell'}^{t+1}$ capture sequential interactions involving only profiles $\ell'\in B^{t}_\ell\cup N^t_{\ell,\leftarrow}$---the one-lookahead horizon ends at period $t+1$, making other potential profiles irrelevant for the current optimization---, while the decision variables $\mathbf{w}^{t+1}$ capture non-sequential interactions with support in the set of available edges $E^{t}$. 
The first family of constraints captures the cardinality requirements among profiles in $B^t_\ell\cup N^t_{\ell,\leftarrow}$ and available edges in $E^t$.
The second family of constraints ensures that two users $e=\lrl{\ell, \ell'}\in E^t$ see each other in period $t+1$ (i.e., $w_{e}^{t+1} = 1$) only if none of them saw the other in period $t$ (i.e., $x_{\ell, \ell'}^t + x_{\ell', \ell}^t + w_e^t = 0$).
The third family of constraints imposes that any profile shown to $\ell$ in period $t$ (i.e., $x_{\ell,\ell'}^{t}=1$) cannot be displayed again in period $t+1$.
The last two families of constraints restrict all decision variables to be binary. Note that we do not need to impose that $x_{\ell,\ell'}^{t+1}+x_{\ell',\ell}^{t+1}+w_{e}^{t+1}\leq 1$ since the variables $\mathbf{x}^{t+1}$ are defined only for profiles taken from a user's backlog and, conditional on $w_e^t = 0$ for $e=\lrl{\ell, \ell'} \in E^t$, the events $\lrl{\ell \in B_{\ell'}^t \cup N_{\ell', \leftarrow}^t},\; \lrl{\ell' \in B_{\ell}^t \cup N_{\ell, \leftarrow}^t}$ and $\lrl{e \in E^t}$ are mutually exclusive.

 %
%
Although $f^{t+1}(\cdot)$ can be efficiently evaluated by solving a linear program (see Proposition~\ref{prop: second stage problem solved with LP} in Appendix~\ref{app: proofs}), this function is not submodular; consequently, previous results on submodularity cannot be used to derive approximation guarantees.
\begin{proposition}\label{prop: next problem is not submodular in the general case}
    The function $f^{t+1}(R,\mathbf{x}^t,\mathbf{w}^t)$ defined in~\eqref{eq: general next problem} is not submodular in \(R\). 
\end{proposition}
Despite this negative result, there are special cases in which this function is submodular, so we can extend previous results to derive provable guarantees (see Appendix~\ref{ec:design_dependent_analysis}).
We now turn to the problem in period $t$. Recall that the objective is to maximize the expected number of matches over the entire horizon $\lrl{t, t+1}$.
Formally, let $\Ma^{t+1}(\mathbf{x}^{t}, \mathbf{w}^t)$ be the maximum expected number of matches produced in period $t+1$ given the subsets of profiles displayed in period $t$ captured by $\lrp{\mathbf{x}^t, \mathbf{w}^t}$, i.e.,
\begin{equation}\label{eq: definition of expectation of function f one-lookahead}
    \Ma^{t+1}(\mathbf{x}^t, \mathbf{w}^t):= \E_{\cR \sim\mathbf{z}^t}\Big[f^{t+1}(\cR,\mathbf{x}^t, \mathbf{w}^t)\Big]  = \sum_{{R}\subseteq \vec{E}^t\cup E^t} f^{t+1}(R,\mathbf{x}^t, \mathbf{w}^t)\cdot \probP_{\mathbf{z}^t}\lrp{\cR = R}\;,
\end{equation}
where $\cR \sim \mathbf{z}^t $ represents the (random) set of likes sampled from vector $\mathbf{z}^t$. Note that the summation in~\eqref{eq: definition of expectation of function f one-lookahead} is over all possible realizations of $\cR$, i.e., all possible subsets ${R}\subseteq \vec{E}^t\cup E^t$ among potential arcs and edges.
In addition, the total expected number of matches generated in period $t$ given the decision vectors $\lrp{\mathbf{x}^t, \mathbf{w}^t}$ can be computed as:
\begin{equation}\label{eq: expected number of matches current period}
    \sum_{\ell\in I\cup J}\sum_{\ell'\in B^t_{\ell}}\phi^t_{\ell,\ell'}\cdot x_{\ell,\ell'}^t + \sum_{e\in E} \beta^t_{e}\cdot w^t_{e}.
\end{equation}
The first sum captures the expected number of sequential matches generated by showing backlog profiles, while the second sum corresponds to the non-sequential matches.
Combining~\eqref{eq: definition of expectation of function f one-lookahead} and \eqref{eq: expected number of matches current period}, we can now formalize the most general form of the \emph{one-lookahead problem}.
\vspace{1em}
\begin{problem}\label{def: problem one-lookahead}
For any given period $t\in[T-1]$,
the one-lookahead version of Problem~\ref{def: general problem with multiple periods} is:
\begin{equation}\label{eq: one-lookahead problem}
    \begin{split}
        \max  \quad & \quad \sum_{\ell\in I\cup J}\sum_{\ell'\in B^t_{\ell}}\phi^t_{\ell,\ell'}\cdot x_{\ell,\ell'}^t + \sum_{e\in E} \beta^t_{e}\cdot w^t_{e} + \Ma^{t+1}(\mathbf{x}^t, \mathbf{w}^t) \\        
        \text{s.t.} \quad & \quad \sum_{\ell'\in P^{t}_\ell} x^{t}_{\ell,\ell'} + \sum_{e\in E:\; \ell \in e} w^{t}_{e} \leq \asize_\ell, \hspace{1.7cm} \forall \ell\in \cIJ,                           \\
        & \quad x_{\ell, \ell'}^{t} + x_{\ell', \ell}^{t} + w^{t}_{e} \leq 1 \hspace{3cm}  \forall e=\{\ell,\ell'\}\in E^{t},  \\
        & \quad x_{\ell,\ell'}^{t}\in\{0,1\},\hspace{4.1cm} \forall (\ell,\ell')\in\vec{E}^t,   \\
        & \quad  w_e^{t}\in\{0,1\},\hspace{4.4cm} \forall e\in E^t,
    \end{split}
\end{equation}
\end{problem}
\vspace{1em}
%
where the first and second families of constraints enforce the cardinality constraints and that profiles are displayed either sequential or non-sequential (but not both), respectively. 

Beyond the lack of submodularity of $f^{t+1}(\cdot)$, analyzing and solving Problem~\ref{def: problem one-lookahead} exactly or approximately remains challenging for two reasons. 
First, the feasible region in~\eqref{eq: general next problem} depends on both $R$ and the decision vectors $\lrp{\mathbf{x}^t,\mathbf{w}^t}$; as a result,
the function $\Ma^{t+1}$ is highly non-linear and lacks structural properties that would simplify its analysis. 
Second, $\Ma^{t+1}(\mathbf{x}^t,\mathbf{w}^t)$ is defined as an expectation over all possible realizations of ${R}$, which themselves depend on the decisions $\lrp{\mathbf{x}^t,\mathbf{w}^t}$. Consequently, evaluating this function requires sampling techniques (e.g., Monte-Carlo simulation methods) that may be computationally intensive for realistic instances. To overcome these challenges, the next section focuses on constructing a tractable upper bound for the one-lookahead problem, forming the foundation of our Integral Dating Heuristic.

\subsection{Integral Dating Heuristic}\label{sec:integral_DH}

The Dating Heuristic (DH)~\citep{rios2021} is a one-lookahead method to solve Problem~\ref{def: general problem with multiple periods} that has been used in practice due to its simplicity and effectiveness. 
This heuristic involves solving a linear program that provides an upper bound for the one-lookahead Problem~\ref{def: problem one-lookahead}, and then rounds its solution prioritizing profiles that may produce sequential matches. However, a linearized version of Problem~\ref{def: problem one-lookahead} is only possible under some assumptions on the like probabilities in period $t+1$. As previously noted, these probabilities are history-dependent and thus depend on the decision variables $\lrp{\mathbf{x}^t,\mathbf{w}^t}$ and the realized likes and dislikes in period $t$, making the objective of Problem~\ref{def: problem one-lookahead} non-linear in the decision variables. To address this, \cite{rios2021} use $\phi^{t}_{\ell,\ell'}$ as a proxy for $\phi^{t+1}_{\ell,\ell'}$ and include a penalty function in the objective to account for the history effect.
For completeness, we formally describe the original DH (without penalty) in Algorithm~\ref{alg: DH heuristic} in Appendix.

Let us briefly discuss some technical and practical aspects of using $\phi^{t}_{\ell,\ell'}$ as a proxy for $\phi^{t+1}_{\ell,\ell'}$. 
From a technical stadpoint, for the one-lookahead problem with horizon $\{t,t+1\}$, the probability $\phi^{t}_{\ell,\ell'} = \phi_{\ell,\ell'}(\sigma_\ell^t)$ is a fixed parameter since it depends on the given history $\sigma_\ell^t$, which is part of the input of the problem. By assuming that the like probabilities remain the same in $t+1$, the non-linearity in the objective of $\Ma^{t+1}(\mathbf{x}^t,\mathbf{w}^t)$ is eliminated, allowing us to design a linear formulation that brings the decision variables $\lrp{\mathbf{x}^{t+1},\mathbf{w}^{t+1}}$ from the following period into period $t$ in their ``expected form''; for example, $x_{\ell,\ell'}^{t+1}$ would represent the probability that $\ell$ sees $\ell'$ in $t+1$ as a sequential show. 
From a practical standpoint, the appropriateness of this approximation depends on the user behavior. If like probabilities tend to decrease over time or as more information is acquired, then $\phi^{t}_{\ell,\ell'}\geq \phi^{t+1}_{\ell,\ell'}$, making the proxy an optimistic approach estimate of future outcomes. 
Similarly, if $t$ and $t+1$ correspond to short time intervals—hours or a single day—the user behavior may not change drastically, and $\phi^{t}_{\ell,\ell'}$ can serve as a reasonable approximation.
However, if the user behavior is non-monotone---e.g., users increase their interest in certain profiles or abruptly change their preferences after reaching a limit on matches---then using $\phi^{t}_{\ell,\ell'}$ as a proxy may overestimate expected results in $t+1$, potentially leading to suboptimal decisions in period $t$. 

We now describe our main proposed method: the \emph{Integral Dating Heuristic (DH-int)}. For each period $t$, let $\lrp{\mathbf{x}^t,\mathbf{w}^t}$ be the decision variables in the initial period of the horizon---equivalent to those in Problem~\eqref{eq: one-lookahead problem}---and let $\lrp{\mathbf{x}^{t+1},\mathbf{w}^{t+1}}$ be the set of variables that account for the next period decisions. Each variable $x_{\ell,\ell'}^{t+1}\in[0,1]$ represents the probability that $\ell$ sees $\ell'$ in period $t+1$ as a part of a sequential interaction, while each variable $w_{e}^{t+1}\in\{0,1\}$ with $e=\{\ell,\ell'\}$ captures that both $\ell$ and $\ell'$ see each other in period $t+1$. To overcome the non-linearity previously discussed, DH-int considers a proxy of $\Ma^{t+1}(\mathbf{x}^t,\mathbf{w}^t)$ obtained by replacing the probabilities $\phi^{t+1}_{\ell,\ell'}$ with a proxy $\hat{\phi}_{\ell,\ell'}^{t+1}$ that does not depend on the realized history in period $t$,\footnote{For instance, we can simply take $\hat{\phi}_{\ell,\ell'}^{t+1} = {\phi}_{\ell,\ell'}^{t}$ as in \citep{rios2021}. We numerically evaluate the robustness of our approach across different behaviors (including cases not satisfying our monotonicity assumption) in Section~\ref{sec: experiments}.} and linearly combining them with the decision variables $\lrp{\mathbf{x}^{t+1},\mathbf{w}^{t+1}}$.
Then, DH-int aims to maximize the expected number of matches obtained in the two-period horizon, including non-sequential and sequential matches in period $t$ and the proxy of the matches generated in period $t+1$, by solving the following mixed-integer program:
\begin{equation}\label{eq: upper bound}
  \begin{split}
    \max  \quad & \sum_{\ell\in I\cup J}\sum_{\ell'\in B^t_{\ell}}\phi^t_{\ell,\ell'}\cdot x_{\ell,\ell'}^t +\sum_{e\in E^t} \beta^t_{e}\cdot w^t_{e}+ \sum_{\ell\in I\cup J}\sum_{\ell'\in P^t_{\ell}}\hat{\phi}^{t+1}_{\ell,\ell'}\cdot x_{\ell,\ell'}^{t+1}+ \sum_{e\in E^t} \hat{\beta}^{t+1}_{e}\cdot w^{t+1}_{e} \\
    \text{s.t.} \quad & \quad \sum_{\ell'\in P^{t}_\ell} x^{\tau}_{\ell,\ell'} + \sum_{e\in E:\; \ell \in e} w^{\tau}_{e} \leq \asize_\ell, \hspace{1.7cm} \forall \ell\in \cIJ,\; \tau\in\{t,t+1\}\\
    & \quad x_{\ell, \ell'}^{t} + x_{\ell', \ell}^{t} + w^{t}_{e} + w_e^{t+1} \leq 1 \hspace{1.8cm} \forall e=\{\ell,\ell'\}\in E^{t}, \\
    & \quad x^{t+1}_{\ell,\ell'} \leq \phi^{t}_{\ell',\ell}\cdot x^{t}_{\ell',\ell}, \hspace{3.4cm} \forall \ell\in \cIJ,\; \ell' \in P^t_{\ell}\setminus B_\ell^t\\
    & \quad x^{t}_{\ell,\ell'} + x^{t+1}_{\ell,\ell'}\leq 1, \hspace{3.7cm} \forall \ell\in \cIJ,\; \ell' \in  B_\ell^t\\
    & \quad \mathbf{x}^t\in\{0,1\}^{\vec{E}^t},\; \mathbf{x}^{t+1}\in[0,1]^{\vec{E}^t},
    \;  \mathbf{w}^t, \; \mathbf{w}^{t+1}\in\{0,1\}^{E^t}
  \end{split}
\end{equation}

Similar to Problem~\eqref{eq: one-lookahead problem}, the first and second families of constraints impose the cardinality constraints for each period and that profiles are displayed either sequential or non-sequential in period $t$. Furthermore, the latter constraints guarantee that non-sequential shows in period $t+1$ only hold if none of the users saw the other in period $t$.
The third family of constraints limits sequential shows of non-backlog profiles in period $t+1$, ensuring that  $\ell$ sees $\ell'\in P^t_{\ell}\setminus B^t_\ell$ only if  $\ell'$ saw $\ell$ in period $t$ (i.e., $x_{\ell',\ell}^t = 1$) and liked them. Finally, the fourth family of constraints forces that backlog profiles are shown at most once in the two-period horizon.
Note that the feasible region in~\eqref{eq: upper bound} is different than that in the original implementation of DH; however, in~\ref{app: dating heuristic}, we show that these regions are equivalent. More importantly, the key difference between DH-int and DH lies in the integrality of the first period decisions, which serves multiple purposes:
(i) it allows us to derive our guarantees in a simpler way as it does not need the rounding step used in \citep{rios2021}, which may lead to inefficiencies; 
(ii) its connection to the one-lookahead Problem~\ref{def: problem one-lookahead} is direct; and
(iii) the search for integral solutions does not generate an extra computational burden.

Based on the optimal solution \( \lrp{\mathbf{x}^t,\mathbf{w}^t, \mathbf{x}^{t+1}, \mathbf{w}^{t+1}}\) of Problem~\eqref{eq: upper bound}, DH-int 
constructs the subsets to display in period $t$ by setting $S_\ell^t=\{(\ell,\ell')\in\vec{E}^t:\; x^t_{\ell, \ell'} = 1\}\cup\{e\in E^t:\; w_e^t=1\}$ for each $\ell \in \cIJ$. Then, the algorithm updates the potentials, backlogs and history  according to~\eqref{eq: updating formulas} and \eqref{eq:def_history} based on the realized like/dislike decisions. 
We formalize this procedure in Algorithm~\ref{alg: integral DH heuristic}.
{\small
\begin{algorithm}[htp!]
      \caption{Integral Dating Heuristic (DH-int)}\label{alg: integral DH heuristic}
      \begin{algorithmic}[1]
      \Require An instance of Problem~\ref{def: general problem with multiple periods}.
      \Ensure A feasible subset of profiles to display to each user in each period.
      \For{$t\in[T]$}
      \State Solve Problem~(\ref{eq: upper bound}) and let \(\lrp{\mathbf{x}^t, \mathbf{w}^t, \mathbf{x}^{t+1},\mathbf{w}^{t+1}}\) be the optimal solution.
      \State \multiline{For each user \(\ell \in \cIJ\), display subset $S_\ell^t=\{(\ell,\ell')\in\vec{E}^t:\; x^t_{\ell, \ell'} = 1\}\cup\{e\in E^t:\; w_e^t=1\}$.}
      \State Observe like/dislike decisions. Update potentials, backlogs and history following~\eqref{eq: updating formulas} and \eqref{eq:def_history}.  
        \EndFor
      \end{algorithmic}
\end{algorithm} 
}

%% file: 4_analysis_two_periods.tex
\section{Analysis for the Two-period Model}\label{sec: analysis two period model}
In the next two sections, we establish that DH-int achieves strong performance guarantees across different platform designs. To simplify the analysis and extract insights that will be useful in the general case, we begin with a two-period version of the model. Section~\ref{sec: two-period model} formalizes the setting. Section~\ref{sec: guarantee for DH} presents our analysis of DH-int. Finally, Section~\ref{sec: design_dependent_analysis} discusses guarantees for polynomial-time algorithms that do not rely on mixed-integer programs.

\subsection{Setting}\label{sec: two-period model}
We start by describing the two-period model, specifying some key assumptions, and discussing how the model described in the previous section adapts to this setting.

\paragraph{Assumptions.} Due to the lack of submodularity stated in Proposition~\ref{prop: next problem is not submodular in the general case} and the discussion on history-dependent like probabilities in Section~\ref{sec:integral_DH}, in the remainder of this section, we make some assumptions to ease the analysis for the case when $T=2$. First, we assume {for simplicity} that there is no initial backlog, i.e.,  $B^1_\ell=\emptyset$.\footnote{{Our approach can be easily extended to consider non-empty initial backlogs.}} This makes the feasible region in \eqref{eq: general next problem} not dependent on $B^1_\ell$ and, consequently, on the first-stage decision variables $x_{\ell,\ell'}^1$. Moreover, this implies that the backlog sets in $t=2$ only include profiles in $N^1_{\ell,\leftarrow}$. 
Second, we allow like probabilities $\phi^2_{\ell,\ell'}$ to be time-dependent but we assume that they do not depend on the realized history from period $t=1$. With this assumption, we avoid the non-linear structure in the objective function of~\eqref{eq: general next problem}. Finally, we assume that non-sequential matches in period $t=2$ are forbidden, which means that variables $w_e^2$ are set to 0. This implies that the feasible region in~\eqref{eq: general next problem} does not depend on the first-stage variables $w_e^1$. These assumptions simplify the analysis and allow us to leverage submodular optimization techniques and other technical tools, such as the correlation gap~\citep{agrawal2010correlation}, to solve the problem.
Moreover, this is without practical loss because most matches materialize sequentially.
In~\ref{subsec: simultaneous matches in second period}, we discuss the general case with non-sequential matches in the second stage, discuss its challenges, and provide a guarantee for when like probabilities are small.

\paragraph{Formulation.} 
Given the assumptions above, the random set $\cR$ that results from the first period decision variables $\mathbf{x}^1,\mathbf{w}^1$ and the like decisions becomes $\cR = \{(\ell,\ell')\in \vec{E}^1:\; \Phi_{\ell,\ell'}^1=1,\; x_{\ell,\ell'}^1=1\}$. 
Moreover, its distribution reduces to
\begin{equation}\label{eq:backlog_distribution two periods}
        \probP_{\mathbf{x}^1}\left(\cR = R \right) = 
        \prod_{(\ell,\ell')\in R} \plike_{\ell,\ell'}^1\cdot x_{\ell,\ell'}^1 \cdot 
        \prod_{(\ell,\ell')\notin R} \lrp{1-\plike_{\ell',\ell}^1 \cdot x_{\ell',\ell}^1}. 
\end{equation}
where $R$ is the realized subset of arcs in $\vec{E}^1$. {Note that this set does not include the set of edges capturing pairs of users that see and like each other in the first period as captured in~\eqref{eq: definition of set R}. This is because non-sequential interactions are excluded in the second period, only profiles taken from the backlog can be displayed in the second period and, consequently, we only need to keep track of users' backlogs to obtain the maximum expected number of matches that can be produced in the latter period. As a result,} Problem~\eqref{eq: general next problem} can be re-written as:
\begin{equation}\label{eq:second_stage_problem}
        f(R)\;:=\; \max_{\mathbf{x}^2\in\{0,1\}^{\vec{E}^1}}\Bigg\{ \sum_{(\ell,\ell')\in\vec{E}^1} \phi_{\ell, \ell'}^{2}\cdot x^{2}_{\ell, \ell'}: \ \sum_{\ell'\in  P^1_\ell} x^{2}_{\ell,\ell'} \leq \asize_\ell,\;\forall \ell\in \cIJ, \;
        x_{\ell,\ell'}^{2}\leq \indicator_{\{(\ell',\ell)\in R\}},\ \forall (\ell,\ell')\in\vec{E}^1 \Bigg\}
\end{equation}
where we remove the dependency of $f$ on $\lrp{\mathbf{x}^1, \mathbf{w}^1}$ to simplify the exposition.
In this problem, the variables $x^2_{\ell,\ell'}$ are set to zero if $(\ell',\ell)\notin R$. Namely, if $\ell'$ did not like $\ell$ or did not see $\ell$ in the first period, then $\ell$ will not see profile $\ell'$ in the second stage. In other words, the matches  in the second period are generated from sequential shows responding to a sequential like in the first period. More importantly, the feasible region in~\eqref{eq:second_stage_problem} is clearly a partition matroid---the  edges incident to each $\ell\in I\cup J$ form a partition over the set of available edges, which is determined by the backlogs in the second family of constraints---and function $f$ is submodular, as we formalize in the next lemma.
\begin{lemma}\label{lemma: f monotone submodular} The function $f$ defined in~\eqref{eq:second_stage_problem} is monotone and submodular.
\end{lemma}
{Problem~\eqref{eq:second_stage_problem} can be solved in polynomial time since we are maximizing a linear function over a single matroid.} 
Nevertheless, the two-period problem is significantly more challenging since the set $R$ depends on the first-period decisions. Formally, the function in \eqref{eq: definition of expectation of function f one-lookahead} with $t=1$ corresponds to  
    \begin{equation}\label{eq: definition of expectation of function f with sequential matches}
    \Ma^2(\mathbf{x}^1, \mathbf{w}^1):= \E_{\mathbf{\cR} \sim \mathbf{\phi}^1 \cdot\mathbf{x}^1}\lrc{f(\cR)}  = \sum_{R\subseteq\vec{E}^1} f(R)\cdot \probP_{\mathbf{x}^1}\lrp{\cR=R}\;,
    \end{equation}
    where $\phi^1 \cdot \mathbf{x}^1 $ denotes the vector with components $\phi^1_{\ell,\ell'} \cdot x^1_{\ell,\ell'} $ and $\cR \sim \mathbf{\phi}^1 \cdot \mathbf{x}^1 $ represents the (random) set of arcs sampled from vector $\phi^1 \cdot \mathbf{x}^1$. 
    Note that $\Ma^2(\mathbf{x}^1, \mathbf{w}^1)$ is also monotone and submodular in $\mathbf{x}^1$ (see Lemma~\ref{lemma:general_function_submodular_monotone} in Appendix~\ref{os:proofs_2periods}).
    Finally, since the initial backlog is empty, the total expected number of matches in the first period are only generated from non-sequential shows, i.e., $\sum_{e\in E^1} \beta^1_{e}\cdot w^1_{e}$.
    \vspace{0.25em}
    \begin{problem}\label{def:two_period_problem}
    Given a set of admissible policies $\Pi$, the two-period version of Problem~\ref{def: general problem with multiple periods} is:
    \begin{equation}\label{eq: objective two-period problem}
    \begin{split}
      \max  \quad       & \sum_{e\in E^1} \beta^1_{e}\cdot w^1_{e} + \Ma^{2}(\mathbf{x}^1, \mathbf{w}^1) \\
        \text{s.t.} \quad & \quad \sum_{\ell'\in P^{1}_\ell} x^{1}_{\ell,\ell'} + \sum_{e\in E:\; \ell \in e} w^{1}_{e} \leq \asize_\ell, \hspace{1.7cm} \forall \ell\in \cIJ,                           \\
        & \quad x_{\ell, \ell'}^{1} + x_{\ell', \ell}^{1} + w^{1}_{e} \leq 1 \hspace{3cm}  \forall e=\{\ell,\ell'\}\in E^1  \\
        & \quad \mathbf{x}^1\in\{0,1\}^{\vec{E}^1},\quad  \mathbf{w}^{1}\in\{0,1\}^{E^1}.
    \end{split}
    \end{equation}
    \end{problem}
    \vspace{0.25em}
In the remainder of this section, we denote the optimal value of Problem~\ref{def:two_period_problem} as $\OPT^\Pi$ where $\Pi$ is the set of admissible policies which captures an specific platform design. We remove the dependency on $\Pi$ when clear from context. Also, we say that a policy guarantees an $\alpha$-approximation with $\alpha\in[0,1]$ if its objective value is at least an $\alpha$ fraction of $\OPT^\Pi$.

Several approaches that perform well in the matching and assortment optimization literature fail to provide meaningful guarantees for Problem~\ref{def:two_period_problem}.  
In Appendix~\ref{app:suboptimal_methods}, we show that the guarantees of two classic policies  in the literature (a local greedy method and finding a perfect matching in each step) go to zero as the size of market grows. Instead, in the next section, we show that our method (DH-int) has a stronger provable guarantee that holds across all platform designs that can be captured by Problem~\ref{def:two_period_problem}. 

    

\subsection{Analysis of DH-int}\label{sec: guarantee for DH}
For our two-period model, DH-int solves the following version of \eqref{eq: upper bound}:
\begin{equation}\label{eq:upper_bound_twoperiods}
      \begin{split}
        \max  \quad & \sum_{e\in E^1} \beta^1_{e}\cdot w^1_{e}+ \sum_{\ell\in I\cup J}\sum_{\ell'\in P^1_{\ell}}{\phi}^{2}_{\ell,\ell'}\cdot x_{\ell,\ell'}^{2} \\
        \text{s.t.} \quad & \quad \sum_{\ell'\in P^{1}_\ell} x^{1}_{\ell,\ell'} + \sum_{e\in E^1:\; \ell \in e} w^{1}_{e} \leq \asize_\ell, \hspace{1.7cm} \forall \ell\in \cIJ,\\
        & \quad \sum_{\ell'\in P^{1}_\ell} x^{2}_{\ell,\ell'} \leq \asize_\ell, \hspace{3.7cm} \forall \ell\in \cIJ,\\
        & \quad x_{\ell, \ell'}^{1} + x_{\ell', \ell}^{1} + w^{1}_{e}  \leq 1 \hspace{2.8cm} \forall e=\{\ell,\ell'\}\in E^{1}, \\
        & \quad x^{2}_{\ell,\ell'} \leq \phi^{t}_{\ell',\ell}\cdot x^{1}_{\ell',\ell}, \hspace{3.4cm} \forall \ell\in \cIJ,\; \ell' \in P^1_{\ell}\\
        & \quad \mathbf{x}^1\in\{0,1\}^{\vec{E}^1},\; \mathbf{x}^{2}\in[0,1]^{\vec{E}^1},
        \;  \mathbf{w}^1, \in\{0,1\}^{E^1}
      \end{split}
    \end{equation}
Due to our assumptions on history independence, note that the proxy probabilities in period $t=2$ are $\phi^2_{\ell,\ell'}$. Moreover, variables $w^2_e$ are not present since we forbid non-sequential shows in $t=2$. Finally, we assume for simplicity that there is no initial backlog, i.e., $B^1_\ell=\emptyset$ for all $\ell\in I\cup J$.
Problem~\eqref{eq:upper_bound_twoperiods} is a mixed-integer programming relaxation of the two-period Problem~\ref{def:two_period_problem}; consequently, it is unclear whether it can be solved in polynomial time.\footnote{Although the objective function is linear, the problem is challenging because the feasible region is defined by the intersection of several partition matroids and involves both continuous and binary decision variables.} Nevertheless, as we show in Theorem~\ref{thm: guarantee for DH}, the implementation of DH-int (in Algorithm~\ref{alg:DH-int-2periods}) has a strong approximation guarantee. 
{\small
\begin{algorithm}
      \caption{Integral Dating Heuristic (DH-int) for $T=2$}\label{alg:DH-int-2periods}
      \begin{algorithmic}[1]
      \Require An instance of Problem~\ref{def:two_period_problem}.
      \Ensure A feasible subset of profiles to display to each user in each period.
      \State Solve Problem~(\ref{eq:upper_bound_twoperiods}) and let \(\lrp{\mathbf{x}^1, \mathbf{w}^1, \mathbf{x}^{2}}\) be an optimal solution.
      \State \multiline{For each user \(\ell \in \cIJ\), display subset $S_\ell^1=\{(\ell,\ell')\in\vec{E}^1:\; x^1_{\ell, \ell'} = 1\}\cup\{e\in E^1:\; w_e^1=1\}$.}
      \State Observe the realized set of arcs $R$. Solve~\eqref{eq:second_stage_problem} and let $\mathbf{x}^2$ be an optimal solution. 
      \State \multiline{For each user \(\ell \in \cIJ\), display subset $S_\ell^2=\{(\ell,\ell')\in\vec{E}^1:\; x^2_{\ell, \ell'} = 1\}.$}
      \end{algorithmic}
\end{algorithm}
}
\begin{theorem}\label{thm: guarantee for DH}
    Algorithm~\ref{alg:DH-int-2periods} achieves a $(1-1/e)$-approximation for Problem \ref{def:two_period_problem}.
\end{theorem}
	The proof of Theorem~\ref{thm: guarantee for DH} (in Appendix~\ref{os:proof_maintheorem_2periods}) relies on a novel connection between formulation~\eqref{eq:upper_bound_twoperiods} and a different relaxation that aims to find the distribution of backlogs (under any possible distribution) that maximizes the expected total number of matches while satisfying constraints on the marginal probabilities. To see this, consider an optimal solution $(\mathbf{x}^{1,\star}, \mathbf{w}^{1,\star}, \mathbf{x}^{2,\star})$ of \eqref{eq:upper_bound_twoperiods}. {Note that the fractional point $\mathbf{x}^{2,\star}$ is a feasible solution for the second stage problem:}
	\begin{align}\label{eq:P1}
		F(\mathbf{x}^{1,\star}, \mathbf{w}^{1,\star}) := \max & \quad \sum_{\ell \in \cIJ} \sum_{\ell'\in \potentials^1_\ell}\phi_{\ell, \ell'}^2 \cdot y_{\ell, \ell'}\\
		s.t.&\quad y_{\ell, \ell'}\leq \phi^1_{\ell', \ell}\cdot x^{1,\star}_{\ell', \ell}, \hspace{2.25cm}  \forall \ell \in \cIJ, \;\ell' \in \potentials^1_\ell \notag\\
		&\quad \sum_{\ell' \in \potentials^1_\ell} y_{\ell, \ell'}\leq K_\ell, \hspace{2.5cm}\forall \ell \in \cIJ\notag\\
		&\quad y_{\ell, \ell'}\geq 0,  \hspace{3.7cm}\forall \ell \in \cIJ, \ell' \in \potentials^1_\ell.\notag
	\end{align}		
    {Therefore, the optimal value of \eqref{eq:upper_bound_twoperiods} is upper bounded by $\sum_{e\in E^1} \beta^1_{e}\cdot w^{1,\star}_{e}+F(\mathbf{x}^{1,\star}, \mathbf{w}^{1,\star})$.}
	To provide a bound for $F(\mathbf{x}^{1,\star}, \mathbf{w}^{1,\star})$, note that a second upper bound for the expected number of matches in the second period can be obtained for any vectors of first-period decisions $\lrp{\mathbf{x}^1, \mathbf{w}^1}$ by solving:
	\begin{align}\label{eq:P2}
		G(\mathbf{x}^1, \mathbf{w}^1) := \max & \quad \sum_{\ell \in \cIJ} \sum_{B \subseteq \potentials_\ell^1} f_\ell(B)\cdot\lambda_{\ell,B}\\
		s.t.&\quad \sum_{B\subseteq \potentials_\ell^1}\lambda_{\ell,B} = 1 \hspace{3.75cm} \forall \ell \in \cIJ\notag\\
		&\quad \sum_{B\subseteq \potentials_\ell^1: \ell'\in B_\ell} \lambda_{\ell,B} = \phi^1_{\ell, \ell'}\cdot x^1_{\ell, \ell'} , \hspace{1.4cm} \forall \ell \in \cIJ, \ell'\in \potentials_\ell^1 \notag\\
		&\quad \lambda_{\ell,B}\geq 0, \hspace{4.5cm} \forall \ell \in \cIJ, \; B\subseteq \potentials_\ell^1.\notag
	\end{align}
	In this formulation, which we refer to as the \emph{distribution problem}, the decision variables $\lambda_{\ell, B}$ represent the probability that the set $B\subseteq \potentials_\ell^1$ is the backlog of user $\ell$ in the second period. In other words, if $R$ is the set of realized arcs, then for each user $\ell$, the backlog would be $B=\{\ell': (\ell',\ell)\in R\}$; we omit the subindex $\ell$ in $B$ because it is just a generic set. On the other hand,
    function $f_\ell(B)$ returns the maximum expected number of matches that user $\ell$ can achieve given a backlog $B$ (as defined in~\eqref{eq:second_stage_problem}, but for each user $\ell$). The first and last families of constraints ensure that the variables $\lrl{\lambda_{\ell, B}}_{B\subseteq \potentials_\ell^1}$ correctly define a probability distribution, while the second family of constraints guarantees that the marginal probabilities are consistent with the first-period decisions. 
	
	Based on these two formulations, the next step in the proof shows that $G(\mathbf{x}^{1,\star}, \mathbf{w}^{1,\star})$ is at least $F(\mathbf{x}^{1,\star}, \mathbf{w}^{1,\star})$ and that a feasible solution for the latter can be constructed as 
	\(\lambda_{\ell,B} = \PP_{\mathbf{x}_\ell^{1,\star}}(B)\) (similarly to the distribution defined in~\eqref{eq:backlog_distribution two periods} but for a given $\ell$). Noticing that the expected number of matches in the second period produced by this feasible solution, \(\lambda_{\ell,B} = \PP_{\mathbf{x}_\ell^{1,\star}}(B)\) for all $\ell\in I\cup J$, coincides with $\Ma^2(\mathbf{x}^{1,\star}, \mathbf{w}^{1,\star})$, we combine strong-duality with the correlation gap in~\citep{agrawal2010correlation} to show that 
	\[\Ma^2(\mathbf{x}^{1,\star}, \mathbf{w}^{1,\star}) \geq (1-1/e)\cdot G\lrp{\mathbf{x}^{1,\star}, \mathbf{w}^{1,\star}} \geq (1-1/e)\cdot F\lrp{\mathbf{x}^{1,\star}, \mathbf{w}^{1,\star}}.\]
	Finally, recalling that $\sum_{e\in E^1}\beta_e^1\cdot w^{1,\star}_e + F(\mathbf{x}^{1,\star}, \mathbf{w}^{1,\star}) \geq \OPT$ (since~\eqref{eq:upper_bound_twoperiods} provides an upper bound of $\OPT$ and the left-hand side coincides with its optimal solution), we conclude that 
	\[
	\sum_{e\in E^1}\beta_e^1\cdot w^{1,\star}_e + \Ma^2(\mathbf{x}^{1,\star}, \mathbf{w}^{1,\star}) \geq (1-1/e)\cdot \OPT.
	\]

\subsection{Analysis of Polynomial-Time Approximation Algorithms}\label{sec: design_dependent_analysis}
Despite DH-int's strong approximation guarantee for the two-period model, one disadvantage is that it is not a polynomial-time algorithm since it solves a mixed-integer program, creating scalability threats for large markets. As a tractable alternative, in the Appendix~\ref{ec:design_dependent_analysis}, we present polynomial-time approximation algorithms based on submodular techniques that provide weaker constant-factor approximation guarantees that depend on (i) the timing of matches (only sequential or adding non-sequential), and (ii) the sequence of interactions (one or two-directional). Specifically,
\begin{enumerate}
    \item[(i)] \emph{Timing of Matches.} A policy $\pi$ restricts to sequential matches if no pair of users see each other in the same period. Formally, for any $(i,j)\in I\times J$ and period $t$, $i\in S_j^{t, \pi}$ implies $j\notin S_i^{t, \pi}$, and vice versa. Under such a policy, matches can only occur if users see and like each other in different periods. In contrast, $\pi$ allows non-sequential matches if no such restriction applies.    
    \item[(ii)] \emph{Sequence of Interactions.} A policy $\pi$ is one-directional if only one side can initiate the path towards a matches. Suppose $I$ is the initiating side. Then, in any period $t$, $\batch_i^{t, \pi} \subseteq \potentials_i^t$ for each $i\in I$, while $\batch_j^{t, \pi} \subseteq \backlog_j^t \cup \{i\in \potentials_j^t : j\in \batch_i^{t,\pi}\}$ for $j\in J$.\footnote{By the updating formulas in~(\ref{eq: updating formulas}), $i\in \potentials_j^t$ and $j\in S_i^\tau$ for some $\tau < t$ imply $i\in \backlog_j^t$.}  
    In words, users in $I$ can see any profile in their set of potentials, while users in $J$ can only see profiles in their backlog or those currently viewing them. A policy is two-directional if any side can initiate, i.e., $\batch_\ell^{t, \pi} \subseteq \potentials_\ell^t$ for all $\ell \in I \cup J$.  
\end{enumerate}
While most dating platforms employ a two-directional design, some adopt a one-directional one by restricting which side can initiate interactions. For example, Coffee Meets Bagel’s ``Ladies Choice'' allows only men to initiate, so women see profiles of men who have previously liked them. Some platforms have also introduced non-sequential interactions: Tinder’s ``Hot Takes’’ pairs users for limited-time chats before mutual likes, while Filteroff and The League offer short video-based speed-dating sessions that require participants to be present at the same time.

A key feature of our model is that it can easily accommodate these different platform designs by simply modifying the feasible region in Problem~\ref{def:two_period_problem}. More importantly, by using submodular optimization techniques, we can provide constant-factor approximation guarantees for each of these design choices, as formalized in Theorem~\ref{thm: submodular guarantees}. We defer the analysis and proofs to the Appendix~\ref{ec:design_dependent_analysis}.
\begin{theorem}\label{thm: submodular guarantees}
Algorithm~\ref{alg:greedy} (see Appendix~\ref{ec:design_dependent_analysis}) achieves design-dependent constant factor guarantees for Problem \ref{def:two_period_problem}.
\end{theorem}

%% file: 5_analysis_multiple_periods.tex
\section{Analysis for the Multiple Periods Model}\label{sec: analysis T period model}

    In this section, we provide guarantees for the general $T$-period version of the model. To simplify the analysis, we assume that like probabilities are time-invariant throughout this section, i.e., $\phi_{\ell, \ell'}^t = \phi_{\ell, \ell'}$ and $\beta_{e}^t = \beta_{e}$ for all $t\in [T], \; (\ell,\ell')\in \vec{E}^1$ and $e\in E^1$. 
    Our main result is that, provided that we exclude non-sequential matches in the last period, a {semi}-adaptive variant of DH-int---called \emph{DHT}---achieves a constant-factor approximation guarantee for any time horizon when we restrict to either \emph{semi-adaptive} or \emph{adaptive} policies, as defined next.
    \vspace{0.25em}
    \begin{definition}
        We say that a policy is \emph{semi-adaptive} if
        (i) it non-adaptively selects profiles for the initiating side, i.e., when selecting profiles to start a sequential or a non-sequential interaction; and (ii) it adaptively selects profiles for the responding side, i.e., when selecting profiles from the backlog.
        In contrast, we say that a policy is \emph{adaptive} if it selects the profiles to display in both sides (initiating and responding) using information from previous stages.
    \end{definition}
    \vspace{0.25em}
    Consider the following extension of the upper bound in~\eqref{eq:upper_bound_twoperiods}:
    \begin{equation}\label{eq: problem DH multiple-periods}
        \begin{split}
            \max \quad & \; \sum_{e\in E}  w_{e} \cdot \beta_{e} + \sum_{\ell\in I\cup J}\sum_{\ell' \in \potentials^1_\ell} y_{\ell, \ell'} \cdot \plike_{\ell, \ell'} \\
            \text{s.t.} \quad 
        & \; y_{\ell, \ell'} \leq x_{\ell', \ell} \cdot \plike_{\ell', \ell}, \hspace{4.6cm} \forall \ell \in I\cup J, \; \ell' \in \potentials^1_\ell, \\
        & \; \ x_{\ell, \ell'} + x_{\ell', \ell} + w_{e} \leq 1, \hspace{3.85cm} \forall e\in E^1, \\ 
        & \; \sum_{\ell' \in \potentials^1_\ell} x_{\ell, \ell'} + \sum_{e\in E^1:\; \ell \in e} w_{e} + \sum_{\ell' \in \potentials^1_\ell} y_{\ell, \ell'} \leq K_\ell\cdot T, \quad \forall \ell \in I\cup J, \\      
        & \; \mathbf{x} \in \lrl{0,1}^{\vec{E}^1},\; \mathbf{w} \in \lrl{0,1}^{E^1},\; \mathbf{y}\in [0,1]^{\vec{E}^1}.
        \end{split}
    \end{equation}
    The decision variables $x_{\ell, \ell'}$ indicate whether user $\ell$ sees $\ell'$'s profile throughout the $T$-period horizon, initiating their interaction. 
    The decision variables $y_{\ell, \ell'}$ represent the probability that $\ell$ sees $\ell'$ as the respondent/follower of the interaction, which happens only conditional on $\ell'$ liking $\ell$ in the first place. Finally, the decision variables $w_{e}$ with $e=\{\ell,\ell'\}$ capture whether $\ell$ and $\ell'$ both see each other in the same period. Note that, when $T=2$, the third constraint in~\eqref{eq: problem DH multiple-periods} is equivalent to the aggregation of the third and fourth constraints in \eqref{eq:upper_bound_twoperiods}.
    As we show in Lemma~\ref{lemma: UB semi-adaptive}, \eqref{eq: problem DH multiple-periods} provides an upper bound for any semi-adaptive policy.

    \begin{lemma}\label{lemma: UB semi-adaptive}
        Problem~\eqref{eq: problem DH multiple-periods} is an upper bound for~\ref{def: general problem with multiple periods} when $\Pi$ is restricted to semi-adaptive policies.
    \end{lemma}


    Based on this result, we can design a $T$-period version of DH-int, formalized in Algorithm~\ref{alg: multiple-period DH}. This algorithm starts by solving Problem~\eqref{eq: problem DH multiple-periods}, whose optimal solution $\lrp{\mathbf{x}^\star, \mathbf{w}^\star, \mathbf{y}^\star}$ provides an upper bound for the maximum expected number of matches that any semi-adaptive policy can obtain in $T$ periods (by Lemma~\ref{lemma: UB semi-adaptive}). 
    {The algorithm then fixes this solution and, in each period $t$, iterates over users to determine which profiles to display based solely on $\mathbf{x}^\star$ and $\mathbf{w}^\star$. Importantly, no re-optimization is performed in response to realized likes and dislikes from earlier periods, in contrast to DH and DH-int, which adapt their decisions as new information becomes available.}
    The algorithm {implements this solution by} first exhausting the profiles initiating an interaction, i.e., $x_{\ell,\ell'}=1$ (Steps 6-8). If there is still space left in $S^t_\ell$, the algorithm displays profiles that are part of non-sequential interactions (i.e., $w^\star_{e} = 1$) as long as the cardinality constraint is not binding for both users (Steps 9-12). If the cardinality constraint is still non-binding, the algorithm computes the optimal subset of profiles to display to $\ell$ based on the profiles in their backlog (Step~\ref{alg-step: solving f}). Finally, the algorithm updates the potentials and backlogs based on the observed like/dislike decisions (Step~\ref{alg-step: updating potentials and backlogs}).
    {\small
    \begin{algorithm}[h]
      \caption{$T$-periods DH-int (DHT)}\label{alg: multiple-period DH}
      \begin{algorithmic}[1]
      \Require An instance of Problem~\ref{def: general problem with multiple periods}.
      \Ensure A subset of profiles to display in each period.
      \State Solve \eqref{eq: problem DH multiple-periods} and let $\lrp{\mathbf{x}^\star, \mathbf{w}^\star, \mathbf{y}^\star}$ be the optimal solution.
      \State Define $X_\ell = \lrl{\ell' \in \potentials_\ell^1: x^\star_{\ell, \ell'} = 1}$ for each $\ell \in \cIJ$ and $W_\ell = \lrl{e\in E^1: \ell\in e, \ w^\star_{e} =1}$.
      \For{$t\in[T]$}
        \State Initialize $S_\ell^t = \emptyset$ for each $\ell \in \cIJ$.
        
        \For{$\ell\in\cIJ$}
            \While{$|S_\ell^t|<K_\ell$ and $X_\ell\neq\emptyset$}
                \State Let $\ell'$ be $\phi_{\ell',\ell}\in\argmax\{\phi_{a,\ell}: \; a\in X_\ell\}$ breaking ties arbitrarily. \label{alg-step: selecting profiles to initiate interactions with highest probs}
                \State $S_\ell^t = S_\ell^t \cup \lrl{\ell'}$, $X_\ell = X_\ell\setminus\{\ell'\}$. 
            \EndWhile
            \While{$|S_\ell^t|<K_\ell$ and $W_\ell\neq\emptyset$}
                \State Choose any $\ell'\in W_\ell$.
                \If{$|S_{\ell'}^t| < \asize_{\ell'}$}
                    \State $S_\ell^t = S_\ell^t \cup \lrl{\ell'}, \; S_{\ell'}^t = S_{\ell'}^t \cup \lrl{\ell}$, $W_\ell = W_\ell\setminus\{e\}$, $W_{\ell'} = W_{\ell'}\setminus\{e\}$
                \EndIf
            \EndWhile
            \If{$\lra{S_\ell^t} < \asize_\ell$}
                \State \multiline{Solve $f^t_\ell(B_\ell^t)$ (analogously to~\eqref{eq:second_stage_problem}) considering $\asize_\ell - \lra{S_\ell^t}$ as the right-hand side of the cardinality constraints. Let $\mathbf{z}$ be the optimal solution.} \label{alg-step: solving f}
                \State $S_{\ell}^t = S_{\ell}^t \cup \lrl{\ell' \in B_\ell^t: z_{\ell, \ell'} = 1}$. 
            \EndIf            
        \EndFor
        \State For each $\ell\in \cIJ$, display the subset of profiles $S_\ell^t$. 
        \State Observe like/dislike decisions. Update the sets of potentials and the backlogs following~\eqref{eq: updating formulas}.  \label{alg-step: updating potentials and backlogs}     
      \EndFor
      \end{algorithmic}
    \end{algorithm}      
    }

     In Theorem~\ref{theorem:approx_nonadaptive_policy}, we show that Algorithm~\ref{alg: multiple-period DH} achieves a $1-1/e$ approximation guarantee for Problem~\ref{def: general problem with multiple periods} when $\Pi$ is restricted to semi-adaptive policies. We defer the proof to Appendix~\ref{os:proofs_Tperiods}.
    \begin{theorem}\label{theorem:approx_nonadaptive_policy}
        Algorithm~\ref{alg: multiple-period DH} achieves a $(1-1/e)$-approximation guarantee for Problem \ref{def: general problem with multiple periods} for any platform design and when $\Pi$ is restricted to semi-adaptive policies.
    \end{theorem}
    Note that we could obtain a similar guarantee for the best adaptive policy by solving the linear relaxation of Problem~\eqref{eq: problem DH multiple-periods}, as this would provide an upper bound for Problem~\ref{def: general problem with multiple periods} when restricted to adaptive policies. However, the analysis of this case is more complex since the optimal solution may be fractional and, thus, deriving the guarantee would require rounding techniques. Specifically, the second set of constraints creates a negative correlation among the variables $(\mathbf{x},\mathbf{w},\mathbf{y})$, and it is unclear how to design a randomized rounding method with meaningful guarantees under this setting. However, this is possible when policies are restricted to one-directional interactions and sequential matches, as we show in Theorem~\ref{theorem:approx_nonadaptive_policy2}. 
    \begin{theorem}\label{theorem:approx_nonadaptive_policy2}
        There exists a semi-adaptive policy that achieves a $(1-1/e)$-approximation guarantee for Problem~\ref{def: general problem with multiple periods} when $\Pi$ allows for adaptive policies, but it is restricted to one-directional interactions and sequential matches.
    \end{theorem}
    \begin{remark}
        The guarantee in Theorem~\ref{theorem:approx_nonadaptive_policy2} is with respect to the best possible adaptive policy.
        Our non-adaptive policy solves a relaxed version of Formulation~\eqref{eq: problem DH multiple-periods} with integrality constraints removed. The resulting fractional solution is converted to a feasible solution using the dependent randomized rounding method of~\cite{gandhi06}, which is then fed into the version of Algorithm~\ref{alg: multiple-period DH} adapted to this setting.    
    \end{remark}
    

%% file: 6_experiments.tex

\section{Experiments}\label{sec: experiments}

In this section, we evaluate the performance of the analyzed algorithms (see Table~\ref{tab: summary comparison algorithms} in Appendix~\ref{ec:experiments} for a summary comparison) for different platform designs and compare it with relevant benchmarks.

\subsection{Data}
Our empirical analysis relies on a dataset provided by our industry partner. The sample comprises all heterosexual users located in Houston, TX, who logged in between February 14 and August 14, 2020. For each user, the dataset contains the observable characteristics displayed on their profiles, including age, height, location, race, and religion. It also includes an attractiveness score---hereafter referred to as score---constructed from the number of likes and evaluations a profile has previously received. In addition, the dataset records the complete set of evaluations performed by each user during the observation window, including the decision made (like or dislike), the alternative profiles displayed, and the corresponding timestamps. Thus, we can compute additional outcomes for each user, including their average like probability, backlog size, among others. Table~\ref{tab: descriptives instance} in Appendix~\ref{ec:experiments} reports descriptive statistics of the sample, including the number of users, their average score, the average number of potential matches available, their average backlog size, and their average like probabilities.
Taken together, these features yield a panel dataset that enables a comprehensive characterization of user profiles and evaluation decisions.

\subsection{Benchmarks}
We compare DH and its variants against the following benchmarks:\footnote{We also tested other benchmarks (e.g., Naive and Random). Since the results reported significantly outperform these other benchmarks, we decided to omit them and focus on the results of the algorithms proposed above.}
\begin{enumerate}
  \item Greedy: for each user, select the subset of profiles that maximizes their expected number of matches, i.e.,
        \(
        S_\ell^t = \argmaxA_{S\subseteq \potentials_\ell \setminus \bigcup_{\tau=1}^{t-1} S_\ell^\tau: \lra{S}\leq K_\ell} \lrl{\sum_{\ell'\in S} \plike^t_{\ell, \ell'}\cdot \lrp{\ind{\ell'\in B_\ell^t} + \plike^t_{\ell' \ell}\cdot \ind{\ell'\notin B_\ell^t}} }. \)
  \item Perfect Matching (PM): in each period \(t\in [T]\), solve the perfect matching problem (including possible initial backlogs).
        We formalize this method in Appendix~\ref{app:suboptimal_methods}.
  \item Partner: select the subset of profiles to display following our industry partner's algorithm.
\end{enumerate}
As for variants of DH, we consider its original implementation (as in~\cite{rios2021}), its integer implementation (DH-int, as described in Section~\ref{sec:integral_DH}), and also the non-adaptive heuristic that adapts DH-int to solve one problem considering the entire horizon (DHT, as described in Algorithm~\ref{alg: multiple-period DH}). In all these variants, we assume that like probabilities remain fixed throughout the horizon and only update them once moving to the next period, i.e., whenever solving the problem in a given period $t$, we assume that proxies are $\hat{\phi}_{\ell, \ell'}^\tau = \phi_{\ell, \ell'}^t = \phi_{\ell, \ell'}\lrp{\sigma_\ell^t}$ for all $\ell\in \cIJ, \; \ell' \in \potentials_\ell^t$, and $\tau > t$ (for DHT, we use $\phi_{\ell, \ell'}^t = \phi_{\ell, \ell'}^1 = \phi_{\ell, \ell'}(\sigma_\ell^1)$ for all $t \in [T]$). Note that if $\phi_{\ell,\ell'}(\cdot)$ is ``decreasing with respect to history''---e.g., the linear and threshold approaches described below---, then using $\phi^{t}_{\ell,\ell'}(\sigma_\ell^t)$ to estimate future behavior is an optimistic approach.
Finally, we compare all these methods with the upper bound (UB) obtained from solving Problem~(\ref{eq: upper bound}) using the most favorable like probabilities depending on the history effect.\footnote{For no history, linear and threshold we assume $M_\ell^t = 0$ for all $t\in [T]$; for signal, we assume that all profiles provide the signal of being part of the backlog, (although some do not belong to it); and for disengage, we assume $\gamma = 0$.}

\subsection{Simulation Setup}\label{sec: simulation setup}
For each platform design and benchmark, we perform 100 simulations each involving a time-horizon of a week (i.e., $T=7$). In each period, we (i) choose the profiles to display to each user considering \(K_\ell = 3\) for all \(\ell\in I\cup J\), (ii) simulate the decisions of the users based on their like probabilities, and (iii) update the state of the system before moving on to the next period. In all these simulations, we assume that each user starts with the actual backlog they had in February 14, 2020.

To test the robustness of the proposed methods to different forms of the history effect, we perform our simulations considering different functional forms for the like probabilities. Specifically, for any pair of users $(\ell, \ell')$ with $\ell \in \cIJ$ and $\ell' \in \potentials_\ell^1$ and any period $t\in [T]$, we assume that: 
\[
\phi^t_{\ell,\ell'} = \phi_{\ell,\ell'}(\sigma^t_\ell) = \frac{\exp(u_{\ell,\ell'}+\gamma(\sigma_\ell^t))}{1+\exp(u_{\ell,\ell'}+\gamma(\sigma_\ell^t))},
\] 
where $u_{\ell,\ell'}$ is the indirect utility that $\ell$ gets from matching with $\ell'$ in absence of history effect (see Appendix~\ref{ec:experiments} for estimates and details),
$\sigma_\ell^t$ is the history of user $\ell$ at the beginning of period $t$, and $\gamma(\sigma)$ is a function that maps histories into their effect in the indirect utilities. We consider the following functional forms:
\begin{itemize}
\item[-] {No History}:  $\gamma(\sigma)=0$.
\item[-] {Linear}: $\gamma(\sigma)=\gamma|M|$.
\item[-] {Disengagement}: $\gamma(\sigma)=\gamma t$
\item[-] {Signaling}: $\gamma(\sigma)=-\gamma\,\mathbf{1}\{\ell'\in\backlog_\ell^t\}$
\item[-] {Threshold}:
$\gamma(\sigma) = \begin{cases}
-\infty & |M|>\gamma,\ \gamma\in\mathbb{Z}_+\\
0 & \text{otherwise}
\end{cases}$
\end{itemize}
%
The first approach assumes no history effect on users’ decisions. The second and third approaches, in contrast, allow the number of matches in the recent past to influence liking behavior. The second approach decreases indirect utilities linearly based on past matches (as in~\cite{rios2021}), whereas the third assumes that users automatically dislike all profiles once they exceed a certain threshold of matches, capturing potential time constraints in app usage. The fourth approach models potential disengagement as users spend more time on the app. Finally, the fifth approach accounts for the possibility that users are more likely to like a profile if they know it is already in their backlog, which could lead to an immediate match---a phenomenon analogous to the signaling effect observed in other speed dating settings~\citep{Lee2014}. For implementation, we consider $\gamma = -0.170$ in the linear case (as estimated in~\cite{rios2021}), $\gamma = 5$ in the threshold case, $\gamma = -0.2$ in the disengagement case and signaling cases.
The results are robust to other values of $\gamma$ and to other functional forms of the history effect.



\subsection{Results}
We now present our main simulation results across different platform designs. Section~\ref{sec: results for two-directional interactions} examines the case of two-directional interactions allowing for non-sequential matches,
comparing policies and assessing how sensitive the outcomes are to the functional form of the history effect. Section~\ref{sec: results for one-directional interactions} then turns to the one-directional case, focusing on DH-int under different initiating sides and exploring the robustness of the results to variations in like probabilities and assortment sizes.

\subsubsection{Two-Directional Interactions.}\label{sec: results for two-directional interactions}

\begin{figure}[htp!]
  \caption{Overall Matches varying History Effect}\label{fig: matches varying history effect}
  \begin{subfigure}{0.49\textwidth}
    \caption{Linear}\label{subfig: overall, linear}
    \includegraphics[width=\textwidth]{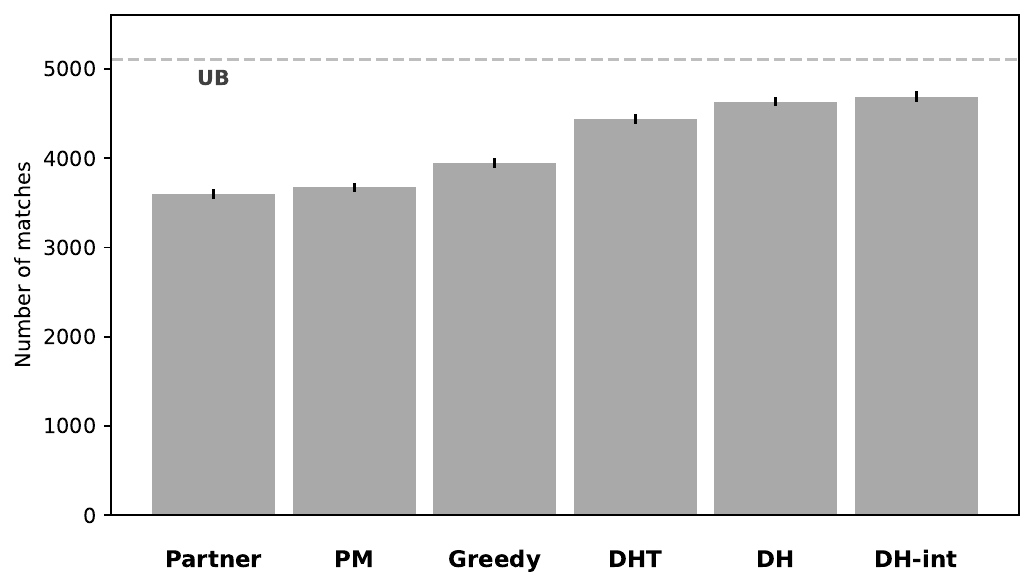}
  \end{subfigure}
  \hspace{0.5cm}
  \begin{subfigure}{0.49\textwidth}
    \caption{Threshold}\label{subfig: overall, threshold}
    \includegraphics[width=\textwidth]{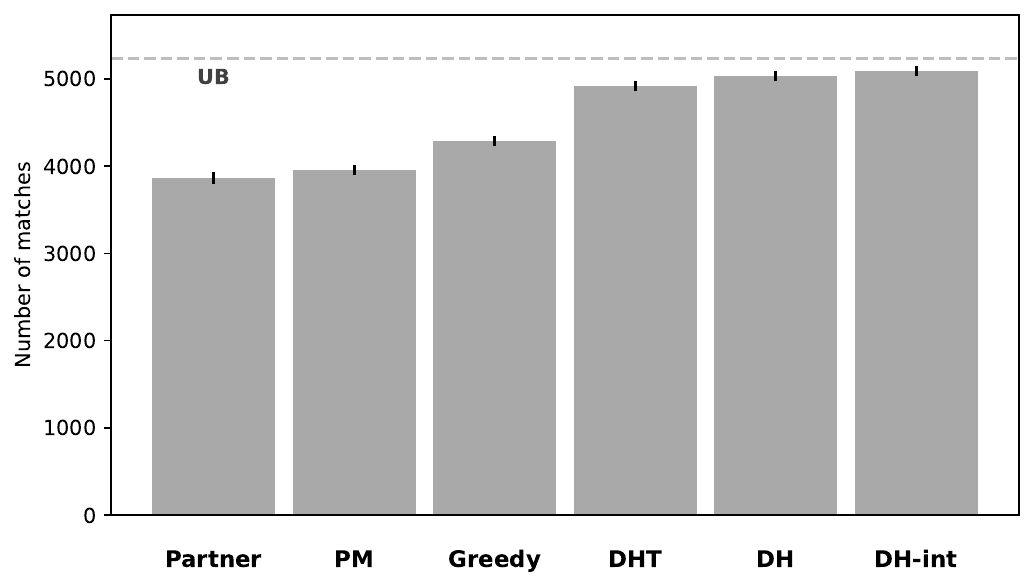}
  \end{subfigure}
  \begin{subfigure}{0.49\textwidth}
    \caption{Disengagement}\label{subfig: overall, disengage}
    \includegraphics[width=\textwidth]{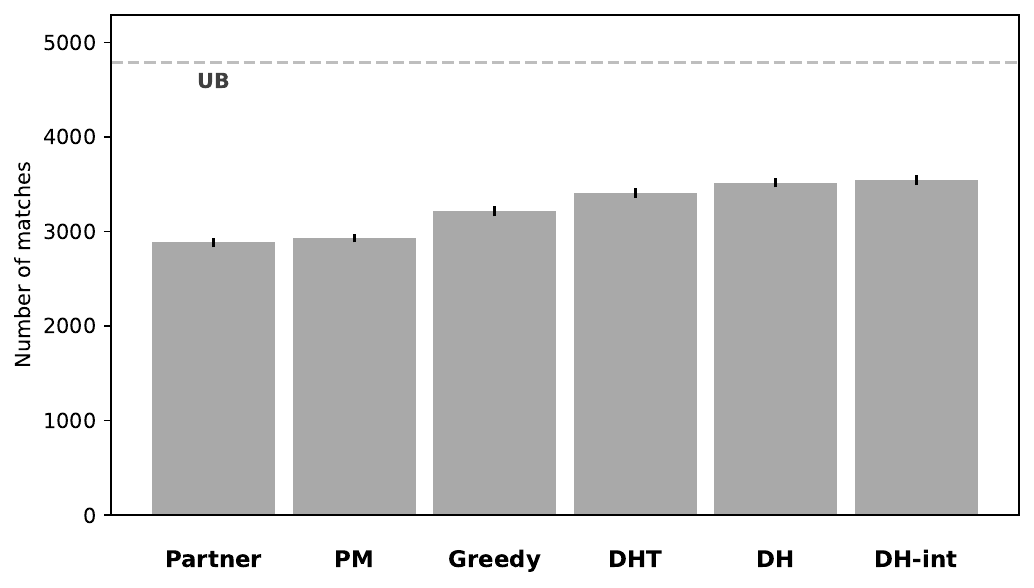}
  \end{subfigure}
  \hspace{0.5cm}
  \begin{subfigure}{0.49\textwidth}
    \caption{Signal}\label{subfig: overall, signal}
    \includegraphics[width=\textwidth]{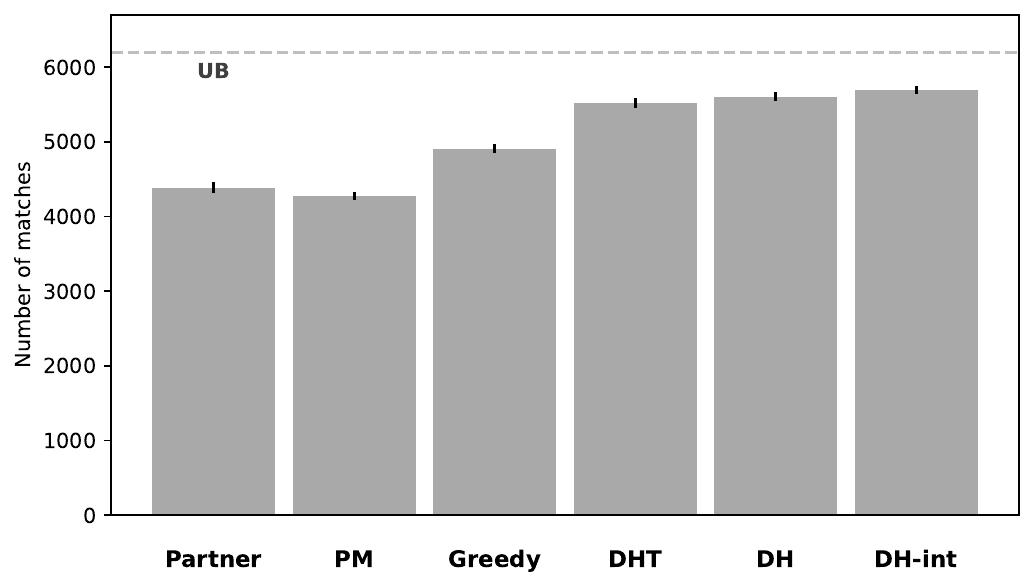}
  \end{subfigure}  
\end{figure}

In Figure~\ref{fig: matches varying history effect}, we report the average (with error bars indicating one standard deviations) number of matches generated by each benchmark across different types of history effects.
First, we observe that the performance of PM and Greedy is considerably better than their worst-case performance. Second, we observe that Greedy outperforms PM and our partner's algorithm for all variants of the history effect. Third, we find that DH and its variants significantly outperform all the other benchmarks, achieving outcomes that are very close to the upper bound. Among them, we note that DH-int performs slightly better than the other DH variants, though the differences are modest. Finally, and more interestingly, we observe that DHT---the semi-adaptive policy that solves a $T$-period variant of the model used by DH---achieves performance comparable to DH and DH-int, despite lacking adaptivity and thus not leveraging realized information on likes and dislikes when selecting subsequent profile subsets.\footnote{We omit the results for no history because they are equivalent to those for threshold.}

\begin{figure}[htp!]
  \caption{Like Rates for Initial and Backlog Interactions}\label{fig: outcomes linear}
  \begin{subfigure}{0.49\textwidth}
    \caption{Initial}\label{subfig: like rate init}
    \includegraphics[width=\textwidth]{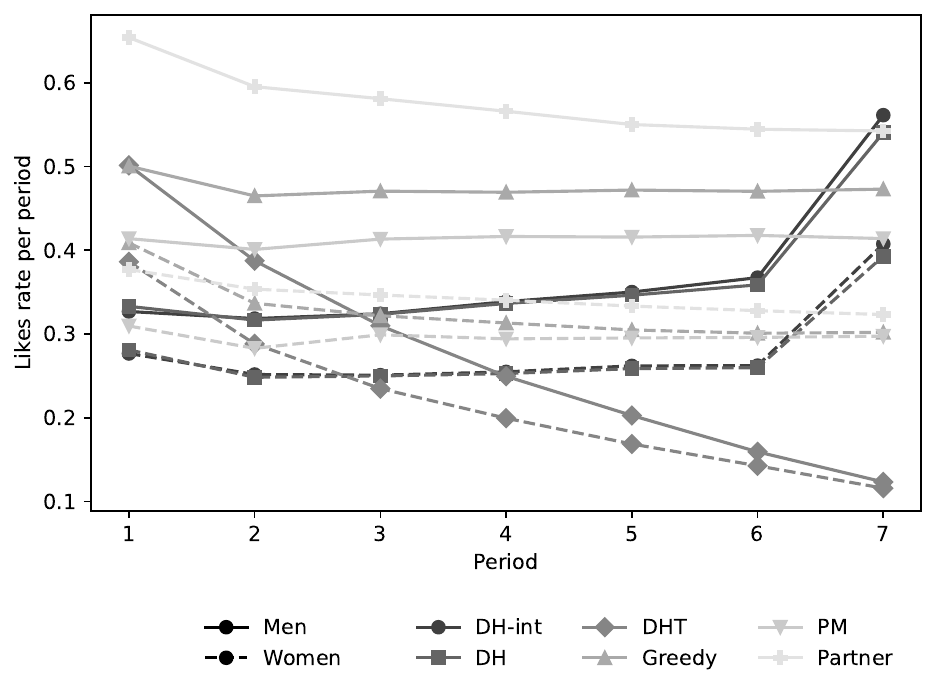}
  \end{subfigure}
  \begin{subfigure}{0.49\textwidth}
    \caption{Backlog}\label{subfig: like rate back}
    \includegraphics[width=\textwidth]{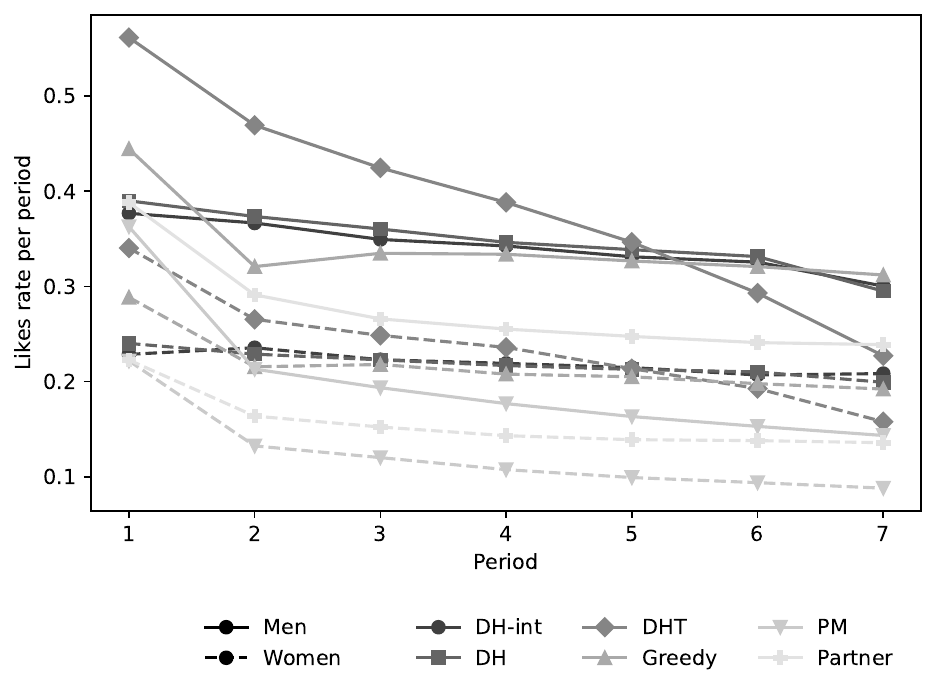}
  \end{subfigure}
\end{figure}

\paragraph{Mechanisms.}
To further investigate the mechanisms underlying these results, Figures~\ref{subfig: like rate init} and~\ref{subfig: like rate back} present the average per-period like rate for profiles initiating a sequential interaction---which we refer to as initial interactions from now on---and backlog profiles, respectively, disaggregated by gender.\footnote{In other words, initial interactions involve any profile not in the user's backlog. Hence, this also includes interactions where both users see each other in the same period.} To ease exposition, we focus the discussion on the linear specification of the history effect; the results are qualitatively similar under the alternative variants described above.

From Figure~\ref{subfig: like rate init}, we observe that the like rate for initial interactions  under DH and its variants is substantially lower than that of the other benchmarks, except in the final period, where DH achieves the highest rate. This pattern arises because DH prioritizes displaying profiles with higher like probabilities in the last period, thereby increasing the likelihood of forming non-sequential matches for users with no backlog. In contrast, Figure~\ref{subfig: like rate back} shows that the like rate for backlog interactions under DH and its variants is consistently higher than that of the other benchmarks and remains relatively stable over time. Taken together, these findings suggest that the gains from DH and its variants stem from carefully balancing the probability of generating an initial like with the enhanced likelihood of securing a subsequent reciprocal like that would result in a match.

\paragraph{Sensitivity to History.}
The results above demonstrate that the effectiveness of DH and its variants is robust across different types of history effects. However, these results rely on a fixed value of the parameter $\gamma$—which varies depending on the type of history effect considered—that determines the magnitude of the effect on like probabilities. To assess the sensitivity of our findings to the strength of the history effect, we conducted simulations for the linear case using a multiplier $\psi \in {0, 0.5, 1.0, 1.5, 2.5, 5.0}$ applied to $\gamma$. In this way, the history effect takes the form $\gamma(\sigma) = \psi \cdot \gamma \cdot \lra{M}$, allowing us to systematically vary its magnitude. Note that $\psi = 0$ leads to the case with no history effect, while $\psi = 1$ is equivalent to the simulations discussed above. Moreover, we focus on DH-int to simplify the exposition; the results are similar for the other variants of DH.

\begin{figure}[htp!]
    \centering
    \caption{Sensitivity to Magnitude of History Effect: Linear Case}
    \label{fig: sensitivity magnitude history effect}    
    \includegraphics[width=0.5\linewidth]{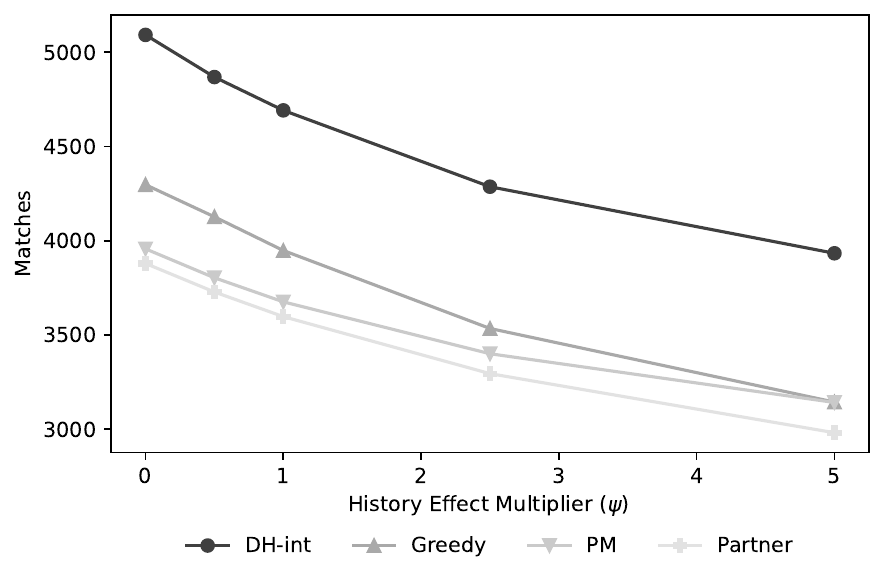}
\end{figure}

As expected, Figure~\ref{fig: sensitivity magnitude history effect} shows that the number of matches decreases with the multiplier $\psi$ across all methods, reflecting the lower like probabilities induced by a stronger history effect. More importantly, DH-int consistently outperforms all benchmarks for every value of $\psi$. Moreover, its relative advantage increases with the intensity of the history effect, rising from about 20\% at $\psi=0.0$ to roughly 25\% at $\psi=5.0$. Together, these findings reinforce the robustness of DH and its variants to both the type and the magnitude of the history effect.

\subsubsection{One-Directional Interactions.}\label{sec: results for one-directional interactions}

\begin{figure}[htp!]
  \caption{Matches per Benchmark by Platform Design}\label{fig: matches by initiating side}
  \begin{subfigure}{0.49\textwidth}
    \caption{Men Initiating}\label{subfig: men initiating}
    \includegraphics[width=\textwidth]{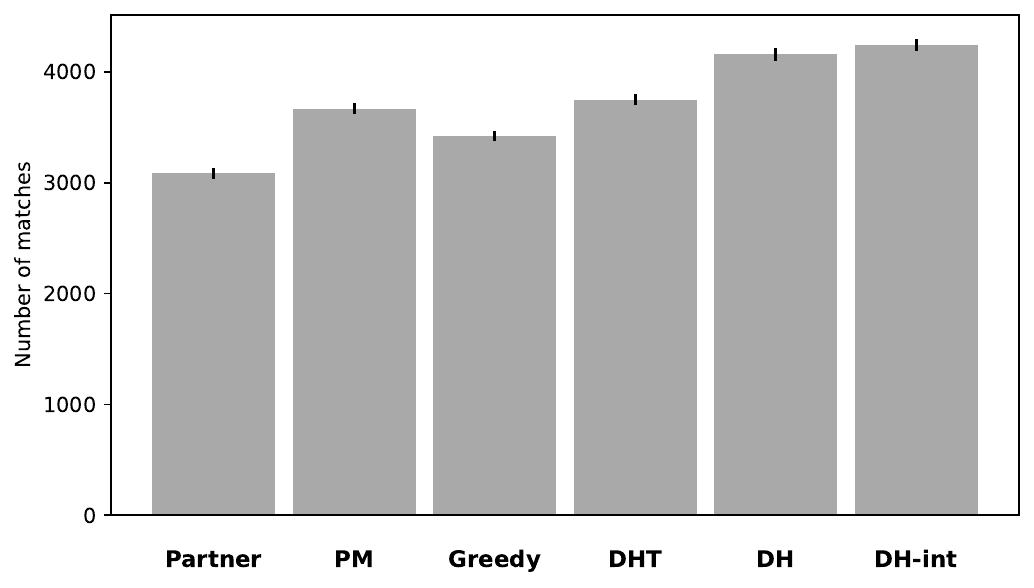}
  \end{subfigure}
  \begin{subfigure}{0.49\textwidth}
    \caption{Women Initiating}\label{subfig: women initiating}
    \includegraphics[width=\textwidth]{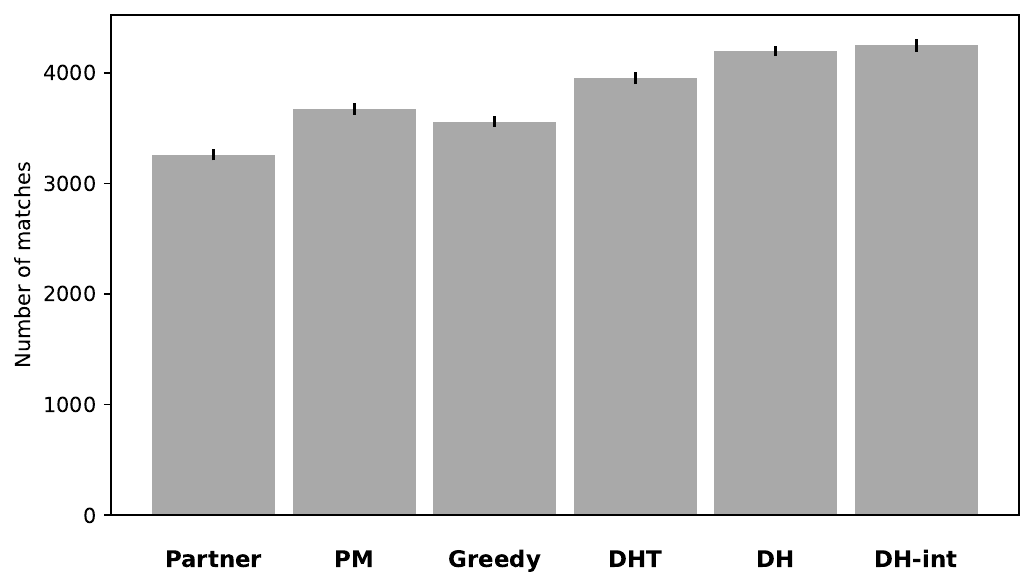}
  \end{subfigure}
\end{figure}

In Figures~\ref{subfig: men initiating} and~\ref{subfig: women initiating}, we report the average (with error bars indicating one standard deviations) number of matches generated by each benchmark under the linear version of the history effect, when men and women are the initiators of each interaction, respectively. As previously discussed, we focus on the linear version of the history effect to ease exposition; the results are similar for other variants of the history effect.

Consistent with the two-directional case, we find that (i) DH and its variants outperform all other benchmarks, and (ii) DH-int delivers a slight advantage over the other variants. Comparing these results with Figure~\ref{subfig: overall, linear}, we also see that one-directional interactions yield a number of matches that is very close to that achieved with two-directional interactions (see Table~\ref{tab:overall_results} in Appendix~\ref{ec:experiments}). 
The intuition is that most matches occur sequentially: while the non-initiating side generates matches from its existing backlog, the initiating side can start new interactions and expand the backlog of the non-initiating side, enabling them to continue evaluating profiles and potentially generating matches in future periods. As a result, the incremental gains from two-directional interactions are limited, provided that users have a backlog or that the time horizon is sufficiently long. Finally, we find that the number of matches produced is generally higher when women initiate interactions, with the improvement being modest for most approaches but substantial under DHT, Greedy and Partner.

%% file: 7_conclusions.tex
\section{Conclusions}\label{sec: conclusions}

In this paper, we study the problem of designing effective policies to determine which subset of profiles to display to each user in each period to maximize the resulting number of matches. We introduced a general theoretical model capturing the dynamic, two-sided nature of the problem, where matches occur only if both users like each other, and users’ past experiences influence their future behavior. By leveraging one-lookahead policies, we propose the Integral Dating Heuristic (DH-int) and provide formal performance guarantees, showing that it achieves a uniform $1-1/e$ approximation across a wide range of platform designs and under various assumptions. {Furthermore, we introduce a semi-adaptive approach,  (DHT) that allows to precompute the profiles to display, reducing the need of re-optimizing in every period. }

Our empirical analysis, conducted using proprietary data from a major U.S.-based dating app, confirms the practical effectiveness of DH-int. Across multiple platform designs and functional forms of the history effects, DH-int consistently outperforms alternative benchmarks, approaching the theoretical upper bound in practice. Our results also demonstrate that careful balancing of initial and follow-up interactions is the key driver of the superior performance of DH-int and its variants, highlighting the importance of incorporating the two-sidedness of the market into policy design. 
Another key observation from our empirical analysis is the robustness of DH-int across a range of history effects. Although the expected total number of matches naturally varies with the different specifications, DH‑int consistently delivers strong performance and small relative gaps to the upper bound for most behaviors; the notable exception is Disengagement, where gaps are larger. This pattern suggests DH‑int is particularly well suited to settings with abrupt, threshold‑type changes.

From a managerial perspective, our findings offer actionable guidance for curated dating platforms. DH-int provides a simple, implementable framework that can substantially improve matching outcomes, increase user satisfaction, and strengthen engagement and retention, all while remaining operationally tractable. More broadly, the insights gained from our study extend beyond dating platforms, offering guidance for other complex, dynamic, two-sided marketplaces—such as freelancing, ride-sharing, and accommodation platforms—where sequencing, allocation, and behavioral dynamics critically shape market outcomes.

%% file: 9_online_supplement.tex
\clearpage
\section{Missing Details in Section~\ref{sec: introduction}}\label{app:intro}

\begin{table}[htpb]
\caption{Comparison of Algorithms}\label{tab: summary comparison algorithms}
\centerline{\scalebox{0.85}{
\begin{tabular}{lccccc}
     \toprule 
     & Pseudo-code & Type & Look-ahead Horizon & Initial-Period Decisions & Later-Period Decisions \\
     \midrule     
     DH & Algorithm~\ref{alg: DH heuristic} & Adaptive & 2 periods & Continuous & Continuous \\
     DH-int & Algorithm~\ref{alg: integral DH heuristic}  & Adaptive & 2 periods & Binary & Continuous \\
     DHT & Algorithm~\ref{alg: multiple-period DH} &  Semi-Adaptive & $T$ periods & Binary & Binary \\
     \bottomrule
\end{tabular}
}\vspace{0.2cm}}
\centerline{\begin{minipage}[t]{16cm}
            \footnotesize \emph{Note:} The column \emph{Initial-Period Decisions} refers to the decision variables $\mathbf{x}^1$, which determine the profiles to show in the initial period of the lookahead horizon. The column \emph{Later-Period Decisions} refer to the decision variables $\mathbf{x}^t$ with $t > 1$, corresponding to decisions made after the initial period of the lookahead horizon.
        \end{minipage}}
\end{table}

%

\section{Missing Details in Section~\ref{sec: model}}\label{app: proofs}
\subsection{Dating Heuristic (DH)}\label{app: dating heuristic}

        In this section, we formalize DH as initially introduced in \citep{rios2021} and formalized in Algorithm~\ref{alg: DH heuristic}.
        For each period \(t\in [T]\), DH considers two steps:
        \begin{enumerate}
          \item {\bf Optimization:} this step involves solving the following linear program:
              \begin{equation}\label{eq: upper bound original DH}
                \begin{split}
                  \max \quad & \;  \sum_{\ell \in I\cup J} \sum_{\ell' \in P^t_\ell} \Bigg[y_{\ell, \ell'}^t\plike^t_{\ell, \ell'} + \frac{1}{2}\cdot w_{\ell, \ell'}^t \plike^t_{\ell, \ell'}\plike^t_{\ell' \ell}\Bigg]+\sum_{\ell \in I\cup J} \sum_{\ell' \in P^t_\ell} \Bigg[y_{\ell, \ell'}^{t+1}\hat{\plike}^{t+1}_{\ell, \ell'} + \frac{1}{2}\cdot w_{\ell, \ell'}^{t+1} \hat{\plike}^{t+1}_{\ell, \ell'}\hat{\plike}^{t+1}_{\ell', \ell}\Bigg] \\
                  \text{s.t.} \quad & \; \sum_{\tau=t}^{t'} y_{\ell, \ell'}^\tau \leq \ind{\ell'\in B^t_\ell} + \sum_{\tau=t}^{t'-1} (x_{\ell' \ell}^\tau - w_{\ell, \ell'}^\tau) \plike^{\tau}_{\ell' \ell}, \hspace{0.5cm} \forall \ell \in I\cup J, \; \ell' \in P^t_\ell, \; t'\in [t,t+1], \\
                  & \; \sum_{\tau=t}^{t+1} x_{\ell, \ell'}^\tau + y_{\ell, \ell'}^\tau \leq 1, \hspace{3.9cm} \forall  \ell \in I\cup J, \; \ell' \in P^t_\ell,\\
                  & \; \sum_{\ell'} x_{\ell, \ell'}^\tau + y_{\ell, \ell'}^\tau \leq K_\ell, \hspace{3.6cm} \forall \ell \in I\cup J, t \in [t,t+1], \\
                  & \; w_{\ell, \ell'}^\tau \leq x_{\ell, \ell'}^\tau, \; w_{\ell, \ell'}^\tau \leq x_{\ell' \ell}^\tau,\; w_{\ell, \ell'}^\tau = w_{\ell' \ell}^\tau, \hspace{1.45cm} \forall \ell \in I\cup J, \; \ell' \in P^t_\ell, \;\tau\in [t,t+1], \\
                  & \; x_{\ell, \ell'}^\tau, \; y_{\ell, \ell'}^\tau, w_{\ell, \ell'}^\tau \in \lrc{0,1}, \hspace{3.4cm} \forall \ell \in I\cup J, \; \ell' \in P^t_\ell, \; \tau\in [t,t+1].
                \end{split}
              \end{equation}
              The decision variables \(y_{\ell, \ell'}^t\) and \(x_{\ell, \ell'}^t\) represent whether \(\ell\) sees profile \(\ell'\) in period \(t\) as part of a backlog and to initiate a sequential match, respectively. The objective is to maximize the expected number of matches obtained in periods \(\lrl{t, t+1}\), including sequential (first term in the objective) and non-sequential matches (second term in the objective). The first family of constraints defines \(y\) and captures the evolution of the backlog. The second family of constraints captures that a profile can be shown at most once, while the third considers the cardinality constraints. Finally, the last family of constraints captures the definition of \(w_{\ell, \ell'}^t\), which accounts for non-sequential matches between \(\ell\) and \(\ell'\) in period \(t\).

          \item {\bf Rounding:} since the optimal decisions \(x^{*,t}, y^{*,t}, w^{*,t}\) of~(\ref{eq: upper bound original DH}) may be fractional, this step involves rounding them in order to decide the profiles to show in the current period. Specifically, the rounding process starts by adding to \(S_\ell^t\) the profiles for which \(y_{\ell, \ell'}^t > 0\) (in decreasing order). Then, if there is space left, the rounding procedure adds to \(S_\ell^t\) the profiles for which \(x_{\ell, \ell'}^t > 0\) (in decreasing order), making sure that the cardinality constraints are satisfied.
        \end{enumerate}
        Notice that these two steps consider the current set of potentials \(P_\ell^t\) and backlog \(B_\ell^t\) for each user \(\ell\in \cIJ\). Then, at the end of each period, the sets of potentials and the backlogs are updated considering the profiles shown and the like/dislike decisions, as shown in~\eqref{eq: updating formulas}.
        
        \begin{algorithm}[h]
        \caption{Dating Heuristic (DH), \citep{rios2021}}\label{alg: DH heuristic}
        \begin{algorithmic}[1]
        \Require An instance of the problem.
        \Ensure A feasible subset of profiles for each user in each period.
        \For{$t\in[T]$}
        \State Solve~\eqref{eq: upper bound original DH}.
         
        \State \multiline{For each \(\ell\), sequentially add profiles \(\ell'\) for which \(y_{\ell, \ell'}^{*,t} > 0\) until the cardinality constraint is binding. If the latter is not binding, add profiles \(\ell'\) for which \(x_{\ell, \ell'}^{*,t} > 0\) until the constraint becomes binding.}
        \State {Update potentials and backlogs following~\eqref{eq: updating formulas}}.
        \EndFor
        \end{algorithmic}
        \end{algorithm}

   \begin{lemma}\label{lemma:equivalent_regions}
      Problem~\eqref{eq: upper bound} is equivalent to Problem~\eqref{eq: upper bound original DH} with variables $\mathbf{y}^{t+1}$ remaining continuous but the rest is constrained to binary values.
    \end{lemma}

    \begin{proof}[Proof of Lemma~\ref{lemma:equivalent_regions}.]
    Let $(\hat{\mathbf{x}}^t,\hat{\mathbf{w}}^t,\hat{\mathbf{y}}^t,\hat{\mathbf{x}}^{t+1},\hat{\mathbf{w}}^{t+1},\hat{\mathbf{y}}^{t+1})$  be a solution in Problem~\eqref{eq: upper bound original DH} where $\hat{\mathbf{y}}^{t+1}$ is the only continuous decision vector and the rest is binary. Then, define $({\mathbf{x}^t},{\mathbf{w}^t},\mathbf{x}^{t+1},\mathbf{w}^{t+1})$ as follows: $\mathbf{x}^{t+1} = \hat{\mathbf{y}}^{t+1}$ and 
 \begin{itemize}
 \item If $\hat{y}_{\ell,\ell'}^t = 1$, then set $x^t_{\ell,\ell'} = 1$, $x^t_{\ell',\ell} = 0$, $w_e^t = 0$ and $w_e^{t+1} = 0$. 
 \item If $\hat{x}^t_{\ell,\ell'} = 1$, $\hat{y}^t_{\ell,\ell'}=\hat{x}^t_{\ell',\ell} =0$ and $\hat{w}^t_{\ell,\ell'} = \hat{w}^t_{\ell',\ell} = 0$, then set ${x}^t_{\ell,\ell'} = 1$, ${x}^t_{\ell',\ell} = 0$ and $w^t_e =w_e^{t+1}= 0$.
 \item If $\hat{x}^t_{\ell,\ell'} = \hat{x}^t_{\ell',\ell} =1$ and $\hat{w}^t_{\ell,\ell'} =  \hat{w}^t_{\ell',\ell} = 1$, then set ${x}^t_{\ell,\ell'} = {x}^t_{\ell',\ell} = 0$, ${w}^t_{e}=1$ and $w_e^{t+1}=0$.
 \item If $\hat{x}^{t+1}_{\ell,\ell'} = \hat{x}^{t+1}_{\ell',\ell} =1$ and $\hat{w}^{t+1}_{\ell,\ell'} = \hat{w}^{t+1}_{\ell',\ell} = 1$, then set ${x}^t_{\ell,\ell'} = {x}^t_{\ell',\ell} = 0$, ${w}^t_{e}=0$ and $w_e^{t+1}=1$.
 \end{itemize}
Clearly $({\mathbf{x}^t},{\mathbf{w}^t},\mathbf{x}^{t+1},\mathbf{w}^{t+1})$ is feasible in Problem~\eqref{eq: upper bound} and it has the same objective value. To show the opposite direction, the construction of variables is analogous.
\end{proof}

\subsection{Missing Proofs}
\begin{proof}[Proof of Proposition~\ref{prop: next problem is not submodular in the general case}]
    For simplicity, consider $T=2$ and no initial backlog. Also, let us write the function $f^{t+1}(R,\mathbf{x}^t,\mathbf{w}^t)$ in \eqref{eq: general next problem} as $f(R)$.
    Let \(I = \lrl{i_1,i_2}, \; J = \lrl{j_1,j_2}\), and the following probabilities:
    \[\plike_{i_1,j_1} = 1, \plike_{j_1,i_1} = \epsilon, \plike_{i_2,j_2} = \epsilon, \plike_{j_2,i_2} = 1, \beta_{i_1,j_2} = 1/2, \beta_{i_2,j_1} = 1/2.\]
    Define \(X(R)\) and \(W(R)\) as the sets of backlog and non-backlog pairs shown in the optimal solution of the second period problem.
    If \(R = \emptyset\), then \(W(R) = \lrl{(i_1,j_2), (i_2,j_1)}\). If we add \((i_2,j_2)\) to $R$, then
    \(W(R \cup {(i_2,j_2)}) = \lrl{(i_2,j_1)}\) and \(X(R \cup {(i_2,j_2)}) = \lrl{(j_2,i_2)}\).
    Hence,
    \[f(R \cup {(i_2,j_2)}) - f(R) = 1+1/2 - (1/2 + 1/2) = 1/2.\]

    On the other hand, if \(R' = \lrl{(j_1,i_1)}\), \(W(R') = \lrl{(i_2,j_1)}\) and \(X(R') = \lrl{(i_1,j_1)}\).
    If we add \((i_2,j_2)\) to \(R'\), then \(W(R'\cup {(i_2,j_2)}) = \lrl{(i_2,j_1)}\) and \(X(R'\cup {(i_2,j_2)}) = \lrl{(i_1,j_1), (j_2,i_2)}\). Then,
    \[f(R' \cup {(i_2,j_2)}) - f(R') = 1+1+1/2 - (1+1/2) = 1.\]
    Hence, we have that
    \[
    R \subset R' \quad  \text{ and } \quad f(R \cup {(i_2,j_2)}) - f(R) < f(R' \cup {(i_2,j_2)}) - f(R'),
    \]
    so we conclude that \(f(R)\) is not submodular. 
\end{proof}

\begin{proposition}\label{prop: second stage problem solved with LP}
          The function $f^{t+1}$ defined in \eqref{eq: general next problem} can be efficiently evaluated by solving a linear program.
\end{proposition}

\begin{proof}[Proof of Proposition~\ref{prop: second stage problem solved with LP}.]
For simplicity, consider $T=2$, therefore $f^{t+1}$ corresponds to $f^2$.
Given a set of realized arcs and edges   \(R\), backlogs $\mathbf{B}=\{B_\ell\}_{\ell\in I\cup J}$ and a set of potentials \(\bm\potentials\), define a bipartite graph with two sides \(U,V\) with \(U = I\cup J\) and \(V = R\cup \lrl{(i,j) \in I\times J\; : \; j\in \potentials_i, i\in \potentials_j}\), i.e., \(U\) contains the set of users and \(V\) the set of arcs that could be displayed in the second period.
  Let \(Q\subseteq U\times V\) be the set of edges. Then, a pair \(\lrp{\ell, (\ell',\ell'')} \in U \times V\) belongs to \(Q\) if and only if
  \[
  (\ell', \ell'') \in R \quad \text{ or } \quad \lrc{ (\ell', \ell'') \notin R, \ell\in \lrl{\ell', \ell''}}.
  \]
  In words, an edge between \(\ell\in U\) and \((\ell', \ell'')\in V\) exists if and only if the edge \((\ell', \ell'')\) is either in the backlog of \(\ell\) or both users \(\ell\) and \(\ell'\) can see each other simultaneously.
  Now, for any pair of nodes \((u,v) \in Q\) such that \(v\in R\), we define a variable \(y_{u,v}\) that is equal to 1 if \(v=(u,u') \in R\) and \(u\) sees \(u'\), and zero otherwise. Similarly, for any pair \((u,v)\in Q\) such that \(v=(u,u') \in V\setminus R\), we define a variable \(x_{u,v}\) that is equal to 1 if \(u\) sees \(u'\), and zero otherwise. Note that here we do a slight abuse of notation and assume that if \(x_{i, (i,j)} = 1\), then \(i\in I\) sees profile \(j\in P_i\setminus B_i\) and, similarly,
  if \(x_{j, (i,j)} = 1\), then \(j\in J\) sees \(i\in P_j\setminus B_j\). For convenience, we will use in these cases that \(i\in (i,j)\) and \(j\in (i,j)\).
  Finally, for any pair \((u,v)\in Q\) such that \(u \in I\) and \(v \in V\setminus R\), let \(w_{u,v}=1\) if both users involved in \(v\) see each other simultaneously, i.e., if \(v = (i,j) \in V\setminus R\), then \(x_{i,(i,j)} = x_{j, (i,j)} = 1\).

  Using these variables, we can formulate the second period problem as follows:
\begin{subequations}
\label{eq: second period problem integer}
\begin{alignat}{3}
\max \quad & \quad \sum_{u\in U} \sum_{v \in R} y_{u,v}\cdot \phi_v^2 + \sum_{i\in I} \sum_{\delta(i) \cap V\setminus R} w_{i,v}\cdot \beta_v^2 \label{eq: second period objective}\\
{s.t.} \quad & \quad \sum_{v \in \delta(u)\cap R} y_{u,v} + \sum_{v \in \delta(u)\cap V\setminus R} x_{u,v} \leq K_u, \quad \forall u \in U \label{eq: second period capacity} \\
& \quad w_{i,(i,j)} - x_{u,(i,j)} \leq 0, \quad \forall (i,j) \in V\setminus R, i \in I, u\in (i,j) \label{eq: second period mccormick1} \\
& \quad w_{i,(i,j)} - x_{u,(i,j)} \leq 0, \quad \forall (i,j) \in V\setminus R, i \in I, u\in (i,j) \label{eq: second period mccormick2} \\
& \quad x_{u,v} \in \lrl{0,1}, \quad v\in V\setminus R, u \in v \\
& \quad y_{u,v} \in \lrl{0,1}, \quad \forall u\in U, v\in \delta(u) \cap R \\
& \quad w_{i,v} \in \lrl{0,1}, \quad \forall v\in V\setminus R, i \in v
\end{alignat}
\end{subequations}
  where \(\beta_v^2 = \phi_v^2 \cdot \phi_{\overline{v}}^2\) is the match probability between the users in the pair \(v\) and \(\delta(u)\) is the set of edges incident to node \(u\in U\).
  Let \(A^2\) be the set of constraints in~(\ref{eq: second period capacity}),~(\ref{eq: second period mccormick1}) and~(\ref{eq: second period mccormick2}).
  Note that each variable appears at most twice in \(A^2\), and that every time they appear they are multiplied by either $1$ or $-1$.
  Thus, to show that the matrix of constraints if totally unimodular, it remains to show that the constraints in \(A^2\) can be separated in two subsets such that (i) if a variable appears twice with different signs, then the constraints belong to the same subset, and (ii) if a variable appears twice with the same sign, then the constraints belong to different subsets.
  Let \(A^2_I\) and \(A^2_J\) be the subsets of constraints of \(A^2\) involving \(u\in I\) and \(u\in J\), respectively. Then, observe that
  \begin{itemize}
    \item each \(w_{i,(i,j)}\) appears in two constraints with the same sign (\(w_{i,(i,j)} - x_{i,(i,j)} \leq 0\) and \(w_{i,(i,j)} - x_{j,(i,j)} \leq 0\)), but these constraints belongs to \(A_I^2\) and \(A_J^2\), respectively.
    \item each \(x_{i,(i,j)}\) appears in two constraints with different signs (\(\sum_{v \in \delta(i)\cap B_i} y_{i,v} + \sum_{v \in \delta(i)\cap V\setminus B_i} x_{i,v} \leq K_i\) and \(w_{i,(i,j)} - x_{i,(i,j)} \leq 0\)), but these constraints belong both to the same subset \(A_I^2\). Similarly, \(x_{j,(i,j)}\) appears in two constraints with different signs, but both constraints belong to \(A_J^2\).
  \end{itemize}
  Hence, using Hoffman's sufficient condition, we conclude that the constraints in \(A^2\) can be written as the product of a totally unimodular matrix and our vector of decisions variables. Finally, since the right-hand side of the constraints are integral, we conclude that the feasible region of the problem is an integral polyhedron, and thus we can solve its linear relaxation. 
\end{proof}


\section{Missing Details in Section~\ref{sec: analysis two period model}}\label{os:proofs_2periods}
\begin{proof}[Proof of Lemma~\ref{lemma: f monotone submodular}.]
Let us construct a partition of the edges $\vec{E}^1$: Define part $\mathcal{E}_\ell$ for each $\ell\in \cIJ$ as the set $\{(\ell,\ell')\in\vec{E}^1: (\ell,\ell')\in R\}$. The budget for each part is $K_\ell$. This shows that the feasible region corresponds to a partition matroid. Then, using Proposition 3.1 in \cite{fisher78}, we conclude that $f(R)$ is monotone and submodular, because the objective is a linear function.
\end{proof}

\begin{definition}[Monotonicity and Submodularity]\label{def:monotonicity}
    A generic non-negative set function $g:\{0,1\}^{\mathcal{U}}\to \RR_+$ defined over the space of elements $\mathcal{U}$ is {\it monotone} if, for every $\mathbf{x}\in\{0,1\}^{\mathcal{U}} $ and $e\in\mathcal{U}$, we have $g(\mathbf{x}+\ind{e})\ge g(\mathbf{x})$, where $\ind{e}\in\{0,1\}^{\mathcal{U}}$ is the indicator vector with value 1 in component $e$ and zero elsewhere. Moreover, $g$ is \emph{submodular} if for every
$e,e'\in\mathcal{U}$, we have $g(\mathbf{x}+\ind{e}+\ind{e'})-g(\mathbf{x}+\ind{e'})\le g(\mathbf{x}+\ind{e})-g(\mathbf{x})$.
\end{definition}


\begin{lemma}\label{lemma:general_function_submodular_monotone}
The function $\Ma^2(\mathbf{x}, \mathbf{w})$ is non-negative monotone and submodular in $\mathbf{x}$.
\end{lemma}

\begin{proof}[Proof of Lemma~\ref{lemma:general_function_submodular_monotone}.]
Since $\Ma^2(\mathbf{x}, \mathbf{w})$ does not directly depend on $\mathbf{w}$, we will drop it from the notation to ease exposition and simply write $\Ma^2(\mathbf{x})$.

Recall that $f$ is a non-negative, monotone submodular function due to Lemma~\ref{lemma: f monotone submodular}. Since $f$ is non-negative, then $\Ma^2$ is also non-negative.
For $\mathbf{x}\in\{0,1\}^{\vec{E}^1}$, let $\vec{E}^1(\mathbf{x}) = \{e\in\vec{E}^1: \ x_e = 1, x_{\overline{e}}= 0\}$ where $\overline{e}$ denotes the inverted arc $e$, i.e., if $e=(\ell,\ell')$ then $\overline{e}=(\ell',\ell)$. \footnote{Note that, when $\mathbf{x}$ is feasible in Problem \ref{def:two_period_problem}, for any pair $i\in I$, $j\in\potentials^1_i$ we cannot have $(i,j), (j,i)\in\vec{E}^1(\mathbf{x})$ because of one of the constraint of the problem.}  Then, we can re-write $\Ma^2$ as follows
\[
\Ma^2(\mathbf{x}) = \sum_{R\subseteq \vec{E}^1(\mathbf{x})}f(R)\prod_{e\in R}\phi_e^1\prod_{e\in\vec{E}^1(\mathbf{x})\setminus R}(1-\phi_e^1)
\]
On the other hand, for any $e\notin\vec{E}^1(\mathbf{x})$ we have
\begin{align*}
\Ma^2(\mathbf{x}+\ind{e}) &= \phi_e^1 \cdot\sum_{R\subseteq \vec{E}^1(\mathbf{x})}f(R\cup e)\prod_{e\in R}\phi_e^1\prod_{e\in\vec{E}^1(\mathbf{x})\setminus R}(1-\phi_e^1) \\
&+ (1-\phi_e^1)\cdot\sum_{R\subseteq \vec{E}^1(\mathbf{x})}f(R)\prod_{e\in R}\phi_e^1\prod_{e\in\vec{E}^1(\mathbf{x})\setminus R}(1-\phi_e^1)
\end{align*}
Therefore,
\[
\Ma^2(\mathbf{x}+\ind{e}) - \Ma^2(\mathbf{x}) = \phi_e^1\sum_{R\subseteq \vec{E}^1(\mathbf{x})}[f(R\cup e)-f(R)]\prod_{e\in R}\phi_e^1\prod_{e\in\vec{E}^1(\mathbf{x})\setminus R}(1-\phi_e^1)
\]
Since $f$ is monotone, then $f(R\cup e)-f(R)\geq 0$ for all $R\subseteq \vec{E}^1(\mathbf{x})$ and $e\notin\vec{E}^1(\mathbf{x})$, which implies $\Ma^2(\mathbf{x}+\ind{e}) - \Ma^2(\mathbf{x})\geq0$.

Now, let us prove that $\Ma^2$ is submodular. Consider any $\mathbf{x}\in\{0,1\}^{\vec{E}^1}$ and $e,e'\notin\vec{E}^1(\mathbf{x})$. Our goal is to show an alternative characterization of submodularity: $\Ma^2(\mathbf{x}+\ind{e})-\Ma^2(\mathbf{x})\geq \Ma^2(\mathbf{x}+\ind{e}+\ind{e'})-\Ma^2(\mathbf{x}+\ind{e'})$. Note that we have the following expression for $\Ma^2(\mathbf{x}+\ind{e}+\ind{e'})$
\begin{align*}
\Ma^2(\mathbf{x}+\ind{e}+\ind{e'}) &= \phi_e^1\phi_{e'}^1 \cdot\sum_{R\subseteq \vec{E}^1(\mathbf{x})}f(R\cup e,e')\prod_{e\in R}\phi_e^1\prod_{e\in\vec{E}^1(\mathbf{x})\setminus R}(1-\phi_e^1) \\
&+ \phi_{e'}^1(1-\phi_e^1)\cdot\sum_{R\subseteq \vec{E}^1(\mathbf{x})}f(R\cup e')\prod_{e\in R}\phi_e^1\prod_{e\in\vec{E}^1(\mathbf{x})\setminus R}(1-\phi_e^1)\\
&+ \phi_{e}(1-\phi_{e'}^1)\cdot\sum_{R\subseteq \vec{E}^1(\mathbf{x})}f(R\cup e)\prod_{e\in R}\phi_e^1\prod_{e\in\vec{E}^1(\mathbf{x})\setminus R}(1-\phi_e^1)\\
&+ (1-\phi_{e'}^1)(1-\phi_e^1)\cdot\sum_{R\subseteq \vec{E}^1(\mathbf{x})}f(R)\prod_{e\in R}\phi_e^1\prod_{e\in\vec{E}^1(\mathbf{x})\setminus R}(1-\phi_e^1)\\
\end{align*}
Analogously, we can compute $\Ma^2(\mathbf{x}+\ind{e})$ and $\Ma^2(\mathbf{x}+\ind{e'})$. By deleting common terms, we can obtain the following
\begin{align*}
&\Ma^2(\mathbf{x}+\ind{e})-\Ma^2(\mathbf{x})- \Ma^2(\mathbf{x}+\ind{e}+\ind{e'})+\Ma^2(\mathbf{x}+\ind{e'}) \\
&= \phi_e^1\phi_{e'}^1\cdot \sum_{R\subseteq \vec{E}^1(\mathbf{x})}[f(R\cup e)-f(R)-f(R\cup e,e')+f(R\cup e')]\cdot\prod_{e\in R}\phi_e^1\prod_{e\in\vec{E}^1(\mathbf{x})\setminus R}(1-\phi_e^1),
\end{align*}
from which submodularity follows due to submodularity of $f$. 
\end{proof}

\subsection{Proof of Theorem~\ref{thm: guarantee for DH}}\label{os:proof_maintheorem_2periods}
Let $(\mathbf{x}^{1,\star},\mathbf{w}^{1,\star},\mathbf{x}^{2,\star})$ be an optimal solution of Problem~\ref{def:two_period_problem}, which we denote as $(\mathbf{x}^{\star},\mathbf{w}^{\star},\mathbf{y}^{\star})$ to minimize notation. Consider the following formulation
\begin{align}\label{eq:P1prime}
F(\mathbf{x}^\star, \mathbf{w}^\star) := \max & \quad \sum_{\ell \in \cIJ} \sum_{\ell'\in \potentials_\ell^1}\phi_{\ell, \ell'}^2 \cdot y_{\ell, \ell'}\\
s.t.&\quad y_{\ell, \ell'}\leq \phi^1_{\ell', \ell}\cdot x^\star_{\ell', \ell} , \hspace{2.5cm} \forall \ell \in \cIJ, \;\ell' \in \potentials_\ell^1 \notag\\
&\quad \sum_{\ell' \in \potentials_\ell^1} y_{\ell, \ell'}\leq K_\ell, \hspace{2.6cm}\forall \ell \in \cIJ\notag\\
&\quad y_{\ell, \ell'}\geq 0,  \hspace{3.7cm}\forall \ell \in \cIJ, \ell' \in \potentials_\ell^1.\notag
\end{align}	
Clearly, $\mathbf{y}^\star$ is an optimal solution in $F(\mathbf{x}^\star, \mathbf{w}^\star)$. Define $\mathcal{X}_\ell = \{\ell'\in\potentials^1_\ell: \ x^\star_{\ell',\ell}=1\}$ and observe that, in the problem above, we can restrict variables $y_{\ell,\ell'}$ to $\ell'\in\mathcal{X}_\ell$ (the others are zero). Also, note that \eqref{eq:P1prime} is separable, i.e., $F(\mathbf{x}^\star, \mathbf{w}^\star)= \sum_{\ell\in I\cup J} F_\ell(\mathbf{x}^\star, \mathbf{w}^\star)$, where 
\begin{equation}\label{eq:P1prime_peruser}
F_\ell(\mathbf{x}^\star, \mathbf{w}^\star):= \max\left\{\sum_{\ell'\in \mathcal{X}_\ell}\phi_{\ell, \ell'}^2 \cdot y_{\ell'}: \ \sum_{\ell'\in\mathcal{X}_\ell}y_{\ell'}\leq K_\ell, \ 0\leq y_{\ell'}\leq \phi^1_{\ell',\ell}, \;\forall \; \ell'\in \mathcal{X}_\ell\right\}.
\end{equation}
Our goal is to compare $F(\mathbf{x}^\star, \mathbf{w}^\star)$ with the \emph{distribution problem}, which considers a larger space of distributions (possibly with correlations) over the arcs in $R$. For $\mathbf{x}^\star,\mathbf{w}^\star$, the {distribution problem} is defined as:         
\begin{align*}
G(\mathbf{x}^\star, \mathbf{w}^\star) := \max \Bigg\{ \sum_{R \subseteq \vec{E}} f(R)\cdot\lambda_{R}: &\sum_{R\subseteq \vec{E}^1}\lambda_{R} = 1, \ \sum_{\substack{R\subseteq \vec{E}^1:\\(\ell',\ell)\in R}}\lambda_{R} = \phi_{\ell',\ell}\cdot x^\star_{\ell',\ell}, \quad\forall (\ell',\ell)\in\vec{E}^1\notag \\
&\hspace{5cm}  \lambda_{R}\geq 0, \quad \forall R \subseteq \vec{E}^1.\Bigg\}\notag
\end{align*}
where $\lambda_R$ denotes the probability that the subset of arcs resulting from sequential likes is $R\subseteq\vec{E}^1$. The second constraint enforces that the marginal probabilities coincide with the probability that $\ell'$ sees and likes $\ell$. Instead of using the formulation above, we will construct an equivalent one such that the objective is separable as we did with $F(\mathbf{x}^\star,\mathbf{w}^\star)$. For this, we use a backlog $B$ for each user instead of the overall set of realized arcs $R$. In fact, the backlog $B\subseteq\potentials^1_\ell$ of each user $\ell\in I\cup J$ corresponds to $B=\{\ell'\in\potentials^1_\ell: (\ell',\ell)\in R\}$. Therefore, we have
\begin{align*}
G(\mathbf{x}^\star, \mathbf{w}^\star) := \max \Bigg\{\sum_{\ell \in \cIJ} \sum_{B \subseteq \potentials_\ell^1} f_\ell(B)\cdot\lambda_{\ell,B}:&\ \sum_{B\subseteq \potentials_\ell^1}\lambda_{\ell,B} = 1, \hspace{2cm} \forall \ell \in \cIJ\notag, \\
& \sum_{\substack{B\subseteq \potentials_\ell^1:\\\ell'\in B}}\lambda_{\ell,B} = \phi_{\ell',\ell}\cdot x^\star_{\ell',\ell}, \hspace{0.5cm} \forall \ell \in \cIJ, \; \ell' \in \potentials_\ell^1, \\  &\lambda_{\ell,B}\geq 0, \hspace{3cm} \forall \ell\in \cIJ, \; B \subseteq \potentials_\ell^1.\Bigg\}\notag
\end{align*}
where $\lambda_{\ell, B}$ can be interpreted as the probability that the backlog of user $\ell$ is $B$. The second family of constraints states that the probability of the backlog of $\ell$ containing $\ell'$ equals the marginal probability that $\ell'$ saw and liked $\ell$ in the first period. 
Similar to $F(\mathbf{x}^\star, \mathbf{w}^\star)$, note that $G(\mathbf{x}^\star, \mathbf{w}^\star)$ is separable since for each user $\ell$ we have that the expected total number of matches that $\ell$ can achieve in the second period given backlog $B$ is $f_\ell(B):=\max\{\sum_{\ell'\in S}\phi^2_{\ell,\ell'}: \ |S|\leq K_\ell, \; S\subseteq B\}$. Also, we can limit variables $\lambda_{\ell,B}$ to subsets $B\subseteq\mathcal{X}_\ell$ because $x^\star_{\ell',\ell} = 0$ implies that $\lambda_{\ell,B} = 0$ for all $B$ such that $\ell'\in B$. Hence, we can write $G(\mathbf{x}^\star, \mathbf{w}^\star) = \sum_{\ell\in I\cup J}G_\ell(\mathbf{x}^\star, \mathbf{w}^\star)$, where
\begin{equation}\label{eq:P2prime_peruser}
G_\ell(\mathbf{x}^\star, \mathbf{w}^\star):= \max\left\{\sum_{B\subseteq \mathcal{X}_\ell}f_\ell(B)\cdot\lambda_B: \ \sum_{B\subseteq\mathcal{X}_\ell}\lambda_B = 1, \ \sum_{\substack{B\subseteq \mathcal{X}_\ell\\\ell'\in B}}\lambda_B = \phi^1_{\ell',\ell}, \;\forall \; \ell'\in \mathcal{X}_\ell,\; \bm\lambda\geq0\right\}.
\end{equation}

To show our main result, consider the dual formulation of~\eqref{eq:P1prime_peruser} and~\eqref{eq:P2prime_peruser} which are given by (resp.): 
\begin{align}\label{eq:P1prime_dual_peruser}
F^D_\ell(\mathbf{x}^\star, \mathbf{w}^\star) := \min & \quad  K_\ell\cdot \theta + \sum_{\ell' \in \mathcal{X}_\ell} \gamma_{\ell'}\cdot \phi^1_{\ell', \ell}\\
s.t. & \quad \theta+\gamma_{\ell'}\geq \phi^2_{\ell,\ell'} \qquad \forall \ell' \in \mathcal{X}_\ell, \notag\\
&\quad \theta,\gamma_{\ell'}\geq0, \hspace{1.25cm} \forall  \ell'\in \mathcal{X}_\ell\notag
\end{align}
and
\begin{align}\label{eq:P2prime_dual_peruser original}
G^D_\ell(\mathbf{x}^\star, \mathbf{w}^\star) := \min & \quad \overline{\theta} + \sum_{\ell' \in \mathcal{X}_\ell} \phi^1_{\ell',\ell}\cdot \overline{\gamma}_{\ell'}\\
s.t.&\quad \overline{\theta} + \sum_{\ell' \in B}\overline{\gamma}_{\ell'} \geq f_\ell(B) \qquad  \forall B \subseteq \mathcal{X}_\ell\notag\\ 
&\quad \overline{\theta}, \; \overline{\gamma}_{\ell'}\in\RR, \hspace{2cm}\forall  \ell'\in \mathcal{X}_\ell.\notag
\end{align}


In Lemma~\ref{lemma: properties dual distribution problem individual}, we characterize some useful properties of the dual problem $G^D_\ell(\mathbf{x}^\star, \mathbf{w}^\star)$. We defer its proof to Appendix~\ref{os:remaining_lemmas_2periods} below.
\begin{lemma}\label{lemma: properties dual distribution problem individual}
Let $G^D_\ell(\mathbf{x}^\star, \mathbf{w}^\star)$ be as defined in~\eqref{eq:P2prime_dual_peruser original}. Then, the following properties hold:
\begin{enumerate}
\item There is an optimal solution $(\overline{\theta}^\star,\overline{\gamma}^\star)$ such that $\overline{\theta}^\star, \; \overline{\gamma}_{\ell'}^\star \geq 0,$ for all $\ell'\in \mathcal{X}_\ell$.
\item The feasible region can be equivalently restricted to constraints for $B \subseteq \mathcal{X}_\ell$ such that $|B|\leq K_\ell$.
\item There is an optimal solution $(\overline{\theta}^\star,\overline{\gamma}^\star)$ and a subset $B^\star\subseteq\mathcal{X}_\ell$ with $|B^\star| = K_\ell$ such that, for some $\alpha^\star\geq 0$,  we have: (i) $\phi^2_{\ell,\ell'}-\overline{\gamma}_{\ell'}^\star = \alpha^\star$ for all $\ell'\in B^\star$, (ii) $\phi^2_{\ell,\ell'}-\overline{\gamma}_{\ell'}^\star \leq \alpha^\star$ for all $\ell'\notin B^\star$, and (iii) $\overline{\theta}^\star = K_\ell\alpha^\star$.\label{sublemma: solution dual distribution individual alpha}
\end{enumerate}
\end{lemma}
Thanks to Lemma~\ref{lemma: properties dual distribution problem individual}, we can show the following:
\begin{lemma}\label{lemma: dual distribution problem is upper bound}
$F(\mathbf{x}^\star, \mathbf{w}^\star) \leq G(\mathbf{x}^\star, \mathbf{w}^\star)$.
\end{lemma} 
This is key in the proof of Theorem~\ref{thm: guarantee for DH} as it shows that the optimal solution of the distribution problem provides an upper bound for the problem solved by DH-int for $T=2$. We defer the proof of Lemma~\ref{lemma: dual distribution problem is upper bound} to Appendix~\ref{os:remaining_lemmas_2periods} below.

\begin{proof}[Proof of Theorem~\ref{thm: guarantee for DH}.]
Let us recall for a moment that the original objective function of our problem can be re-written as
\[
 \sum_{e\in E^1} \beta^1_{e}\cdot w_{e} + \sum_{\ell\in \cIJ} \sum_{B\subseteq \potentials_\ell^1} f_\ell(B)\cdot \probP_{\mathbf{x}^1}(B),
\]
where $\probP_{\mathbf{x}^1}(B) = \prod_{\ell'\in B}\phi^1_{\ell',\ell}\cdot x_{\ell',\ell} \prod_{\ell'\notin B}(1-\phi^1_{\ell',\ell}\cdot x_{\ell',\ell})$\
and $f_\ell(B) = \max_{S\subseteq B} \lrl{\sum_{\ell' \in S} \phi_{\ell, \ell'}^2 : \lra{S} \leq K_\ell }$. 
In particular, given an optimal solution $(\mathbf{x}^\star,\mathbf{w}^\star,\mathbf{y}^\star)=(\mathbf{x}^{1,\star},\mathbf{w}^{1,\star},\mathbf{x}^{2,\star})$ of~\eqref{eq:upper_bound_twoperiods} $(\mathbf{x}^\star,\mathbf{w}^\star,\mathbf{y}^\star)$, we can define $\lambda^{\sf ind}_{\ell,B}= \probP_{\mathbf{x}_\ell^\star}(B)$, which is a feasible solution in \eqref{eq:P2}. The correlation gap, introduced in \citep{agrawal2010correlation}, lower bounds the ratio between the objective value of the independent distribution $\lambda^{\sf ind}_{\ell,B}$ and the optimal value in \eqref{eq:P2}. Formally, let $\lambda^\star_{\ell, B}$ be an optimal solution in \eqref{eq:P2}. Since $f_\ell(\cdot)$ is a monotone submodular function (similar to Lemma~\ref{lemma: f monotone submodular}) for each $\ell$, we know that
\[
\frac{\sum_{\ell\in I\cup J}\sum_{B\subseteq\potentials^1_\ell}f_\ell(B)\cdot \lambda^{\sf ind}_{\ell,B}}{\sum_{\ell\in I\cup J}\sum_{B\subseteq\potentials^1_\ell}f_\ell(B)\cdot \lambda^{\star}_{\ell,B}} \geq 1-1/e.
\]
Note that the numerator is equivalent to $\Ma^2(\mathbf{x}^\star,\mathbf{w}^\star)$ (as defined in~\eqref{eq: definition of expectation of function f with sequential matches}) and the denominator is $G(\mathbf{x}^\star,\mathbf{w}^\star)$.            
Then, thanks to Lemma~\ref{lemma: dual distribution problem is upper bound}, we obtain
\[
\Ma^2(\mathbf{x}^\star, \mathbf{w}^\star)\geq (1-1/e)\cdot G(\mathbf{x}^\star, \mathbf{w}^\star) \geq (1-1/e)\cdot F(\mathbf{x}^\star, \mathbf{w}^\star)
\]
Finally, we conclude the proof by noting that 
\[
\sum_{e\in E^1} \beta^1_{e}\cdot w^\star_{e}+\Ma^2(\mathbf{x}^\star, \mathbf{w}^\star)\geq \sum_{e\in E^1} \beta^1_{e}\cdot w^\star_{e}+(1-1/e)\cdot F(\mathbf{x}^\star, \mathbf{w}^\star)\geq (1-1/e)\cdot \OPT'\geq (1-1/e)\cdot \OPT
\]
where $\OPT'$ is the optimal value of \eqref{eq:upper_bound_twoperiods} and, in the last inequality (i.e., $\OPT'\geq\OPT$), we use Lemma~\ref{lemma:DH_relaxation}, which we prove in Appendix~\ref{os:remaining_lemmas_2periods}.   
\end{proof}   

\subsection{Remaining Technical Lemmas}\label{os:remaining_lemmas_2periods}

\begin{proof}[Proof of Lemma~\ref{lemma: properties dual distribution problem individual}.]
We prove separately each point.
\begin{enumerate}
\item Let $(\overline{\theta}^\star,\overline{\gamma}^\star)$ be an optimal solution.
The non-negativity of $\overline{\theta}^\star_\ell$'s results from taking $B=\emptyset$ in the constraint ($f_\ell(\emptyset) =0$). To prove the non-negativity of $\overline{\gamma}^\star_{ \ell'}$ 
we note that in the optimal solution, there is a set $B^\star\subseteq\mathcal{X}_\ell$ such that $\overline{\theta}^\star_\ell = f_\ell(B^\star) - \sum_{\ell'\in B^\star}\overline{\gamma}^\star_{\ell'}$; which is the set that maximizes $f_\ell(B)-\sum_{\ell'\in B}\overline{\gamma}^\star_{ \ell'}$. Denote by $C_\ell = \{\ell'\in \mathcal{X}_\ell: \ \overline{\gamma}^\star_{ \ell'}<0\}$. 
Let us redefine $\overline{\theta}'_\ell =\overline{\theta}^\star_\ell + \sum_{\ell'\in C_\ell}\phi_{\ell', \ell}\cdot \overline{\gamma}^\star_{\ell'}$. Clearly $\overline{\theta}'_\ell \geq0$ since
\begin{align*}
\overline{\theta}'_\ell &=\overline{\theta}^\star_\ell + \sum_{\ell'\in C_\ell}\phi^1_{\ell', \ell}\cdot  \overline{\gamma}^\star_{\ell'} \\
&= f_\ell(B^\star) - \sum_{\ell'\in B^\star}\overline{\gamma}^\star_{\ell'}+ \sum_{\ell'\in C_\ell}\phi^1_{\ell', \ell}\cdot \overline{\gamma}^\star_{\ell'}\\
&\geq f_\ell(C_\ell) - \sum_{\ell'\in C_\ell}\overline{\gamma}^\star_{\ell'}+ \sum_{\ell'\in C_\ell}\overline{\gamma}^\star_{\ell'} \\
&= f_\ell(C_\ell)\geq 0
\end{align*}
where we use that $B^\star$ is the maximizing set and $\phi^1_{\ell', \ell}\cdot \overline{\gamma}^\star_{\ell'}\geq \overline{\gamma}^\star_{ \ell'}$ as $\overline{\gamma}^\star_{ \ell'}<0$ for $\ell'\in C_\ell$ and $\phi^1_{\ell', \ell}\ \leq 1$. Redefine $\overline{\gamma}'_{ \ell'} = \overline{\gamma}^\star_{\ell'}$ for all $\ell'\notin C_\ell$ and zero otherwise. Note that the objective of $\overline{\theta}',\overline{\gamma}'$ is the same than the one of $\overline{\theta}^\star,\overline{\gamma}^\star$.
Finally, the constraints are satisfied because, for any $B\subseteq \mathcal{X}_\ell$,
\begin{align*}
\overline{\theta}'_\ell+\sum_{\ell'\in B}\overline{\gamma}'_{\ell'} 
&=f_\ell(B^\star)-\sum_{\ell'\in B^\star}\overline{\gamma}^\star_{\ell'}+\sum_{\ell'\in C_\ell}\phi^1_{\ell', \ell} \cdot \overline{\gamma}^\star_{\ell'}+\sum_{\ell'\in B\setminus C_\ell}\overline{\gamma}^\star_{\ell'} \\
&\geq f_\ell(B\cup C_\ell)-\sum_{\ell'\in B\cup C_\ell}\overline{\gamma}^\star_{\ell'}+\sum_{\ell'\in C_\ell}\overline{\gamma}^\star_{\ell'}+\sum_{\ell'\in B\setminus C_\ell}\overline{\gamma}^\star_{\ell'}\\
&= f_\ell(B\cup C_\ell)\\
& \geq f_\ell(B),
\end{align*}
where in the first equality we use that $\sum_{\ell'\in B}\overline{\gamma}'_{ \ell'}= \sum_{\ell'\in B\setminus C_\ell}\overline{\gamma}^\star_{ \ell'}$ since $\overline{\gamma}'_{\ell'} = 0$ for $\ell'\in C_\ell$. The first inequality follows by the optimality of $B^\star$ i.e. $f_\ell(B^\star)-\sum_{\ell'\in B^\star}\overline{\gamma}^\star_{\ell'}\geq f_\ell(B\cup C_\ell)-\sum_{\ell'\in B\cup C_\ell}\overline{\gamma}^\star_{\ell'}$ and that $\sum_{\ell'\in C_\ell}\phi^1_{\ell', \ell}\cdot  \overline{\gamma}^\star_{\ell'}\geq \sum_{\ell'\in C_\ell}\overline{\gamma}^\star_{\ell'}$. The last inequality is due to monotonicity of $f_\ell$. 

\item Let $(\overline{\theta}^\star,\overline{\gamma}^\star)$ be an optimal solution. Then, $\overline{\theta}^\star = \max_{B\subseteq\mathcal{X}_\ell}\left\{f_\ell(B)-\sum_{\ell'\in B}\overline{\gamma}^\star_{\ell'}\right\}$. Let $B^\star$ be the corresponding maximizer, i.e., $\overline{\theta} = f_\ell(B^\star) - \sum_{\ell'\in B^\star}\overline{\gamma}^\star_{\ell'}$, and let $S^{\star}\subseteq B^\star$ be an optimal solution in $f_\ell(B^\star)$, i.e., $|S^\star|\leq K_\ell$ and $f_\ell(B^\star) = \sum_{\ell'\in S^\star}\phi^2_{\ell,\ell'}$. Therefore, $\overline{\theta}^\star = \sum_{\ell'\in S^\star}\phi^2_{\ell,\ell'} - \sum_{\ell'\in B^\star}\overline{\gamma}^\star_{\ell'}$. 

To find a contradiction, suppose that $|B^\star| > K_\ell$. Then, we know that $S^\star \subset B^\star$ and, since $\overline{\gamma}^\star_{\ell'}\geq 0$, we can remove terms in the second sum and potentially increase this difference, i.e.,
\[
\sum_{\ell'\in S^\star}\phi^2_{\ell,\ell'} - \sum_{\ell'\in B^\star}\overline{\gamma}^\star_{\ell'} \leq \sum_{\ell'\in S^\star}\phi^2_{\ell,\ell'} - \sum_{\ell'\in S^\star}\overline{\gamma}^\star_{\ell'},
\]
contradicting the optimality of $B^\star$.  

\item Let $(\overline{\theta}^\star,\overline{\gamma}^\star)$ be an optimal solution. By part 2. of this lemma, we know that there exists a subset $B^\star$ such that $|B^\star|\leq K_\ell$ and
\[
\overline{\theta}^\star = \sum_{\ell'\in B^\star}(\phi^2_{\ell,\ell'}-\overline{\gamma}^\star_{\ell'}) \geq \max_{B\subseteq \mathcal{X}_\ell}\lrl{ f_\ell(B) -\sum_{\ell'\in B}\overline{\gamma}^\star_{\ell'}} .
\]
Observe that $\phi^2_{\ell,\ell'}-\overline{\gamma}^\star_{\ell'} \geq 0, \; \forall \ell' \in \mathcal{X}_\ell$. To see this, we argue by contradiction; suppose there exists $\ell' \in \mathcal{X}_\ell$ for which the difference is strictly negative (i.e., $\overline{\gamma}^\star_{\ell'} > \phi^2_{\ell,\ell'}$), then:
\begin{itemize}
\item If $\ell'\in B^\star$, we can remove it from $B^\star$ and easily show that $\overline{\theta}^\star$ does not satisfy the constraint for  $B^\star\setminus\{\ell'\}$, which contradicts the feasibility of $(\overline{\theta}^\star,\overline{\gamma}^\star)$.

\item If $\ell'\in\mathcal{X}_\ell\setminus B^\star$, we can create a new solution by decreasing $\overline{\gamma}^\star_{\ell'}$ until $\phi^2_{\ell,\ell'}$. This solution is still feasible since, for any $B$ that contains $\ell'$, we have
\[
\overline{\theta}^\star + \sum_{\ell''\in B\setminus \ell'}\overline{\gamma}^\star_{\ell''}  + \overline{\gamma}^\star_{\ell'} =  \overline{\theta}^\star + \sum_{\ell''\in B\setminus \ell'}\overline{\gamma}^\star_{\ell''}  + \phi^2_{\ell,\ell'} \geq \sum_{\ell''\in B\setminus \ell'}\phi^2_{\ell,\ell''}+ \phi^2_{\ell,\ell'} = \sum_{\ell''\in B}\phi^2_{\ell,\ell''}.
\]
The feasibility for any $B$ that does not contain $\ell'$ also holds. More importantly, this solution has a lower objective value, which contradicts the optimality of $(\overline{\theta}^\star,\overline{\gamma}^\star)$. 
\end{itemize}

Finally, note that due to non-negativity we have $\phi^2_{\ell,\ell''}-\overline{\gamma}^\star_{\ell''}\leq \min_{\ell'\in B^\star}\{\phi^2_{\ell,\ell'}-\overline{\gamma}^\star_{\ell'}\}, \; \forall \ell' \in \mathcal{X}\setminus B^\star$, otherwise we could swap the corresponding terms and contradict the optimality of $B^\star$.\\

Given the properties above, without loss of generality, suppose that the indexes $\ell' \in \mathcal{X}_\ell$ are sorted in decreasing order of $\phi^2_{\ell,\ell'}-\overline{\gamma}^\star_{\ell'}$, i.e.,
\[\phi^2_{\ell,1}-\overline{\gamma}^\star_{1}\geq \phi^2_{\ell,2}-\overline{\gamma}^\star_{2} \geq\cdots\geq \phi^2_{\ell,\lra{\mathcal{X}_\ell}}-\overline{\gamma}_{\lra{\mathcal{X}_\ell}}\geq0.\]
By its optimality, we know that $B^\star$ consists of the first $K_\ell$ elements in this ordering, i.e., $B^\star = \lrl{1, \ldots, K_\ell}$. 

We now show that $\phi^2_{\ell,\ell'}-\overline{\gamma}^\star_{\ell'} = \phi^2_{\ell,\ell''}-\overline{\gamma}^\star_{\ell''} = \alpha^\star$ for any $\ell', \ell'' \in B^\star$. To find a contradiction, suppose that this does not hold. Then, there exists $\ell' \in \lrl{1, \ldots, K_\ell - 1}$ such that $\phi^2_{\ell,\ell'}-\overline{\gamma}^\star_{\ell'} > \phi^2_{\ell,\ell'+1}-\overline{\gamma}^\star_{\ell'+1}$. Let $\ell'$ be the smallest index that this happens. Let $\lrp{\overline{\theta}', \overline{\gamma}'}$ be such that $\overline{\theta}' = \overline{\theta}^\star -\ell'\cdot\epsilon$, $\overline{\gamma}_{\ell''}' = \overline{\gamma}^\star_{\ell''}+\epsilon$ for all $\ell''\in\{1,\ldots,\ell'\}$, and $\overline{\gamma}_{\ell''}' = \overline{\gamma}^\star_{\ell''}$ for all $\ell'' \in\{\ell'+1,\ldots,\lra{\mathcal{X}_\ell}\}$, where $\epsilon > 0$ is such that $\phi^2_{\ell,\ell'}-\overline{\gamma}^\star_{\ell'}-\epsilon = \phi^2_{\ell,\ell'+1}-\overline{\gamma}^\star_{\ell'+1}$. In words, this new solution shifts all the terms in $\{1,\ldots,\ell'\}$ to make them equal to the $(\ell'+1)$-th term. 

Clearly, by feasibility of $\overline{\gamma}^\star$ and the fact that $\epsilon > 0$, we know that $\overline{\gamma}'$ is non-negative, while this also holds for  $\overline{\theta}'$ since 
\begin{align*}
\overline{\theta}' &= \overline{\theta}^\star - \ell'\cdot\epsilon \\
&\geq \sum_{\ell''=1}^{\ell'}(\phi^2_{\ell,\ell''}-\overline{\gamma}^\star_{\ell''}) - \ell'\cdot\epsilon \\
&= \ell'\cdot (\phi^2_{\ell,\ell'+1}-\overline{\gamma}^\star_{\ell'+1})  \\
&\geq 0,
\end{align*}
where the first inequality is due to the feasibility of $\overline{\theta}^\star$ and the fact that $\phi^2_{\ell,\ell''}-\overline{\gamma}^\star_{\ell''} \geq 0, \; \forall \ell'' \in \mathcal{X}_\ell$, while the last equality follows from the definition of $\epsilon$.                   

Also, note that the constraints for $B\subseteq \mathcal{X}_\ell$ such that $\{1,\ldots,\ell'\}\subseteq B$ still hold since
\begin{align*}
\overline{\theta}'&=\overline{\theta}^\star-\ell'\cdot \epsilon \\
&= \sum_{\ell'' \in B^\star} (\phi^2_{\ell,\ell''}-\overline{\gamma}^\star_{\ell''}) - \ell'\cdot\epsilon \\
&= \sum_{\ell'' \in B^\star\setminus\{1,\ldots,\ell'\}} (\phi^2_{\ell,\ell''}-\overline{\gamma}^\star_{\ell''}) + \sum_{\ell''=1}^{\ell'}(\phi^2_{\ell,\ell''}-\overline{\gamma}^\star_{\ell''}) - \ell'\cdot\epsilon \\
&= \sum_{\ell'' \in B^\star\setminus\{1,\ldots,\ell'\}} (\phi^2_{\ell,\ell''}-\overline{\gamma}'_{\ell''}) + \sum_{\ell''=1}^{\ell'}(\phi^2_{\ell,\ell''}-\overline{\gamma}'_{\ell''}) \\  
&\geq \sum_{\ell'' \in B\setminus\{1,\ldots,\ell'\}} (\phi^2_{\ell,\ell''}-\overline{\gamma}'_{\ell''}) + \sum_{\ell''=1}^{\ell'}(\phi^2_{\ell,\ell''}-\overline{\gamma}'_{\ell''}),                                               
\end{align*}
where the inequality is because $B^\star$ is composed by top elements. The argument for $B\subseteq \mathcal{X}_\ell$ such that $\{1,\ldots,\ell'\}\cap(\mathcal{X}_\ell\setminus B) \neq \emptyset$ is analogous. Hence, $\lrp{\overline{\theta}', \overline{\gamma}'}$ is a feasible solution of Problem~\eqref{eq:P2prime_dual_peruser}. Finally, note that this new feasible solution leads to the following objective
\begin{align*}
\overline{\theta}' + \sum_{\ell''=1}^{\ell'}\phi^1_{\ell'',\ell}\overline{\gamma}_{\ell''}' + \sum_{\ell''\in\mathcal{X}_\ell\setminus\{1,\ldots,\ell'\}}\phi^1_{\ell'',\ell}\overline{\gamma}'_{\ell''} 
&= \overline{\theta}^\star-\ell'\cdot\epsilon + \sum_{\ell''=1}^{\ell'}\phi^1_{\ell'',\ell}(\overline{\gamma}_{\ell''}^\star+\epsilon) + \sum_{\ell''\in\mathcal{X}_\ell\setminus\{1,\ldots,\ell'\}}\phi^1_{\ell'',\ell}\overline{\gamma}^\star_{\ell''} \\
&= \overline{\theta}^\star + \sum_{\ell''\in\mathcal{X}_\ell}\phi^1_{\ell'',\ell}\overline{\gamma}^\star_{\ell''} - \epsilon\sum_{\ell''=1}^{\ell'}(1-\phi^1_{\ell',\ell}) \\
&< \overline{\theta}^\star + \sum_{\ell''\in\mathcal{X}_\ell}\phi^1_{\ell'',\ell}\overline{\gamma}^\star_{\ell''}, 
\end{align*}
which contradicts the optimality of $(\overline{\theta}^\star,\overline{\gamma}^\star)$. Therefore, we must have $\phi^2_{\ell,\ell'}-\overline{\gamma}_{\ell'}  = \phi^2_{\ell,\ell'+1}-\overline{\gamma}_{\ell'+1}$ for all $\ell' \in \lrl{1, \ldots, K_\ell }$ and, thus, there exists $\alpha^\star\geq 0$ such that $\phi^2_{\ell,\ell'}-\overline{\gamma}^\star_{\ell'} = \alpha^\star$ for all $\ell'\in B^\star$ and $\phi^2_{\ell,\ell'}-\overline{\gamma}^\star_{\ell'} \leq \alpha^\star$ for all $\ell'\notin B^\star$ (by optimality of $B^\star$). Finally,
\(
\overline{\theta}^\star = \sum_{\ell'\in B^\star} \phi^2_{\ell,\ell'}-\overline{\gamma}^\star_{\ell'} = |B^\star|\cdot\alpha^\star = K_\ell\cdot\alpha^\star.
\). Note that the same argument applies if $ \phi^2_{\ell,\ell'}-\overline{\gamma}^\star_{\ell'} = 0$, in which case $\alpha^\star = 0$. In contrast, if $\alpha^*>0$, the argument does not necessarily apply beyond the first $K_\ell$ elements of the sequence (although we do not care about those elements). For example, suppose that the first $K_\ell+1$ terms are such that $\phi^2_{\ell,\ell'}-\overline{\gamma}^\star_{\ell'} = \alpha$ and strictly greater than the $(K_\ell+2)$-th term. Then, to make them all equal, we need to reduce $\overline{\theta}^\star=K_\ell\cdot\alpha$ by $(K_\ell+1)\cdot\epsilon$ and, thus,  the resulting $\overline{\theta}' = \overline{\theta}^\star - (K_\ell+1)\cdot\epsilon$ may be negative.
\end{enumerate}
\end{proof}

\begin{proof}[Proof of Lemma~\ref{lemma: dual distribution problem is upper bound}.]
Note that it is enough to show that $F_\ell(\mathbf{x}^\star, \mathbf{w}^\star) \leq G_\ell(\mathbf{x}^\star, \mathbf{w}^\star)$ for all $\ell \in \cIJ$. Given a solution $\lrp{\mathbf{x}^\star, \mathbf{w}^\star}$, there are two possible cases.
\begin{enumerate}
\item If $|\mathcal{X}_\ell| \leq K_\ell$, we show that $F_\ell(\mathbf{x}^\star, \mathbf{w}^\star) = G_\ell(\mathbf{x}^\star, \mathbf{w}^\star)$. On the one hand, 
observe that the optimal solution in \eqref{eq:P1prime_peruser} is such that $y_{\ell'} = \phi^1_{\ell',\ell}$ for all $\ell'\in\mathcal{X}_\ell$ because $\sum_{\ell'\in\mathcal{X}_\ell}y_{\ell,\ell'} = \sum_{\ell'\in\mathcal{X}_\ell}\phi^1_{\ell',\ell}\leq |\mathcal{X}_\ell| \leq K_\ell$. Therefore, the optimal objective value in \eqref{eq:P1prime_peruser} is $F_\ell(\mathbf{x}^\star, \mathbf{w}^\star) = \sum_{\ell'\in\mathcal{X}_\ell}\phi^2_{\ell,\ell'}\phi^1_{\ell',\ell}$. 

On the other hand, note that $ |\mathcal{X}_\ell| \leq K_\ell$ implies that $f_\ell (B) = \sum_{\ell'\in \mathcal{X}_\ell}\phi^2_{\ell,\ell'}\ind{\ell'\in B}$ for all $B\subseteq\mathcal{X}_\ell$ since the backlog is smaller than the budget and, thus, the optimal decision in the second period is to show the entire backlog. Therefore, the optimal objective value in \eqref{eq:P2prime_peruser} is 
\[
G_\ell(\mathbf{x}^\star, \mathbf{w}^\star) = \sum_{B\subseteq\mathcal{X}_\ell}\lambda_{B}\cdot f_\ell(B) =  \sum_{B\subseteq\mathcal{X}_\ell}\lambda_{B}\cdot \sum_{\ell'\in \mathcal{X}_\ell}\phi^2_{\ell,\ell'}\ind{\ell'\in B} = \sum_{\ell'\in\mathcal{X}_\ell}\phi^2_{\ell,\ell'}\sum_{B: B\ni \ell'}\lambda_{B} = \sum_{\ell'\in\mathcal{X}_\ell}\phi^2_{\ell,\ell'}\phi^1_{\ell',\ell},
\]
where in the last equality we use the feasibility of $\bm\lambda$. Thus, we conclude that $F_\ell(\mathbf{x}^\star, \mathbf{w}^\star) = G_\ell(\mathbf{x}^\star, \mathbf{w}^\star)$ in this case.

\item If $|\mathcal{X}_\ell| > K_\ell$, we show that $F_\ell(\mathbf{x}^\star, \mathbf{w}^\star) \leq G_\ell(\mathbf{x}^\star, \mathbf{w}^\star)$.
To see this, first note that we can combine Lemma~\ref{lemma: properties dual distribution problem individual} (parts 1. and 2.) and the fact that $f_\ell(B) = \sum_{\ell'\in B}\phi^2_{\ell,\ell'} $ for any $|B|\leq K_\ell$ to rewrite the dual of $G_\ell(\mathbf{x}^\star, \mathbf{w}^\star)$ as
\begin{align}\label{eq:P2prime_dual_peruser}
G^D_\ell(\mathbf{x}^\star, \mathbf{w}^\star) := \min & \quad \overline{\theta} + \sum_{\ell' \in \mathcal{X}_\ell} \phi^1_{\ell',\ell}\cdot \overline{\gamma}_{\ell'}\\
s.t.&\quad \overline{\theta} + \sum_{\ell' \in B}\overline{\gamma}_{\ell'} \geq \sum_{\ell'\in B}\phi^2_{\ell,\ell'}  \qquad  \forall B \subseteq \mathcal{X}_\ell, \ |B|\leq K_\ell\notag\\ 
&\quad \overline{\theta}, \; \overline{\gamma}_{\ell'}\geq 0, \hspace{3cm}\forall  \ell'\in \mathcal{X}_\ell.\notag
\end{align}

Let $\lrp{\bar{\theta}^\star, \bar{\gamma}^\star}$ be an optimal solution of~\eqref{eq:P2prime_dual_peruser} and $\alpha^\star$ be the corresponding constant as defined in Lemma~\ref{lemma: properties dual distribution problem individual} (part 3.). We can construct a feasible solution $(\theta,\gamma)$ for $F_\ell^D$ in \eqref{eq:P1prime_dual_peruser} by taking $\theta = \alpha^\star$ and $\gamma_{\ell'}  = \overline{\gamma}_{\ell'}^\star$ for all $\ell'\in\mathcal{X}_\ell$. Clearly, these values are non-negative. Moreover, for each $\ell'\in\mathcal{X}_\ell$,
\[
\theta + \gamma_{\ell'}  = \alpha^\star + \overline{\gamma}_{\ell'}^\star \geq \phi^2_{\ell,\ell'}-\overline{\gamma}_{\ell'}^\star + \overline{\gamma}_{\ell'}^\star = \phi^2_{\ell,\ell'}.
\]
Finally, the objective value of $(\theta,\gamma)$ is $K_\ell\alpha^\star + \sum_{\ell'}\phi^1_{\ell',\ell}\overline{\gamma}_{\ell'}^\star$, which is equal to objective value of $(\overline{\theta}^\star,\overline{\gamma}^\star)$ in~\eqref{eq:P2prime_dual_peruser}. Hence, we know that
\[
G_\ell(\mathbf{x}^\star, \mathbf{w}^\star) = G^D_\ell(\mathbf{x}^\star,\mathbf{w}^\star) = K_\ell\theta + \sum_{\ell' \in \mathcal{X}_\ell}\phi^1_{\ell',\ell}{\gamma}_{\ell'} \geq F^D_\ell(\mathbf{x}^\star,\mathbf{w}^\star) = F_\ell(\mathbf{x}^\star, \mathbf{w}^\star), 
\]
where in the first and last equalities we use strong duality.\footnote{Both duals are always feasible, since we can consider $\theta=\overline{\theta} = 0$ and $\gamma_{\ell'}=\overline{\gamma}_{\ell'} = \phi^2_{\ell,\ell'}$ for all $\ell'\in \mathcal{X}_\ell$.}
Hence, we conclude that $F_\ell(\mathbf{x}^\star, \mathbf{w}^\star) \leq G_\ell(\mathbf{x}^\star, \mathbf{w}^\star)$ in this case.
\end{enumerate}
\end{proof}

\begin{lemma}\label{lemma:DH_relaxation}
Problem~\eqref{eq:upper_bound_twoperiods} is an upper bound of Problem~\ref{def:two_period_problem}.
\end{lemma}  
\begin{proof}[Proof of Lemma~\ref{lemma:DH_relaxation}.]
Let $(\mathbf{x},\mathbf{w})$ be a feasible solution in Problem~\ref{def:two_period_problem}. Define $\hat{\mathbf{x}}^1= \mathbf{x}$, $\hat{\mathbf{w}}^1= \mathbf{w}$ and $\hat{\mathbf{x}}^2$ as follows: for all $\ell\in I\cup J$, $\ell'\in \potentials^1_\ell$
\(
\hat{x}^2_{\ell',\ell} = \sum_{R\subseteq\vec{E}^1} \ind{\ell\in S_{\ell'}(R)}\cdot \probP_{{\hat{\mathbf{x}}^1}}(\cR=R),
\)
where $S_{\ell'}(R)$ is an optimal solution for user $\ell'$ in the problem defined in $f(R)$. Clearly, $\hat{\mathbf{x}}^1, \hat{\mathbf{w}}^1$ satisfy their constraints in Problem~\eqref{eq:upper_bound_twoperiods}, so it remains to show the constraints that involve $\hat{\mathbf{x}}^2$ are also satisfied. First, note that
\[
\begin{split}
      \hat{x}^2_{\ell',\ell} &= \sum_{R\subseteq\vec{E}^1} \ind{\ell\in S_{\ell'}(R)}\cdot \probP_{\hat{\mathbf{x}}^1}(\cR=R) \\
      &= \sum_{R'\subseteq\vec{E}^1\setminus(\ell,\ell')} \ind{\ell\in S_{\ell'}(R'\cup \{(\ell,\ell')\})}\cdot \probP_{\hat{\mathbf{x}}^1}(\cR =R'\cup (\ell,\ell')) \\
      &= \sum_{R'\subseteq \vec{E}^1\setminus (\ell,\ell')} \ind{\ell\in S_{\ell'}(R'\cup (\ell,\ell'))}\cdot \probP_{\hat{\mathbf{x}}^1}(\cR = (\ell,\ell')) \cdot\probP_{\hat{\mathbf{x}}^1}(\cR=R') \\
      &= \hat{x}^1_{\ell,\ell'}\cdot \plike_{\ell,\ell'}\cdot \sum_{R'\subseteq \vec{E}^1\setminus (\ell,\ell')} \ind{\ell\in S_{\ell'}(R'\cup (\ell,\ell'))} \cdot\probP_{\hat{\mathbf{x}}^1}(\cR =R') \leq \hat{x}^1_{\ell,\ell'}\cdot \plike_{\ell,\ell'}^1 
\end{split}
  \]
where the third equality follows by the independent choices made by the users, the fourth uses that $\probP_{\hat{\mathbf{x}}}(\cR = (\ell,\ell')) =  \hat{x}^1_{\ell,\ell'}\cdot \plike_{\ell,\ell'}^1$ and the last inequality follows because the sum is at most 1. We now focus on the cardinality constraints of $\hat{\mathbf{x}}^2$:
\[
\begin{split}
\sum_{\ell\in \potentials^1_{\ell'}} \hat{x}^2_{\ell',\ell} &= \sum_{\ell\in \potentials^1_{\ell'}} \sum_{R\subseteq \vec{E}^1} \ind{\ell\in S_{\ell'}(R)}\cdot \probP_{\hat{\mathbf{x}}^1}(\cR = R) \\
&= \sum_{R\subseteq \vec{E}^1} \sum_{\ell\in \potentials^1_{\ell'}} \ind{\ell\in S_{\ell'}(R)}\cdot \probP_{\hat{\mathbf{x}}^1}(\cR = R) \\
&= \sum_{R\subseteq \vec{E}^1} \lrc{\sum_{k=0}^{K_{\ell'}-1} k\cdot \probP_{\hat{\mathbf{x}}^1}\lrp{\lra{B_{\ell'}} = k } + K_{\ell'}\cdot  \probP_{\hat{\mathbf{x}}^1}\lrp{\lra{B_{\ell'}} \geq K_{\ell'} }  }\cdot \probP_{\hat{\mathbf{x}}^1}(\cR = R) \\
&\leq \sum_{R\subseteq \vec{E}^1} K_{\ell'} \cdot \probP_{\hat{\mathbf{x}}^1}(\cR = R) \\
&= K_{\ell'}, 
\end{split}
\]
where $B_\ell = \{\ell'\in\potentials^1_\ell:\; (\ell',\ell)\in R\}$. 
\end{proof}


%% file: 10_electronic_companion.tex
%
%

\section{Analysis of Suboptimal Methods}\label{app:suboptimal_methods}
A natural approach to solve Problem~\ref{def:two_period_problem} is to adapt commonly used algorithms in the online matching and assortment optimization literature. One such algorithm is the greedy policy, which provides a performance guarantee of $1/2$ for the online matching problem~\citep{kvv}.
In our setting, such a (local) greedy policy would select, for each user and period, the subset of profiles that maximizes their expected number of matches, i.e.,
\[
S_\ell^t = \argmaxA_{S\subseteq \potentials_\ell^t: \lra{S}\leq K_\ell} \lrl{\sum_{\ell'\in S} \plike_{\ell, \ell'}\cdot \lrp{\ind{\ell'\in B_\ell^t} + \plike_{\ell' \ell}\cdot \ind{\ell'\notin B_\ell^t}} }.
\]
In Proposition~\ref{prop: worst case performance greedy}, we show that this policy achieves a worst-case performance arbitrarily close to zero, as it does not account for the potential ``congestion'' that some users may cause on others.
\begin{proposition}\label{prop: worst case performance greedy}
The worst-case approximation guarantee for the local greedy policy is $O(1/n)$, where $n$ is the size of the market.
\end{proposition}
\begin{proof}[Proof.]
Suppose that there are \(n\) users on each
side of the market, i.e., \(I = \lrl{i_1, \ldots, i_n}\) and \(J = \lrl{j_1, \ldots, j_n}\).
In addition, suppose that \(\potentials^1_i = J\) for every $i \in I$,  $\potentials^1_j = I$ for every $j \in J$ and $K_\ell = 1$ for all $\ell\in I\cup J$. Let us set the probabilities: \(\beta^1_{i,j} = 1\) for \(j=j_1\) and for all \(i \in I\),
\(\beta^1_{i,j} = 1-\varepsilon\) for all \(i\in I\) and \(j\neq j_1\), while \(\beta^2_{i,j}=0\) for all \(i\in I\cup J, \, j\in \potentials^1_i\).
In this setting, the Local Greedy policy will choose \(S^1_i = \lrl{j_1}\) for every user \(i\),
and therefore only one match will take place in expectation. In contrast, an optimal solution is
to assign \(S^1_{i_k} = \lrl{j_k}\), which leads to \(1+(n-1) (1-\varepsilon)\) matches
in expectation. Then, the performance of the greedy policy is given by
$1/(1+(n-1)(1-\varepsilon)) \to 0$,
as \(n\rightarrow \infty\) for \(\varepsilon\) sufficiently small.
\end{proof}
An alternative approach is to find a maximum weight perfect matching in each period, where the weight of each edge is the probability of having a match between the users. Formally, the Perfect Matching heuristic solves, in each period $t\in [T]$, the following problem:
\begin{equation}\label{eq: perfect matching}
\begin{split}
\max \quad & \; \sum_{\ell \in I\cup J} \sum_{\ell' \in \potentials_\ell^1} y_{\ell, \ell'}^t \plike^t_{\ell, \ell'} + \frac{1}{2}  w_{\ell, \ell'}^t \beta^t_{\ell, \ell'} \\
st.   \quad & y_{\ell, \ell'}^t \leq \ind{\ell'\in \backlog_\ell^t}, \hspace{9.7em} \forall \ell \in \cIJ , \; \ell' \in \potentials^t_\ell\\
&  x_{\ell, \ell'}^t + y_{\ell, \ell'}^t \leq 1, \hspace{9.8em} \forall \ell \in \cIJ, \; \ell' \in \potentials^t_\ell  \\
&  \sum_{\ell' \in \potentials_\ell^t} x_{\ell, \ell'}^t + y_{\ell, \ell'}^t \leq K_\ell, \hspace{6.8em} \forall  \ell \in \cIJ \\
&  w_{\ell, \ell'}^t \leq x_{\ell, \ell'}^t, \; w_{\ell, \ell'}^t \leq x_{\ell' \ell}^t,\; w_{\ell, \ell'}^t = w_{\ell' \ell}^t, \hspace{0.7em} \forall \ell \in \cIJ, \; \ell' \in \potentials^t_\ell \\
&  x_{\ell, \ell'}^t, y_{\ell, \ell'}^t, w_{\ell, \ell'}^t \in \lrl{0,1}, \hspace{6.2em} \forall \ell \in \cIJ, \; \ell' \in \potentials^t_\ell.
\end{split}
\end{equation}
Then, the method sets \(S_\ell^t = \lrl{\ell' \in \potentials_\ell^t: x_{\ell, \ell'}^t = 1 \text{ or } y_{\ell, \ell'}^t = 1 }\) for each \(\ell\in \cIJ\). Note that, if there is no initial backlog, then the problem can be re-formulated as:
\begin{equation}\label{eq: perfect matching without initial backlog}
\begin{split}
\max \quad & \; \sum_{e\in E} \beta^t_{e} \cdot w_{e}^t \\
st.   \quad &  \sum_{e\in E: \ell\in e} w_{e}^t \leq K_\ell, \hspace{5.5em} \forall \ell \in \cIJ \\
&  \sum_{t\in [T]} w_{e}^t \leq 1, \hspace{7.4em} \forall e\in E \\
&  w_{e}^t \in \lrl{0,1}, \hspace{7.5em} \forall e\in E,\; t\in [T],
\end{split}
\end{equation}
where the second constraint ensures that each edge is used at most once, i.e., that no two users see each other more than once. This is captured in the previous formulation through the set of potentials {~\cite{Chen09} and~\cite{Jeloudar21} consider a similar approach in the matching probing problem with and without commitment and show that it achieves a performance guarantee of 1/4 and 0.43, respectively.}
Nevertheless, as we show in Proposition~\ref{prop:sequential_worst}, this policy has also a worst-case performance arbitrarily close to zero. This result is not surprising, as the perfect matching policy does not exploit the information provided by the realized like decisions.
\begin{proposition}\label{prop:sequential_worst}
The worst-case approximation guarantee for the perfect matching policy is $O(1/n)$, where $n$ is the size of the market		
\end{proposition}
\begin{proof}[Proof.]
Suppose that \(\lra{I} = 2n\), \(\lra{J} = 2\), that $\potentials^1_i = J$ for every $i \in I$, $\potentials^1_j = I$ for every $j \in J$, $K_\ell = 1$ for all $\ell\in I \cup J$ and that \(\plike^t_{i,j} = p\) while \(\plike^t_{j,i} = q\) for all \(i\in I, \; j\in J\), and $t\in \{1,2\}$.
Then, it is easy to see that the sequential perfect match policy leads to \(4pq\) matches in expectation. On the other hand, consider the policy where: 
\begin{enumerate}
\item[(i)] In \(t=1\), \(\lrl{i_1, \ldots, i_n}\) see \(j_1\), \(\lrl{i_{n+1}, \ldots, i_{2n}}\) see \(j_2\), \(j_1\) sees \(i_{2n}\) and \(j_2\) sees \(i_1\).
\item[(ii)] In \(t=2\), \(\lrl{i_1, \ldots, i_n}\) see \(j_2\), \(\lrl{i_{n+1}, \ldots, i_{2n}}\) see \(j_1\), \(j_1\) sees any profile that liked her in \(t=1\), and same for \(j=2\).
\end{enumerate}
Given this policy, the matches \((i_{2n}, j_1)\) and \((i_1, j_2)\) happen with probability \(pq\) each. On the other hand, \(j_1\) matches with someone in \(\lrl{i_1, \ldots, i_n}\) with probability \((1-(1-p)^n)q\), and the same for \(j_2\) matching with someone in \(\lrl{i_{n+1}, \ldots, i_{2n}}\). Then, the total expected number of matches is \(2pq + 2q(1-(1-p)^n)\), which is optimal for this instance. Then, the sequential perfect match policy achieves a performance of
$4pq/(2q(p + 1 - (1-p)^n)) \to 2/(1+n)$ when $p\rightarrow 0$, and since $2/(1+n)\to 0$ when $n\to \infty$, we conclude the proof. 
\end{proof}

Given that our model requires the realization of two random variables to generate a match, the approaches above perform poorly because of their \emph{non-adaptive} nature. In the following section, we overcome this challenge by using submodular maximization techniques.

\section{Analysis of Polynomial-Time Approximation Algorithms}\label{ec:design_dependent_analysis}
As discussed in Section~\ref{sec: design_dependent_analysis}, we now present polynomial-time approximation algorithms based on submodular techniques that provide weaker constant-factor approximation guarantees but that ensure scalability in large markets. Our main result, formalized in Theorem~\ref{thm: submodular guarantees}, shows that the specific guarantees depend on the combination of design choices, i.e., whether matches are simultaneous or sequential and whether the design is one or two directional.


At the core of the proof is Algorithm~\ref{alg:greedy}, which can be adapted to the different platform designs to show the guarantees. In Sections~\ref{subsec: one-directional sequential} and~\ref{subsec: two-directional sequential}, we consider settings where the platform only allows sequential matches, and we vary whether interactions are one or two-directional, respectively. Then, in Sections~\ref{subsec: one-directional simultaneous} and~\ref{subsec: two-directional simultaneous}, we relax the sequential-only restriction and allow for non-sequential matches in the first period. Recall that all of these settings do not consider non-sequential matches in the second period, however, in Section~\ref{subsec: simultaneous matches in second period}, we study a setting with non-sequential matches in the second period under the small probability assumption.
\begin{algorithm}[h]
\caption{Policy based on submodular maximization}\label{alg:greedy}
\begin{algorithmic}[1]
\Require  An instance of the problem and an algorithm $\mathrm{ALG}$ for submodular maximization subject to a design-dependent feasible region.
\Ensure Feasible subsets: $\mathbf{x}^1$, $\mathbf{x}^2$
\State Use {$\mathrm{ALG}$} to obtain an approximate solution $\mathbf{x}^1$ of Problem~\ref{def:two_period_problem}.
\State Update potentials and backlogs.
\State Obtain $\mathbf{x}^2$ by solving to optimality \eqref{eq:second_stage_problem}.
\end{algorithmic}
\end{algorithm}

\subsection{One-directional Interactions and Sequential Matches.}\label{subsec: one-directional sequential}
This setting is similar to that studied in the recent two-sided assortment optimization literature~\citep{ashlagi19,torrico21,aouad21}, where customers first select a supplier from their assortment and, later, the suppliers observe all the customers that chose them and decide whom to serve. Our model departs from this literature in two key aspects. First, we assume that users can like as many profiles as they want and, consequently, they can potentially match with multiple users on the other side. Second, we assume that like probabilities are independent of the subset of profiles displayed, while these papers consider an underlying choice model (e.g., MNL) to compute them. As mentioned above, this independence assumption greatly simplifies the analysis, and it is practical given that assortments minimally affect users' like probabilities~\citep{rios2021}.

Without loss of generality, we assume that interactions can only be initiated by agents in $I$. Hence, in the first period, the platform must choose a subset of profiles for each user in $I$, i.e., $\mathbf{x}^1 = \lrl{x_{i,j}^1}_{i\in I,\; j\in \potentials^1_i}$, that satisfies the constraints defined by the platform, namely, that $\sum_{j\in \potentials^1_i} x_{i,j}^1 \leq \asize_i$ for each $i\in I$. Hence, the feasible region for the first-period decisions can be formulated as:
\begin{align}\label{eq:feasibleregion_onesided_nosimul}
\Bigg\{\mathbf{x}^1\in\{0,1\}^{\vec{E}^1_I}:  &\sum_{j\in\potentials^1_i}x^1_{i,j} \leq \asize_i,\; \forall i\in I\Bigg\}.
\end{align}
Note that the only difference between \eqref{eq:feasibleregion_onesided_nosimul} and the feasible region in Problem~\ref{def:two_period_problem} is that the former enforces $\mathbf{w}^1 = 0$ and $\vec{E}^1_I$ instead of $\vec{E}^1 = \vec{E}^1_I \cup \vec{E}^1_J$, ensuring that only users in $I$ can see profiles in the first period. 
It is easy to see that the feasible region \eqref{eq:feasibleregion_onesided_nosimul} is a partition matroid. Therefore, we can show a $(1-1/e)$-approximation guarantee by combining the continuous greedy algorithm introduced by~\cite{vondrak08} (for submodular maximization under matroid constraints) with the dependent randomized rounding algorithm by~\cite{gandhi06}. We formalize this in Proposition~\ref{prop:approximation_one_directional_alg1}.
\begin{proposition}\label{prop:approximation_one_directional_alg1}
When Problem~\ref{def:two_period_problem} 
 is restricted to one-directional policies with sequential matches, there exists a feasible solution $\mathbf{x}^1$ whose objective value is at least $(1-1/e)\cdot\OPT$.
\end{proposition}
Before showing this proposition, let us recall the definition of the multilinear extension of a set function.
 \begin{definition}[Multilinear Extension]\label{def:multext_submodular}
For a generic non-negative set function $g:\{0,1\}^{\mathcal{U}}\to \RR_+$ defined over a space of elements $\mathcal{U}$, we define its {\it multilinear extension} $G:[0,1]^{\mathcal{U}}\to\RR_+$ by
\(
G(\mathbf{x})=\sum_{S\subseteq \mathcal{U}}g(S)\prod_{e\in S}x_{e}\prod_{e\in \mathcal{U}\setminus S}(1-x_{e}).
\)
\end{definition}
\begin{proof}[Proof of Proposition~\ref{prop:approximation_one_directional_alg1}.]
To formalize our analysis, let $F$ be the multilinear extension of the set function $f$ defined in \eqref{eq:second_stage_problem}, i.e.,
\[
F(\mathbf{x}) = \sum_{R\subseteq \vec{E}^1_I} f(R)\cdot \prod_{e\in R} x_e \cdot\prod_{e\notin R} (1-x_e).
\]
Note that we define this extension over arcs in $\vec{E}^1_I$ since this is a one-directional setting. Also, observe that for any $\mathbf{x}\in\{0,1\}^{\vec{E}^1_I}$, we have $F(\mathbf{\phi}^1\cdot \mathbf{x}^1 ) = \Ma^2(\mathbf{x}^1, \mathbf{w}^1)$,
where $\mathbf{\phi}^1\cdot\mathbf{x}^1 $ denotes the vector with components $\phi^1_{i,j}\cdot x^1_{i,j}$ for all $i\in I, \ j\in J$. As previously discussed, representing $\Ma^2(\mathbf{x}^1, \mathbf{w}^1)$ through the multilinear extension of $f$ has the advantage that the latter can be evaluated in $[0,1]^{\vec{E}^1_I}$ rather than only in $\{0,1\}^{\vec{E}^1_I}$.
Then, consider the following optimization problem:
\begin{subequations}
\begin{align}
\max &  \quad F(\mathbf{z})  \label{eq:objective}\\
\text{s.t.}&  \;\sum_{j\in J: \phi^1_{i,j}>0} \frac{z_{i,j}}{\phi^1_{i,j}} \le K_i \ \quad\text{ for every } i\in I, \label{eq:budget} \\
&\quad 0\le z_{e} \le \phi^1_{e}\hspace{3.2em}\text{ for every } e \in\vec{E}^1_I. \label{eq:box}
\end{align}
\end{subequations}
First, we show the following:
\begin{lemma}\label{lemma:relaxation_nosimul}
The optimal value of \eqref{eq:objective}-\eqref{eq:box} is an upper bound on the optimal value of Problem~\ref{def:two_period_problem} under one-directional interactions and sequential matches.
\end{lemma}

\begin{proof}[Proof.]
Consider a feasible solution $\mathbf{x}^1\in\{0,1\}^{\vec{E}^1_I}$ of Problem~\ref{def:two_period_problem} under one-directional interactions and sequential matches, that is,  $\sum_{j\in\potentials^1_i}x^1_{i,j} \leq K$ for every $i\in I$.
Let $\mathbf{z}=\phi^1 \cdot \mathbf{x}^1$.
Note that, for each $i\in I$, 
\[\sum_{j\in J:\phi^1_{i,j}>0}z_{i,j}/\phi^1_{i,j}=\sum_{j\in J:\phi^1_{i,j}>0}x^1_{i,j}\le K_i,\]
and $z_e=x^1_e\phi^1_e\le \phi^1_e$ for every $e \in\vec{E}^1_I$. Therefore, $\mathbf{z}$ is a feasible solution for the problem \eqref{eq:objective}-\eqref{eq:box}.
Since the objective value of $\mathbf{x}^1$ in~Problem~\ref{def:two_period_problem} is equal to the objective of $\mathbf{z}$ in \eqref{eq:objective}-\eqref{eq:box}, the proof of the lemma follows. 
\end{proof}

Since $f$ is monotone and submodular and $F$ inherits all its properties, we can use Lemma 4.2 in~\cite{vondrak08} (see Corollary~\ref{theorem:continuous-greedy}) to find our desired performance guarantee. Formally,
\begin{corollary}[\cite{vondrak08}]\label{theorem:continuous-greedy}
There exists an efficient algorithm that computes a point $\mathbf{z}$ that satisfies \eqref{eq:budget} and \eqref{eq:box} such that $F(\mathbf{z})\geq (1-1/e)\cdot F(\mathbf{z^*})$, where $\mathbf{z}^*$ is an optimal solution of \eqref{eq:objective}-\eqref{eq:box}.
\end{corollary}
We emphasize that the solution $\mathbf{z}$ in Corollary~\ref{theorem:continuous-greedy} might be fractional, so we need to use a rounding procedure to construct the final solution of our problem. To find our feasible point we use: (i) the continuous greedy algorithm proposed by~\citep{vondrak08}, and (ii) the dependent randomized rounding algorithm by~\cite{gandhi06}. Specifically, the feasible solution can be obtained as follows: 
\begin{enumerate}
\item Compute a solution $\mathbf{z}$ for the problem \eqref{eq:objective}-\eqref{eq:box} using the algorithm from Corollary \ref{theorem:continuous-greedy}.
\item For each $i\in I$ and $j\in J$, set $ \tilde{x}_{i,j}=  z_{i,j}/\phi^1_{i,j}$ when $\phi^1_{i,j}>0$ and zero otherwise.
\item Independently for each user $i\in I$, run the dependent randomized rounding algorithm \citep{gandhi06} on the fractional vector $\mathbf{\tilde{x}}_i\in\RR^J$ to compute an integral random vector $\mathbf{x}_i\in \{0,1\}^J$.
\end{enumerate}

\begin{algorithm}[h]
\caption{Approximation Algorithm for One-Directional Interactions and Sequential Matches}\label{alg:one}
\begin{algorithmic}[1]
\Require  An instance of the problem.
\Ensure A feasible vector $\mathbf{x}$.
\State Compute a solution $\mathbf{z}$ for the problem \eqref{eq:objective}-\eqref{eq:box} using the algorithm from Corollary \ref{theorem:continuous-greedy}.
\State For each $i\in I$ and $j\in J$, set $ \tilde{x}_{i,j}=  z_{i,j}/\phi^1_{i,j}$ when $\phi^1_{i,j}>0$ and zero otherwise.
\State \multiline{Independently for each user $i\in I$, run the dependent randomized rounding algorithm \citep{gandhi06} on the fractional vector $\tilde{\mathbf{x}}_i$ to compute an integral random vector $\mathbf{x}_i\in \{0,1\}^J$.}
\end{algorithmic}
\end{algorithm}

Given the fractional solution $\mathbf{z}$ satisfying the guarantee in Corollary~\ref{theorem:continuous-greedy}, let $\mathbf{\tilde{x}}_i\in [0,1]^J$ be the fractional vector such that the $j$-th entry is equal to $ z_{i,j}/\phi^1_{i,j}$ when $\phi^1_{i,j}>0$ and zero otherwise.
Observe that thanks to constraint \eqref{eq:box} we have $\mathbf{\tilde{x}}_i\in [0,1]^J$ for each $i\in I$.
Then, independently for each user $i\in I$, by the algorithm in \citep{gandhi06} it is possible to efficiently compute an integral random vector $\mathbf{x}_i\in \{0,1\}^J$ satisfying the following conditions:
\begin{enumerate}
\item  $\sum_{j\in J}{x}_{i,j}\le \lceil \sum_{j\in J}\tilde{x}_{i,j}\rceil$, and \label{eq:hitting}
\item  $\E[{x}_{i,j}]=\tilde{x}_{i,j}$ for each $i\in I$ and $j\in J$. \label{eq:expected}
\end{enumerate}
Thanks to condition \eqref{eq:hitting} of the randomized rounding algorithm, for each $i\in I$ we have
$
\sum_{j\in J}{x}_{ij}\le \left\lceil \sum_{j\in J}\tilde{x}_{ij}\right\rceil=\lceil \sum_{j\in J:\phi_{ij}>0}  z_{ij}/\phi_{ij}^1\rceil\le K_i,
$
where the last inequality holds since $K_i$ is integral and $\tilde{\mathbf{z}}$ satisfies constraint \eqref{eq:budget}.
Therefore, our algorithm gives a feasible solution for Problem~\ref{def:two_period_problem}.
We now analyze the approximation guarantee.
\begin{align*}
\E_{\mathbf{x}\sim\tilde{\mathbf{x}}}[\Ma^2({\mathbf{x}})] = \E_{\mathbf{x}\sim\tilde{\mathbf{x}}}[\E_{\cR \sim \mathbf{\phi}^1 {\mathbf{x}}}\lrc{f(\cR)}]
&= \E_{\mathbf{x}\sim\tilde{\mathbf{x}}}\Bigg[\sum_{R\subseteq \vec{E}^1_I} f(R)\cdot \probP_{{\mathbf{x}}}\lrp{{\cR} = R}\Bigg]\\
&= \sum_{R\subseteq \vec{E}^1_I}f(R)\cdot \E_{\mathbf{x}\sim\tilde{\mathbf{x}}} \left[\prod_{e\in R}\plike^1_e{x}_e\prod_{e\notin R}(1-\plike^1_e{x}_e)\right] \\
&= \sum_{R\subseteq \vec{E}^1_I}f(R)\cdot \prod_{e\in R}\plike^1_e\cdot\E_{\mathbf{x}\sim\tilde{\mathbf{x}}}[{x}_e]\cdot\prod_{e\notin R}(1-\plike^1_e\cdot \E_{\mathbf{x}\sim\tilde{\mathbf{x}}}[{x}_e]) \\
&=\sum_{R\subseteq \vec{E}^1_I}f(R)\cdot\prod_{e\in R}\plike^1_e\cdot \frac{ z_{e}}{\phi^1_{e}}\prod_{e\notin R}(1-\plike^1_e\cdot \frac{ z_{e}}{\phi^1_{e}}) = F(\mathbf{z}),
\end{align*}
where the second equality comes from the fact that ${\mathbf{x}}_i$ is independent from ${\mathbf{x}}_{i'}$ for every $i,i'\in I$ with $i\ne i'$, and the third equality comes from condition \ref{eq:expected} of the randomized rounding procedure.
Finally, Lemma \ref{lemma:relaxation_nosimul} states that $\OPT'\geq \OPT$, where $\OPT'$ is the optimal value of Problem~ \eqref{eq:objective}-\eqref{eq:box}, so we conclude the proof of Proposition~\ref{prop:approximation_one_directional_alg1} by using Corollary \eqref{theorem:continuous-greedy}. 
\end{proof}

\subsection{Two-directional Interactions and Sequential Matches.}\label{subsec: two-directional sequential}
To capture this platform design, we now assume that the first period decisions, represented by $\mathbf{x}^1 \in \lrl{0,1}^{\vec{E}^1}$, include arcs in both directions, i.e., $\vec{E}^1 = \vec{E}^1_I\cup \vec{E}^1_J$.
Then, the feasible region for the first-period decisions can be characterized as:
\begin{equation}\label{eq:feasibleregion_generalsequential}
\Bigg\{\mathbf{x}^1\in \lrl{0,1}^{\vec{E}^1}:\sum_{ \ell'\in\potentials^1_\ell}x^1_{\ell,\ell'} \leq K_\ell,\; \forall \ell\in I\cup J, \text{ and } x^1_{\ell,\ell'}+x^1_{\ell',\ell}\leq 1, \ \forall  \ell\in I\cup J, \; \ell'\in\potentials^1_\ell	\Bigg\}.
\end{equation}
The first family of constraints ensures our cardinality requirements, while the second ensures that no pair of users see each other in the first period, preventing non-sequential matches. Note that we set $\mathbf{w}^1=0$ since non-sequential matches are not allowed.

Following a similar strategy as for the one-directional case, we can exploit the submodularity of the objective function and the structure of the feasible region to find an approximation guarantee. Specifically, the feasible region in~\eqref{eq:feasibleregion_generalsequential} corresponds to the intersection of two matroids as the following lemma shows.
\begin{lemma}\label{lemma:intersection_2partitions}
The feasible region in~\eqref{eq:feasibleregion_generalsequential} corresponds to the intersection of two matroids.
\end{lemma}
\begin{proof}[Proof.]
Our ground set of elements is $\mathcal{E} = \vec{E}^1_I\cup\vec{E}^1_J$.
The first partition consists of the following parts: $\mathcal{E}_\ell = \{e: e=(\ell,\ell') \ \text{for every} \ \ell'\in\potentials^1_\ell\}$ for all $\ell \in I\cup J$.
It is easy to check that $\mathcal{E} = \cup_{\ell\in I\cup J}\mathcal{E}_\ell$ and $\mathcal{E}_\ell\cap \mathcal{E}_{\ell'}=\emptyset$ for every $\ell,\ell'$ such that $\ell \neq \ell'$. Finally, the budget for each part $\mathcal{E}_\ell$ is $K_\ell$. Now, let us construct the second partition matroid.
For every pair $i\in I$ and $j\in J$, we define a part $\mathcal{E}_{i,j}$ as the set $\{(i,j), (j,i)\}$. Indeed this forms a partition of $\mathcal{E}$. Finally, the budget for each part $\mathcal{E}_{i,j}$ is 1. 
\end{proof}   
Given this result, we can use~\citet{fisher78} to devise a feasible solution with an approximation factor of $1/3$, or we can apply the local search algorithm in~\citep{lee2009submodular} for an improved guarantee of $1/(2+\epsilon)$\footnote{Both methods were originally designed for submodular maximization under the intersection of matroids.}. We formalize this result in the next proposition.
\begin{proposition}\label{prop:approximation_two_directional_alg2}
When Problem~\ref{def:two_period_problem} 
 is restricted to two-directional policies with sequential matches, there exists a feasible solution $\mathbf{x}^1$ whose objective value is at least $\frac{1}{2+\epsilon}\cdot\OPT$, for any $\epsilon>0$.
\end{proposition}
\begin{proof}[Proof.]
Consider Problem~\ref{def:two_period_problem} with two-directional interactions and sequential matches. Since $\mathbf{x}^1$ is a 0-1 vector and $\Ma^2(\cdot)$ is monotone submodular over elements in $\vec{E}^1$, we know that a vanilla greedy algorithm achieves a $1/(1+r)$-approximation for the problem of maximizing a monotone submodular function over the intersection of $r$ matroids \cite{fisher78} and local-search guarantees a factor of $1/(r+\epsilon)$ for any fixed $\epsilon>0$~\citep{lee2009submodular}.
By Lemma~\ref{lemma:intersection_2partitions}, we know that $r=2$, and, thus, the guarantee is $1/3$ for greedy and $1/(2+\epsilon)$ for local-search.
\end{proof} 

\begin{remark}
Note that the approximation factor worsens relative to Proposition~\ref{prop:approximation_one_directional_alg1} because the family of constraints that prevents non-sequential matches correlates the decisions for each pair of users.
An alternative approach would be to remove these constraints and penalize non-sequential matches. However, this approach may affect the submodularity of $f$ in the second period. Another possibility is to consider the approach used in Proposition~\ref{prop:approximation_one_directional_alg1}. However, the correlation between the decision variables prevents from using the dependent randomized rounding for each user.
\end{remark}
\subsection{One-directional Interactions and Non-Sequential Matches.}\label{subsec: one-directional simultaneous}
As in Section~\ref{subsec: one-directional sequential}, we assume (without loss of generality) that interactions are one-directional with $I$ as the initiating side. However, we now allow the platform to select pairs of users that will simultaneously see each other in the first period and, thus, can potentially generate a non-sequential match. To accomplish this, we use the variables $\mathbf{w}^1 \in \lrl{0,1}^{E^1}$ defined over the set of undirected edges between $I$ and $J$, and we also set $x^1_{\ell,\ell'}=0$ for all $\ell\in J$, $\ell'\in\potentials^1_\ell$ since interactions are only initiated by users in $I$, i.e., $\mathbf{x}^1\in \lrl{0,1}^{\vec{E}^1_I}$. Then, the feasible region for the first-period decisions becomes:
\begin{align}\label{eq:feasibleregion_simultaneous_onesided}
    \Bigg\{\mathbf{x}^1\in\{0,1\}^{\vec{E}^1_I}, \ \mathbf{w}^1\in\{0,1\}^{E^1}: &  \ x^1_{i,j} + w^1_{e} \leq 1, \; \forall i\in I, \; j\in \potentials^1_i, \; e=\{i,j\}\notag\\
    &\sum_{j\in \potentials^1_i}x^1_{i,j} +\sum_{e\in E: i\in e} w^1_{e}\leq K_i, \; \forall  i\in I, \notag\\
    &\sum_{e\in E^1: j\in e} w^1_e \leq K_j, \; \forall  j\in J.\Bigg\}
\end{align}
The first family of constraints guarantees that no profile targets both a sequential (i.e., $x_{i,j}^1 = 1$) and a non-sequential match (i.e., $w_{e}^1 = 1$), while the second and third families of constraints ensure that the subsets of profiles to display satisfy the cardinality requirements for sides $I$ and $J$, respectively. Note that each user $j\in J$ only sees profiles involving users $i\in I$ for which $w_{e}^1=1$ and, thus, users in $J$ cannot initiate sequential matches.

Similar to the previous case, we can show that the feasible region in \eqref{eq:feasibleregion_simultaneous_onesided}  corresponds to the intersection of a partition with a laminar matroid.
\begin{definition}[Laminar Matroid]
A family $\mathcal{X}\subseteq 2^{\mathcal{E}}$ over a ground set of elements $\mathcal{E}$ is called laminar if for any $X,Y\in\mathcal{X}$ we either have $X\cap Y=\emptyset$, $X\subseteq Y$ or $X\subseteq Y$. Assume that for each element $u\in\mathcal{E}$ there exists some $A\in\mathcal{X}$ such that $A\ni u$. For each $A\in\mathcal{X}$ let $c(A)$ a positive integer associated with it. A laminar matroid $\mathcal{I}$ is defined as $\mathcal{I} = \{A\subseteq \mathcal{E}: \ |A\cap X|\leq c(X) \ \forall X\in\mathcal{X}\}$.
\end{definition}

\begin{lemma}\label{lemma:laminarmatroid}
The feasible region in~\eqref{eq:feasibleregion_generalsequential} corresponds to the intersection of a partition and a laminar matroid.
\end{lemma}

\begin{proof}[Proof.]
Our ground set of elements is $\mathcal{E} = \vec{E}^1_I \cup E^1$.
First, let us define our laminar family $\mathcal{X}$. For every pair $i\in I,j\in J$ consider $X_{i,j} = \{(i,j),\{i,j\}\}$, also for every $i\in I$ consider $Y_i=\{(i,\ell):\ell\in\potentials^1_i\}\cup\{\{i,\ell\}:\ell\in\potentials^1_i\}$. Indeed, this is a laminar family, two sets of type $X$ do not intersect and two sets of type $Y$ also do not intersect. Sets of type $X$ and $Y$ intersect only if they correspond to the same $i\in I$ in which case $X_{i,j}\subseteq Y_i$. Finally, for every $X_{i,j}$ we have $c(X_{i,j}) = 1$ and for each $Y_i$ we have $c(Y_i)=K_i$. Therefore, for this laminar family $\mathcal{X}$ and values $c(\cdot)$ we have that $\mathcal{I} = \{A\subseteq \mathcal{E}: \ |A\cap X|\leq c(X) \ \forall X\in\mathcal{X}\}$ coincides with feasible region $Q^1$. The partition matroid is defined by the following parts: for each $j\in J$ consider $\mathcal{E}_j=\{\{i,j\}: \ i\in\mathcal{P}_j\}$ and $\mathcal{E}_0 = \mathcal{E}\setminus \bigcup_{j\in J}\mathcal{E}_j$. The budget for each part $\mathcal{E}_j$ is $K_j$ and for $\mathcal{E}_0$ is $|\mathcal{E}_0|$. 
\end{proof}

Therefore, since the objective function is still monotone and submodular, we can use the local search algorithm proposed by~\cite{lee2009submodular} to find a feasible solution with an approximation factor of $1/(2+\epsilon)$, as we formalize in Proposition~\ref{prop:approximation_one_directional_nonseq}.

\begin{proposition}\label{prop:approximation_one_directional_nonseq}
When Problem~\ref{def:two_period_problem} 
is restricted to one-directional policies with non-sequential matches in the first period, there exists a feasible solution $\mathbf{x}^1,\mathbf{w}^1$ whose objective value is at least $\frac{1}{2+\epsilon}\cdot\OPT$, for any $\epsilon>0$.
\end{proposition}
\begin{proof}[Proof.]
Similar to the proof of Proposition~\ref{prop:approximation_two_directional_alg2}.
\end{proof}

\subsection{Two-directional Interactions and Non-Sequential Matches.}\label{subsec: two-directional simultaneous}
Under this setting, we note that the feasible region of Problem~\ref{def:two_period_problem}  corresponds to the intersection of three partition matroids. 
\begin{lemma}\label{lemma:feasibleregion_extendiblesystem}
The feasible region of Problem~\ref{def:two_period_problem}  corresponds to the intersection of three partition matroids when restricted to two-directional interactions allowing non-sequential matches in the first period.
\end{lemma}

\begin{proof}[Proof.]
Our ground set of elements is $\mathcal{E} = \vec{E}^1_I\cup\vec{E}^1_J \cup E^1$.
As in Lemma~\ref{lemma:intersection_2partitions}, the first partition consists of the following parts: $\mathcal{E}_i = \{\{(i,j),\{i,j\}\}: \ \text{for every} \ j\in\potentials^1_i\}$ for all $i \in I$ and $\mathcal{E}_0 =  \mathcal{E}\setminus\cup_{i\in I}\mathcal{E}_i$.
It is easy to check that $\mathcal{E} = \mathcal{E}_0\cup\bigcup_{i\in I}\mathcal{E}_\ell$ and $\mathcal{E}_\ell\cap \mathcal{E}_{\ell'}=\emptyset$ for every $\ell,\ell'\in I\cup\{0\}$ such that $\ell \neq \ell'$. Finally, the budget for each part $\mathcal{E}_i$ is $K_i$ for every $i\in I$ and $|\mathcal{E}_0|$ for $\mathcal{E}_0$. Analogously, we can construct the second partition matroid by considering parts: $\mathcal{E}_{i,j}= \{\{(j,i),\{i,j\}\}: \ \text{for every} \ i\in\potentials^1_j\}$ for all $j \in J$ and $\mathcal{E}_0 =  \mathcal{E}\setminus\cup_{j\in J}\mathcal{E}_j$.
The final partition is composed by the following parts: For every pair $i\in I$ and $j\in J$, we define a part $\mathcal{E}_{i,j}$ as the set $\{(i,j), (j,i), \{i,j\}\}$. Indeed this forms a partition of $\mathcal{E}$. Finally, the budget for each part $\mathcal{E}_{i,j}$ is 1.
\end{proof}

Hence, we can use the local search algorithm in~\cite{lee2009submodular} to find a feasible solution with an approximation factor of $1/(3+\epsilon)$, as we formalize in the next proposition.

\begin{proposition}\label{prop:approximation_two_directional_nonseq}
When Problem~\ref{def:two_period_problem} is restricted to two-directional policies with non-sequential matches in the first period, there exists a feasible solution $(\mathbf{x}^1,\mathbf{w}^1)$ whose objective value is at least $\frac{1}{3+\epsilon}\cdot\OPT$, for any $\epsilon>0$.
\end{proposition}
\begin{proof}[Proof.]
Similar to the proof of Proposition~\ref{prop:approximation_two_directional_alg2}.
\end{proof}
Note that the guarantee in the proposition above is worse than that of DH-int, however, the former leads to a polynomial-time algorithm.

\subsection{One-Directional Interactions and Non-Sequential Matches in Both Periods}\label{subsec: simultaneous matches in second period} 
As we show in Proposition~\ref{prop: next problem is not submodular in the general case}, even when $T=2$, the function that returns the optimal number of matches attained in the second period (from sequential and non-sequential interactions) given a family of realized arcs and edges $R$ is not submodular when we enable non-sequential matches in the second period. Thus, the analysis is significantly more complex as we can no longer rely on submodular optimization techniques and the correlation gap. Nevertheless, we can still provide a performance guarantee for the case when non-sequential matches are allowed in both periods and when like probabilities on the initiating side are small, as we formalize in the next theorem. 
\begin{theorem}\label{thm:limitcase_simul_both}
Consider the case with one-directional sequential matches starting from side \(I\) and non-sequential matches in both periods. Suppose that probabilities are time-independent, i.e., $\plike_{\ell\ell'}^1=\plike_{\ell\ell'}^2=\plike_{\ell\ell'}$ and that $\phi_{ij} \leq 1/n$ for any $i\in I, j\in J$, where $n = |I|$.
Denote by $\text{OPT}_1$ the optimal solution when $\Pi$ is restricted to policies with non-sequential matches only in the first period and by $\text{OPT}_2$ the optimal value when $\Pi$ allows for policies to implement non-sequential matches in both periods. Then, we have
\[
\text{OPT}_1\geq \left(\frac{1}{2e}-o(1)\right)\cdot \text{OPT}_2.
\]
More importantly, any $\gamma$-approximation algorithm for $\text{OPT}_1$ leads to a $\gamma\left(\frac{1}{2e}-o(1)\right)$ approximation for $\text{OPT}_2$.
\end{theorem}
Theorem~\ref{thm:limitcase_simul_both} implies that allowing non-sequential matches in the second period does not arbitrarily improve the optimal expected number of matches when the initiating side is sufficiently picky and the market is sufficiently large. In other words, the majority of matches come from sequential interactions. As a result, we conjecture that the guarantees in the Appendices above may extend to the more general case. We emphasize that analyses on large markets with small probabilities are common in the literature; see~\cite{mehta2012online,goyal2023online} for some examples from the online matching literature. 


To show Theorem~\ref{thm:limitcase_simul_both}, recall function $\Ma^2$ defined for $(\mathbf{x}^1,\mathbf{w}^1)$ as
\[
\Ma^2(\mathbf{x}^1,\mathbf{w}^1) = \sum_{R\subseteq \vec{E}^1} f(R)\cdot \PP_{\mathbf{x}^1}(\cR=R).
\]
This sum is over all possible arcs in $\vec{E}^1$. Note that $\PP_{\mathbf{x}^1}(\cR=R)$ ``restricts'' the valid sets $R$ that come from the first-period profiles $\mathbf{x}^1$, i.e., this probability will be zero for any $R$ that contains an element $(\ell,\ell')\in\vec{E}^1$ with $x^1_{\ell,\ell'}= 0$. This means that, effectively, we are summing over subsets of $\vec{E}^1(\mathbf{x}^1)=\{(\ell,\ell')\in\vec{E}^1: \ x^1_{\ell,\ell'}=1\}$.

\begin{lemma}\label{lemma:redefinition_Mf}
The function $\mathcal{M}^2$ can be reformulated as
\[
\Ma^2(\mathbf{x},\mathbf{w}) = \sum_{R\subseteq\vec{E}^1}f(R\cap\vec{E}^1(\mathbf{x}))\cdot \PP(\cR=R),
\]
where $\PP(\cR=R)=\prod_{e\in R}\phi^1_e\prod_{\vec{E}^1\setminus R}(1-\phi^1_e)$.
\end{lemma}
Note that in the result above the distribution of $\cR$ does not depend on the first-period decisions.
\begin{proof}[Proof of Lemma~\ref{lemma:redefinition_Mf}.]
First, note that for any $L\subseteq \vec{E}^1(\mathbf{x})$ and any $R\subseteq\vec{E}^1$ such that $L = R\cap \vec{E}^1(\mathbf{x})$ we have $f(R) = f(L)$. Then,
\begin{align*}
\Ma^2(\mathbf{x},\mathbf{w}) &= \sum_{R\subseteq\vec{E}^1}f(R\cap\vec{E}^1(\mathbf{x}))\cdot \PP(\cR=R)\\
& = \sum_{L\subseteq \vec{E}^1(\mathbf{x})}\sum_{\substack{R\subseteq\vec{E}^1:\\ R\cap\vec{E}^1(\mathbf{x})=L}}f(R\cap\vec{E}^1(\mathbf{x})) \cdot \PP(\cR=R)\\
& = \sum_{L\subseteq \vec{E}^1(\mathbf{x})}\sum_{\substack{R\subseteq\vec{E}^1:\\R\cap\vec{E}^1(\mathbf{x})=L}}f(L) \cdot \PP(\cR=R)\\
& = \sum_{L\subseteq \vec{E}^1(\mathbf{x})}f(L) \sum_{\substack{R\subseteq\vec{E}^1:\\R\cap\vec{E}^1(\mathbf{x})=L}} \PP(\cR=R)\\
& = \sum_{L\subseteq \vec{E}^1(\mathbf{x})}f(L) \prod_{e\in L}\phi^1_ex_e\prod_{e\in \vec{E}^1(\mathbf{x})\setminus L}(1-\phi^1_ex_e).
\end{align*}
\end{proof}

Let $\tilde{f}\lrp{R,\mathbf{x}^1, \mathbf{w}^1}$ be the solution of~\eqref{eq: general next problem} with $t+1=2$. Following a similar argument as in the proof of Lemma~\ref{lemma:redefinition_Mf}, we can define $\tilde{\Ma}^2\lrp{\mathbf{x}^1, \mathbf{w}^1}$ as
\[
\tilde{\Ma}^2(\mathbf{x}^1,\mathbf{w}^1) = \sum_{R\subseteq\vec{E}^1}\tilde{f}\lrp{R\cap\vec{E}^1(\mathbf{x}),\mathbf{x}^1, \mathbf{w}^1}\cdot \PP(\cR=R)
\]
Now, we are ready to prove our main result.
\begin{proof}[Proof of Theorem~\ref{thm:limitcase_simul_both}.]
For simplicity, we show the result for the case when $K_\ell = 1$ for all $\ell\in I\cup J$; for the general case, the analysis is similar since we can do it per pair of users.
Let $(\mathbf{x}^{1,\star},\mathbf{w}^{1,\star})$ and $(\tilde{\mathbf{x}}^{1,\star},\tilde{\mathbf{w}}^{1,\star})$ be the optimal solutions of the two-period version of Problem~\ref{def:two_period_problem} excluding and enabling non-sequential matches in the second period, leading to $\text{OPT}_1$ and $\text{OPT}_2$ expected matches, respectively.
To ease the exposition, we will drop super-indices of the variables.

Our goal is to lower bound the following ratio
\[
\frac{\text{OPT}_1}{\text{OPT}_2} = \frac{\sum_{e\in E^1} w_e\cdot \beta_e+ \Ma^2(\mathbf{x}, \mathbf{w})}{\sum_{e\in E^1} \tilde{w}_e\cdot \beta_e+ \tilde{\Ma}^2(\tilde{\mathbf{x}},\tilde{\mathbf{w}})}.
\]
By using Lemma~\ref{lemma:redefinition_Mf}, we redefine both objectives to consider backlog distributions that are independent of the decisions in the first period, i.e.,
\begin{equation}\label{eq:ratio_analysis}
\frac{\sum_{e\in E^1} w_e\cdot \beta_e+ \Ma^2(\mathbf{x}, \mathbf{w})}{\sum_{e\in E^1} \tilde{w}_e\cdot \beta_e+ \tilde{\Ma}^2(\tilde{\mathbf{x}},\tilde{\mathbf{w}})}
=  \frac{\sum_R\left(\sum_{e\in E^1} w_e \beta_e+ f(R\cap\vec{E}^1(\mathbf{x}))\right)\PP(\cR=R)}{\sum_R\left(\sum_{e\in E^1} \tilde{w}_e \cdot \beta_e
+ \tilde{f}(R\cap\vec{E}^1(\tilde{\mathbf{x}}), \tilde{\mathbf{x}},\tilde{\mathbf{w}})\right)\PP(\cR=R)}.
\end{equation}
\vspace{0.5em}
To lower bound this ratio, we now proceed to lower bound the expected contribution of each pair \((i,j) \in I\times J\). Given that we are assuming one directional interactions with sequential matches going from $I$ to $J$, the expected contribution of the pair $(i,j)$ in $\OPT_1$ depends on the first-period decisions:
\begin{itemize}
\item If $x_{i,j} = 1$, then the contribution of the pair $(i,j)$ in the numerator is 
\[
\Delta_{ij} = \phi_{ij}\phi_{ji}\sum_{R: (i,j)\notin R}\ind{i\in S_j(R,\mathbf{x})}\cdot\PP(\cR=R)
\]
where $S_j(R,\mathbf{x})$ is an optimal solution for user $j$ in objective $f$ of the second period for $R$ and $\mathbf{x}$. $\Delta_{ij}$ is the product of: the probability that $i$ liked $j$, the probability that $i$ was shown to $j$ in the second stage and the probability that $j$ liked $i$. Observe that if no other agent $i'\in I$ such that $x_{i',j} = 1$ liked $j$, then $i$ would be part of the optimal solution of the second stage. In other words, the event in which no one (except $i$) likes $j$ implies the event of $i$ being part of the optimal solution. Therefore,
\[
\Delta_{ij} = \phi_{ij}\phi_{ji}\sum_{R: (i,j)\notin R}\ind{i\in S_j(R,\mathbf{x})}\cdot\PP(\cR=R)\geq \phi_{ij}\phi_{ji}\prod_{\ell\neq i: x_{\ell,j}=1}(1-\phi_{\ell,j})\geq \phi_{ij}\phi_{ji}\left(1-\frac{1}{n}\right)^{n-1}
\]
where the last inequality is due to our assumption.
\item If $w_{i,j} = 1$, then the contribution of the pair $i,j$ in the numerator is $\phi_{ij}\phi_{ji}$.
\end{itemize}
Now, let us compare the contribution of $(i,j)$ to $\text{OPT}_1$ relative to its contribution to $\text{OPT}_2$.
As before, we have different cases depending on the solutions $\tilde{\mathbf{x}}$, $\tilde{\mathbf{w}}$, and their second-period responses:
\begin{itemize}
\item If $\tilde{x}_{i,j'} = 1$ for some $j'\in \potentials_i$ with $j'\neq j$ (when $j'=j$ is analogous). In this case, the platform shows $j'$ to $i$ in the first period (initiating a sequential interaction) instead of $j$, hoping to get an extra non-sequential match in the second period.
Therefore, the contribution of $\tilde{x}_{i,j'}$ to $\text{OPT}_2$ would potentially involve two terms.
First,  between $i$ and $j'$, we have
\[\tilde{\Delta}_{ij'}=\phi_{ij'}\phi_{j'i}\sum_{R: (i,j')\notin R}\ind{i\in S_{j'}(R,\tilde{\mathbf{x}}, \tilde{\mathbf{w}})}\cdot\PP(\cR=R)
\]
where $S_{j'}(R,\tilde{\mathbf{x}}, \tilde{\mathbf{w}})$ is an optimal solution for user $j'$ in the objective $\tilde{f}$ of the second period with set $R$, $\tilde{\mathbf{x}}$ and $\tilde{\mathbf{w}}$. In the worst case, there is also a non-sequential match between $i$ and $j$ in the second period, which contributes (in expectation) 
\[
\tilde{\Delta}_{ij}=\phi_{ij}\phi_{ji}\sum_{R}\ind{i\in S_{j}(R,\tilde{\mathbf{x}}, \tilde{\mathbf{w}})}\cdot\PP(\cR=R).
\]
Therefore, we have the following contribution in the denominator
\begin{align*}
\tilde{\Delta}_{ij'} +\tilde{\Delta}_{ij}&= \phi_{ij'}\phi_{j'i}\sum_{R: (i,j')\notin R}\ind{i\in S_{j'}(R,\tilde{\mathbf{x}}, \tilde{\mathbf{w}})}\cdot\PP(\cR=R) + \phi_{ij}\phi_{ji}\sum_{R}\ind{i\in S_{j}(R,\tilde{\mathbf{x}}, \tilde{\mathbf{w}})}\cdot\PP(\cR=R)\\&\leq 2\phi_{ij}\phi_{ji},
\end{align*}
where the inequality is due to two facts: (i) $\tilde{\Delta}_{ij'}\leq \phi_{ij}\phi_{ji}$, otherwise in solution $\mathbf{x}$ with value $\text{OPT}_1$ we could show $j'$ to $i$ instead of $j$  and obtain a better solution, which would contradict the optimality of $\mathbf{x}$; (ii) $\sum_{R}\ind{i\in S_{j}(R,\tilde{\mathbf{x}}, \tilde{\mathbf{w}})}\cdot\PP(\cR=R)\leq 1$.
\item If $\tilde{w}_{i,j'} = 1$ for some $j'\in\potentials_i$; $j'$ is potentially different than $j$. In this case, the platform decided to show $j'$ to $i$ simultaneously in the first period. Then, the contribution in the worst case is the following
\begin{align*}
\tilde{\Delta}_{ij'} +\tilde{\Delta}_{ij}=\phi_{ij'}\phi_{j'i} + \phi_{ij}\phi_{ji}\sum_{R}\ind{i\in S_{j}(R,\tilde{\mathbf{x}}, \tilde{\mathbf{w}})}\cdot\PP(\cR=R)\leq 2\phi_{ij}\phi_{ji}.
\end{align*}
where the inequality follows as before, i.e., we must have $\phi_{j'i}\leq \phi_{ji}$ (recall both $\phi_{ij'}$ and $\phi_{ij}$ are at most $1/n$), since otherwise 
we can change the solution in $\mathbf{x}$ for $\OPT_1$ and get a better objective as $j'$ would be part of the second-period optimal solution ($j'$ would see $i$ in response) whenever $j$ is.
\item If $\tilde{x}_{i,j} = \tilde{w}_{i,j}= 0$ for all $j\in\potentials_i$. This case is similar to the previous one, but now the contribution would be composed only by a second-period term $\tilde{\Delta}_{ij}$, which is no worse (in terms of denominator) than the other cases.
\end{itemize}
Therefore, for any pair $(i,j)$, the ratio between each contribution is at least
\[
\frac{\Delta_{ij}}{\tilde{\Delta}_{ij'}+\tilde{\Delta}_{ij}}\geq\frac{\phi_{ij}\phi_{ji}\left(1-\frac{1}{n}\right)^n}{2\phi_{ij}\phi_{ji}} \geq\frac{1}{2e}-o(1)
\]
Finally,
\[
\frac{\sum_R\left(\sum_{e\in E^1} w_e \beta_e+ f(R\cap\vec{E}^1(\mathbf{x}))\right)\PP(\cR=R)}{\sum_R\left(\sum_{e\in E^1} \tilde{w}_e \cdot \beta_e
+ \tilde{f}(R\cap\vec{E}^1(\tilde{\mathbf{x}}), \tilde{\mathbf{x}},\tilde{\mathbf{w}})\right)\PP(\cR=R)}\geq \min_{i,j,j'}\left\{\frac{\Delta_{ij}}{\tilde{\Delta}_{ij'}+\tilde{\Delta}_{ij}}\right\}\geq \frac{1}{2e}-o(1).
\]
Note that we may be comparing more terms than we actually need, but we were looking just for a lower bound.
To conclude the second part of the theorem, we observe that we can always obtain more matches (in expectation) when we allow non-sequential matches in both periods rather than in the first period only. Formally, for any feasible solution $(\mathbf{x}^{1},\mathbf{w}^{1})$ in Problem~\eqref{def:two_period_problem}, then
\[
\Ma^2(\mathbf{x}^{1},\mathbf{w}^{1}) \leq \tilde{\Ma}^2(\mathbf{x}^{1},\mathbf{w}^{1}).
\]
Therefore, if we consider $\mathbf{x}^1,\mathbf{w}^1$ a $\gamma$-approximate solution for Problem~\eqref{def:two_period_problem}. Then, we have
\[
\tilde{\Ma}^2(\mathbf{x}^{1},\mathbf{w}^{1}) + \sum_{e\in E^1}\beta_ew^1_e\geq \Ma^2(\mathbf{x}^{1},\mathbf{w}^1) + \sum_{e\in E^1}\beta_ew^1_e \geq\gamma\cdot\text{OPT}_1\geq \gamma\cdot\left(\frac{1}{2e}-o(1)\right)\cdot\text{OPT}_2. 
\]
\end{proof}

\section{Missing Proofs in Section~\ref{sec: analysis T period model}}\label{os:proofs_Tperiods}
\begin{proof}[Proof of Lemma~\ref{lemma: UB semi-adaptive}.]
Consider an optimal semi-adaptive policy $\pi^\star$. Since $\pi^\star$ decides non-adaptively the profiles that initiate sequential interactions and the profiles of non-sequential interactions, then define $x_{\ell,\ell'}=1$ if $\pi^\star$ shows profile $\ell'$ to $\ell$ in some stage (as the initiation of a sequential interaction), and zero otherwise. Similarly, let $w_e=1$ if the users in $e$ mutually see each other in some period, and zero otherwise. Let $\Omega$ be the set of sample paths. For every $\omega\in\Omega$, let $y_{\ell,\ell'}^\omega =1$ if the policy displays profile $\ell'$ to user $\ell$ from the backlog in some period in path $\omega$, and zero otherwise. Finally, define $y_{\ell,\ell'} = \sum_{\omega\in\Omega}y_{\ell,\ell'}^\omega\cdot\PP(\omega)$, where $\PP(\omega)$ is the probability of sample path $\omega$, which is completely determined by the like/dislike realizations.

Let us prove that these variables are feasible in Problem~\eqref{eq: problem DH multiple-periods}. Clearly, we have $x_{\ell,\ell'}+x_{\ell',\ell}+w_e\leq 1$ since the policy shows a profile either as a sequential interaction or as a non-sequential interaction. For every $\ell\in I\cup J$ and $\ell'\in\potentials^1_\ell$ we have
\begin{equation*}
    y_{\ell,\ell'} 
    = \sum_{\omega\in\Omega}y_{\ell,\ell'}^\omega\cdot\PP(\omega)
    \leq  \sum_{\omega\in\Omega}\ind{\ell'\text{ liked }\ell\text{ in }\omega}\cdot x_{\ell',\ell}\cdot\PP(\omega)
    =  x_{\ell',\ell}\cdot\sum_{\omega\in\Omega}\ind{\ell'\text{ liked }\ell\text{ in }\omega}\cdot\PP(\omega)
    \leq  \phi_{\ell',\ell}x_{\ell',\ell}.
\end{equation*}
Finally, we have that for every $\ell\in I\cup J$
\begin{align*}
\sum_{\ell'\in\potentials^1_\ell}x_{\ell,\ell'}+\sum_{e\in E^1:\ell\in e}w_e + \sum_{\ell'\in\potentials^1_\ell}y_{\ell,\ell'} &= \sum_{\omega\in\Omega}\Big[\sum_{\ell'\in\potentials^1_\ell}x_{\ell,\ell'}+\sum_{e\in E^1:\ell\in e}w_e + \sum_{\ell'\in\potentials^1_\ell}y_{\ell,\ell'}^\omega\Big]\cdot\PP(\omega)\\
&\leq K_\ell\cdot T\cdot\sum_{\omega\in\Omega}\PP(\omega) = K_\ell\cdot T,
\end{align*}
where in the inequality we used that, for every path, we show at most $K_\ell$ in every period and there are $T$ periods. 
\end{proof}

\begin{proof}[Proof of Theorem~\ref{theorem:approx_nonadaptive_policy}.]
The first part of the analysis is similar to that of Theorem~\ref{thm: guarantee for DH}, i.e., we compare the duals of two formulations and use the concept of correlation gap. With a slight abuse of notation, we will use the same notation as in the proof of Theorem~\ref{thm: guarantee for DH}. Given an optimal solution $(\mathbf{x}^\star,\mathbf{w}^\star,\mathbf{y}^\star)$ of Problem~\eqref{eq: problem DH multiple-periods}, we define
\begin{align}\label{eq:P1prime T periods}
F(\mathbf{x}^\star, \mathbf{w}^\star) := \max & \quad \sum_{\ell \in \cIJ} \sum_{\ell'\in \potentials_\ell^1}\phi_{\ell, \ell'} \cdot y_{\ell, \ell'}\\
s.t.&\quad y_{\ell, \ell'}\leq x^\star_{\ell', \ell}\cdot \phi_{\ell', \ell}, \hspace{6.8em} \forall \ell \in \cIJ, \;\ell' \in \potentials_\ell^1 \notag\\
&\quad \sum_{\ell' \in \potentials_\ell^1} y_{\ell, \ell'}\leq \tilde{K}^T_\ell, \hspace{2.5cm}\forall \ell \in \cIJ\notag\\
&\quad y_{\ell, \ell'}\geq 0,  \hspace{3.7cm}\forall \ell \in \cIJ, \ell' \in \potentials_\ell^1.\notag
\end{align}
where $\tilde{K}^T_\ell := \asize_\ell\cdot T - \sum_{\ell' \in \potentials_\ell^1} x_{\ell, \ell'}^\star - \sum_{e\in E: \ell \in e} w_e^\star$ for each $\ell \in \cIJ$. We now define our second formulation:
\begin{align}\label{eq:P2 T periods}
G(\mathbf{x}^\star, \mathbf{w}^\star) := \max & \quad \sum_{\ell \in \cIJ} \sum_{L \subseteq \potentials_\ell^1} f^{[T]}_\ell(L)\cdot\lambda_{\ell,L}\\
s.t.&\quad \sum_{B\subseteq \potentials_\ell^1}\lambda_{\ell,L} = 1 \hspace{3.75cm} \forall \ell \in \cIJ\notag\\
&\quad \sum_{L\subseteq \potentials_\ell^1} \lambda_{\ell,L} = \phi_{\ell', \ell}\cdot x^\star_{\ell', \ell} , \hspace{5em} \forall \ell \in \cIJ, \ell'\in \potentials_\ell^1 \notag\\
&\quad \lambda_{\ell,L}\geq 0, \hspace{4.5cm} \forall \ell \in \cIJ, \; L\subseteq \potentials_\ell^1.\notag
\end{align}  
where 
\begin{equation}\label{eq:f_Tperiods}
f^{[T]}_\ell(L) := \max\left\{\sum_{\ell'\in L}\phi_{\ell, \ell'}\cdot z_{\ell,\ell'}: \ \sum_{\ell'\in L} z_{\ell, \ell'}\leq \tilde{K}^T_\ell, \ z_{\ell, \ell'}\leq \ind{\ell'\in L}, \ z_{\ell, \ell'}\geq0\right\}.
\end{equation}

Consider an optimal solution $(\mathbf{x}^\star,\mathbf{w}^\star,\mathbf{y}^\star)$ of Problem~\eqref{eq: problem DH multiple-periods}. Let us denote by $\mathbf{L}=\{L_\ell\}_{\ell\in I\cup J}$ the random set of likes that result from showing profiles in $\mathbf{x}^\star$, i.e., 
\[
\mathbf{L} = \lrl{\ell': x^\star_{\ell,\ell'}=1, \  \Phi_{\ell, \ell'}^t = 1 \text{ in some period } t\in[T], \ \forall \ell\in I\cup J}.
\] Given the independence of users' decisions, the distribution of $\mathbf{L}$ is such that for every $\ell'\in I\cup J$
\[
\lambda^{\sf ind}_{\ell', L} = \prod_{\ell\in L}\phi_{\ell,\ell'}x^\star_{\ell,\ell'}\prod_{\ell\notin L}(1-\phi_{\ell,\ell'}x^\star_{\ell,\ell'}).
\]
Since the function $f_\ell^{[T]}$ defined in \eqref{eq:f_Tperiods} is monotone and submodular, then \citep{agrawal2010correlation} shows that
\begin{equation}\label{eq:correlation_gap_Tperiods}
\frac{\sum_{\ell\in I\cup J} \sum_{L\subseteq\potentials^1_\ell}f^{[T]}_\ell(L)\cdot\lambda^{\sf ind}_{\ell, L}}{\sum_{\ell\in I\cup J} \sum_{L\subseteq\potentials^1_\ell}f^{[T]}_\ell(L)\cdot\lambda^{\star}_{\ell, L}}\geq 1-\frac{1}{e},
\end{equation}
where $\lambda^\star$ is an optimal solution of Problem~\eqref{eq:P2 T periods}. Note that
\[
\sum_{\ell\in I\cup J} \sum_{L\subseteq\potentials^1_\ell}f^{[T]}_\ell(L)\cdot\lambda^{\star}_{\ell, L} = G(\mathbf{x}^\star,\mathbf{w}^\star) \geq F(\mathbf{x}^\star,\mathbf{w}^\star),
\]
where the inequality can be shown equivalently as Lemma~\ref{lemma: dual distribution problem is upper bound}. 
Therefore, we can conclude that 
\begin{equation}\label{eq:aux_semiadaptive}
\sum_{\ell\in I\cup J} \sum_{L\subseteq\potentials^1_\ell}f^{[T]}_\ell(L)\cdot\lambda^{\sf ind}_{\ell, L}\geq (1-1/e)\cdot F(\mathbf{x}^\star,\mathbf{w}^\star).
\end{equation}
The final step in our proof is to show that the expected number of sequential matches achieved by Algorithm~\ref{alg: multiple-period DH} is lower bounded by the term on the left in \eqref{eq:aux_semiadaptive}.

Recall the definitions in Algorithm~\ref{alg: multiple-period DH}, i.e., $X_\ell = \{\ell'\in\potentials^1_\ell: \; x_{\ell,\ell'}^\star = 1\}$ and $W_\ell=\{e\in E^1: \; \ell\in e, \; w_e^\star = 1\}$ for all $\ell\in I\cup J$. Note that Algorithm~\ref{alg: multiple-period DH} exhausts all the profiles in $X_\ell\cup W_\ell$ for every $\ell\in I\cup J$. Given this, we can first conclude that the expected number of non-sequential matches achieved by our method is
\(
\sum_{e\in E^1}\beta_e w^\star_e
\),
which coincides with the non-sequential part of the optimal objective in Problem~\eqref{eq: problem DH multiple-periods}.

We now focus on analyzing the expected number of sequential matches obtained by our method. For every $t\in[T]$, let $X^t_{\ell} = X_\ell\cap S^t_\ell$ be the set of profiles displayed to $\ell$ in period $t$ as the initiating side of a sequential interaction. Similarly, let $W^t_\ell = W_\ell\cap S_\ell^t$ be the set of profiles shown to $\ell$ in period $t$ as part of non-sequential interaction. Consider a sample path of Algorithm~\ref{alg: multiple-period DH} and fix the realization of $\Phi_{\ell,\ell'}^t$ for every $\ell'\in X^t_{\ell}$ and $t\in[T]$. This implies that the backlog in each period $t\in [T]$ for each user $\ell\in I\cup J$ is a deterministic set $B^t_\ell$. Denote by $Z^{t}_\ell$ the set of profiles shown from the backlog $B_\ell^t$ of user $\ell$ in stage $t$; note that $Z^t_\ell$ is deterministic (as $B_\ell^t$  is) and results from solving the problem in Step~\ref{alg-step: solving f} in Algorithm~\ref{alg: multiple-period DH}.
Finally, as introduced in Section~\ref{sec: model}, let $L_\ell^t$ be the set of users that liked $\ell$ in period $t$, i.e., $L_\ell^t = \lrl{\ell': \ell \in X^t_{\ell'}, \Phi^t_{\ell', \ell} = 1}$.
In a slight abuse of notation, let $L_\ell = \bigcup_{t\in[T]} L_\ell^t$ be the total set of profiles that liked $\ell$ during $T$ periods. Note that, for any $t\in [T]$, the backlog satisfies $B_\ell^t = \bigcup_{\tau < t} L_\ell^\tau\setminus \lrp{\bigcup_{\tau < t} Z_\ell^\tau}$. 

For every $\ell\in I\cup J$, the total expected number of sequential matches achieved by Algorithm~\ref{alg: multiple-period DH} is
\begin{equation}\label{eq:aux_semiadaptive2}
\sum_{t\in[T]} f^t_{\ell}(B_\ell^t) = \sum_{t\in[T]}\sum_{\ell'\in Z_\ell^t}\phi_{\ell',\ell}
\end{equation}
where $f^t_\ell$ is as defined analogous to~\eqref{eq:second_stage_problem} considering $K_\ell - |X_\ell^t| - |W_\ell^t|$ as the right-hand side of the capacity constraint. Since we include profiles in each $S^t_\ell$ from $X_\ell$ in decreasing order of $\phi_{\ell',\ell}$ (Step~\ref{alg-step: selecting profiles to initiate interactions with highest probs}), then $\bigcup_{t\in[T]}Z_\ell^t$ contains the profiles with the highest values $\phi_{\ell',\ell}$ available in $L_\ell$. This is crucial when we compare our method with the value of Problem~\eqref{eq:f_Tperiods} for the same set of likes $L_\ell$.

For the total set of likes $L_\ell$ that $\ell$ received, let $\tilde{Z}_\ell$ be an optimal set in problem $f^{[T]}_\ell(L_\ell)$ as defined in \eqref{eq:f_Tperiods}. Construct a partition of $\tilde{Z}_\ell$ in $T$ sets $\tilde{Z}^1_\ell,\ldots,\tilde{Z}^T_\ell$ as follows: (i) profiles in $\tilde{Z}_\ell$ are added to the sets in decreasing order of values $\phi_{\ell,\ell'}$, (ii) the set $\tilde{Z}^t_\ell$ is filled before continuing to $\tilde{Z}^{t+1}_\ell$, (iii) for each $t\in[T]$ we have $|\tilde{Z}^t_\ell|\leq K_\ell - |X_\ell^t| - |W_\ell^t|$ and a profile $\ell'$ is added to $\tilde{Z}^{t}_\ell$ only if was not added before and $\ell'\in B_\ell^t$. This is possible because $L_\ell = \bigcup_{t\in[T]} B_\ell^t$ and
\[
\sum_{t\in[T]}|\tilde{Z}^t_\ell|\leq \sum_{t\in[T]}K_\ell - |X_\ell^t| - |W_\ell^t| = K_\ell\cdot T - \sum_{t\in[T]} |X_\ell^t| - \sum_{t\in[T]}|W_\ell^t| = K_\ell\cdot T - \sum_{\ell'\in\potentials^1_\ell} x^\star_{\ell,\ell'} - \sum_{e\in E:\ell\in e}w^\star_e = \tilde{K}^T_\ell,
\]
where $\tilde{K}^T_\ell$ is the right-hand side of the cardinality constraint in $f^{[T]}_\ell(L_\ell)$. Note that the second equality follows as the profiles in $X_\ell\cup W_\ell$ are exhausted.
\begin{claim}
For every $\ell\in I\cup J$, $t\in[T]$ we have $Z_\ell^t = \tilde{Z}^t_\ell$.
\end{claim}
\begin{proof}[Proof of Claim.]
First, note that $Z_\ell^t = \emptyset$ if and only if $B_\ell^t=\emptyset$ or $K_\ell = |X_\ell^t| +|W_\ell^t| = |S^t_\ell|$, which means that $\tilde{Z}^t_\ell = \emptyset$ due to the construction above. 
Second, suppose that $Z_\ell^t \neq \emptyset$ which means that $B_\ell^t\neq\emptyset$ and $K_\ell > |X_\ell^t| +|W_\ell^t|$. Since $Z_\ell^t$ is an optimal set for the problem defined by $f^t_\ell(B_\ell^t)$, then $Z_\ell^t$ contains the profiles with highest values $\phi_{\ell,\ell'}$ that are available in $B^t_\ell$. Due to Step~\ref{alg-step: selecting profiles to initiate interactions with highest probs} of Algorithm~\ref{alg: multiple-period DH}, the profiles in  $Z_\ell^t$ have higher values than profiles in $Z_\ell^{t+1}$. Therefore, because of our construction was over the same set of likes $L_\ell$ then the sets $\tilde{Z}_\ell^t$ must coincide with $Z_\ell^t$.  
\end{proof}
The above claim implies that the total expected number of sequential matches is such that
\[
\sum_{t\in[T]} f^t_{\ell}(B_\ell^t) = \sum_{t\in[T]}\sum_{\ell'\in Z_\ell^t}\phi_{\ell,\ell'} =\sum_{t\in[T]}\sum_{\ell'\in \tilde{Z}_\ell^t}\phi_{\ell,\ell'} = \sum_{\ell'\in \tilde{Z}_\ell}\phi_{\ell,\ell'} = f^{[T]}_\ell(L_\ell).
\]
Taking expectation over the set of likes $L_\ell$ and summing over $\ell\in I\cup J$ we obtain
\[
\E\left[\sum_{\ell\in I\cup J}\sum_{t\in[T]} f^t_{\ell}(B_\ell^t)\right] = \sum_{\ell\in I\cup J}\sum_{L_\ell\subseteq\potentials^1_\ell}f^{[T]}_\ell(L_\ell)\cdot\lambda^{\sf ind}_{\ell, L_\ell} \geq (1-1/e)\cdot F(\mathbf{x}^\star,\mathbf{w}^\star)
\]
By summing non-sequential and sequential matches, we get
\begin{align*}
\sum_{e\in E^1}\beta_ew^\star_e+\E\left[\sum_{\ell\in I\cup J}\sum_{t\in[T]} f^t_{\ell}(B_\ell^t)\right] &\geq \sum_{e\in E^1}\beta_ew^\star_e + (1-1/e)\cdot F(\mathbf{x}^\star,\mathbf{w}^\star) \\
&\geq (1-1/e)\cdot\left[\sum_{e\in E^1}\beta_ew^\star_e + F(\mathbf{x}^\star,\mathbf{w}^\star)\right]\\
&\geq (1-1/e)\cdot \OPT_{\eqref{eq: problem DH multiple-periods}}\\
&\geq (1-1/e)\cdot \OPT^\Pi,
\end{align*}
where $\OPT_{\eqref{eq: problem DH multiple-periods}}$ denotes the optimal value of Problem~\eqref{eq: problem DH multiple-periods} and $\OPT^\Pi$ is the optimal value achieved by the best semi-adaptive policy. In the last inequality, we use Lemma~\ref{lemma: UB semi-adaptive}.
\end{proof}

\begin{proof}[Proof of Theorem~\ref{theorem:approx_nonadaptive_policy2}.]
Fix the initiating side to be $I$. Since we restrict the space of policies to one-directional interactions and sequential matches, then our main relaxation is Problem~\eqref{eq: problem DH multiple-periods} with $x_{j,i}=0$ for all $j\in J$, $i\in\potentials_j^1$ (i.e., users in $J$ do not initiate interactions), $y_{i,j}=0$ for all $i\in I$ and $j\in\potentials^1_i$ (i.e., users in $I$ are not followers in an interaction), $w_e= 0$ for all $e\in E^1$ (i.e., there are no non-sequential interactions) and relaxing the integrality for the remaining variables, i.e., $x_{i,j}\in[0,1]$ for all $i\in I$, $j\in\potentials^1_i$. The resulting formulation is
\begin{equation}\label{eq: problem DH multiple-periods2}
\begin{split}
    \max \quad & \; \sum_{i\in I}\sum_{j\in \potentials^1_i} y_{j,i} \cdot \plike_{j,i} \\
    \text{s.t.} \quad 
& \; y_{j, i} \leq x_{i, j} \cdot \plike^{1}_{i, j}, \qquad \forall i \in I, \; j \in \potentials^1_i, \\     
& \; \sum_{j \in \potentials^1_i} x_{i, j}  \leq K_i\cdot T, \quad \forall i \in I, \\      
& \; \sum_{i \in \potentials^1_j} y_{j, i}  \leq K_j\cdot T, \quad \forall j \in J, \\ 
& \; \mathbf{x} \in [0,1]^{\vec{E}^1_I}, \; \mathbf{y} \in [0,1]^{\vec{E}^1_J}
\end{split}
\end{equation}

Our main method is similar to Algorithm~\ref{alg: multiple-period DH}, but we include a rounding procedure since an optimal solution $(\mathbf{x}^\star,\mathbf{y}^\star)$ of Problem~\eqref{eq: problem DH multiple-periods2} may be fractional. Specifically, we use the dependent randomized rounding procedure introduced in \citep{gandhi06}, as we did to prove Proposition~\ref{prop:approximation_one_directional_alg1} (see Appendix~\ref{subsec: one-directional sequential}). This method outputs a random binary vector $\tilde{\mathbf{x}}\in\{0,1\}^{\vec{E}^1_I}$ such that: (i) $\sum_{j \in \potentials^1_i} \tilde{x}_{i, j}  \leq K_i\cdot T$ for all $i\in I$ with probability 1; (ii) $\E[\tilde{x}_{i,j}] = x^\star_{i,j}$ for all $i\in I$, $j\in\potentials^1_j$.
Finally, based on $\tilde{x}$, we proceed as in Algortihm~\ref{alg: multiple-period DH} to decide the subsets of profiles to show in each period.

Fix the random binary vector $\tilde{\mathbf{x}}$ obtained from the rounding. The proof from here is analogous to that of Theorem~\ref{theorem:approx_nonadaptive_policy}. Therefore, we obtain that the total expected number of sequential matches achieved by our method is such that
\[
\E\left[\sum_{j\in J}\sum_{t\in[T]} f^t_{j}(B_j^t)\right] = \sum_{j\in  J}\sum_{L_j\subseteq\potentials^1_j}f^{[T]}_j(L_j)\cdot\lambda^{\sf ind}_{j, L_j} 
\]
where $f^{[T]}_j(L_j)$ is defined as in \eqref{eq:f_Tperiods} but with $\tilde{K}^T_j = K_j\cdot T$, the expectation is over the randomness of likes/dislikes and 
\[
\lambda^{\sf ind}_{j,L_j} = \prod_{i\in L_j}\phi_{i,j}\tilde{x}_{i,j} \prod_{i\notin L_j}(1-\phi_{i,j}\tilde{x}_{i,j}).
\]
By taking expectation over the randomness of $\tilde{\mathbf{x}}$, we obtain
\begin{align*}
\E_{\tilde{\mathbf{x}}\sim\mathbf{x}^\star,\mathbf{B}}\left[\sum_{j\in J}\sum_{t\in[T]} f^t_{j}(B_j^t)\right] &= \sum_{j\in  J}\sum_{L_j\subseteq\potentials^1_j}f^{[T]}_j(L_j)\cdot \E_{\tilde{\mathbf{x}}\sim\mathbf{x}^\star}[\lambda^{\sf ind}_{j, L_j}]\\
&=\sum_{j\in  J}\sum_{L_j\subseteq\potentials^1_j}f^{[T]}_j(L_j)\cdot\prod_{i\in L_j}\phi_{i,j}x^\star_{i,j} \prod_{i\notin L_j}(1-\phi_{i,j}x^\star_{i,j}).
\end{align*}
where in the last equality we used the independence of like/dislike decisions of each user $i$ and $\E_{\tilde{\mathbf{x}}\sim\mathbf{x}^\star}[\tilde{x}_{i,j}] = x^\star_{i,j}$ which is the property of the rounding algorithm.

To finish the proof, note that thanks to the correlation gap we have
\[
\sum_{j\in  J}\sum_{L_j\subseteq\potentials^1_j}f^{[T]}_j(L_j)\cdot\prod_{i\in L_j}\phi_{i,j}x^\star_{i,j} \prod_{i\notin L_j}(1-\phi_{i,j}x^\star_{i,j}) \geq (1-1/e)\cdot F({\mathbf{x}^\star}), 
\]
where $F(\mathbf{x}^\star)$ is \eqref{eq:P1prime T periods} with $\tilde{K}^T_j = K_j\cdot T$ and $y_{i,j} = 0$ for all $i\in I$, $j\in\potentials^1_i$. Finally, $F({\mathbf{x}^\star})\geq \OPT^\Pi$ follows from Lemma~\ref{lemma:UB_Tperiods_onedirectional} (see below), where $\OPT^\Pi$ is the optimal value achieved by an adaptive policy with one-directional interactions.
\end{proof}




\begin{lemma}\label{lemma:UB_Tperiods_onedirectional}
Problem~\eqref{eq: problem DH multiple-periods2} is a relaxation of Problem~\ref{def: general problem with multiple periods} for any adaptive policy with one-directional sequential interactions.
\end{lemma}

\begin{proof}[Proof.]
Consider an optimal adaptive policy $\pi^\star$ with one-directional interactions. Let $\Omega$ be the set of sample paths. For every $\omega\in\Omega$, let $x_{i,j}^\omega =1$ if the policy displays profile $j$ to user $i$ in some period in path $\omega$, and zero otherwise, and let $x_{i,j} = \sum_{\omega\in\Omega}x_{i,j}^\omega\cdot\PP(\omega)$, where $\PP(\omega)$ is the probability of sample path $\omega$. Similarly, let $y_{i,j}^\omega =1$ if the policy displayed profile $i$ to user $j$ from the backlog in some period in path $\omega$, and zero otherwise.  Finally, define $y_{j,i} = \sum_{\omega\in\Omega}y_{j,i}^\omega\cdot\PP(\omega)$. Let us prove that these variables are feasible in Problem~\eqref{eq: problem DH multiple-periods2}.

First, for every $i\in I$ we have
\[
\sum_{j\in\potentials^1_i}x_{i,j}=\sum_{\omega\in\Omega}\sum_{j\in\potentials^1_i}x_{i,j}^\omega\cdot\PP(\omega)\leq K_i\cdot T\cdot \sum_{\omega\in\Omega}\PP(\omega) = K_i\cdot T,
\]
where the first inequality is because in every sample path we have $\sum_{j\in\potentials^1_i}x_{i,j}^\omega\leq K_i\cdot T$ as there are $T$ periods and at most $K_i$ profiles are shown per period. In a similar fashion, we can prove that $\sum_{i\in\potentials^1_j}y_{j,i}\leq K_j\cdot T$.

Finally, for every $j\in J$ and $i\in\potentials^1_j$ we have
\begin{align*}
y_{j,i} &= \sum_{\omega\in\Omega}y_{j,i}^\omega\cdot\PP(\omega)\\
&\leq  \sum_{\omega\in\Omega}\ind{i\text{ liked }j\text{ in }\omega}\cdot x_{i,j}^\omega\cdot\PP(\omega)\\
&\leq  \phi_{i,j}\sum_{\omega\in\Omega_{\text{-}ij}}x_{i,j}^\omega\cdot\PP_{\text{-}ij}(\omega)\\
&=  \phi_{i,j}x_{i,j},
\end{align*}
where in the first inequality we use that $y_{j,i}^\omega\leq \ind{i\text{ liked }j\text{ in }\omega}\cdot x_{i,j}^\omega$. In the following inequality, we remove the independent event that $i$ liked $j$ from the sample path space $\Omega_{\text{-}ij}$ and from the probability $\PP_{\text{-}ij}(\omega)$. In the last equality, we use that the decision of $\pi^\star$ to show $j$ to $i$ is independent of the realization of $\Phi_{i,j}$ since the policy doesn't have access to realizations beforehand. 
\end{proof}

\section{Missing Details in Section~\ref{sec: experiments}}\label{ec:experiments}

To estimate the probability that each user \(i\) likes a profile \(j\in \potentials_i^t\), we use a logit model with user-fixed effects:
\begin{equation}\label{eq: model specification preliminary evidence}
y_{ijt} = \alpha_i + \lambda_t + X_{i,j}' \beta + \epsilon_{ijt}.
\end{equation}
The dependent variable, \(y_{ijt}\), is a latent variable that is related to the like decision \(\Phi_{i,j}^t\) according to
\[
\Phi_{i,j}^t = \begin{cases}
1 & \text{if } y_{ijt} > 0, \\
0 & \text{otherwise.}
\end{cases}
\]
We control for users' unobserved heterogeneity
by including user fixed-effects, \(\alpha_i\). We also control for period-dependent unobservables by including period fixed-effects, \(\lambda_t\). The third term on
the right-hand side, \(X_{i,j}' \beta\), controls for observable characteristics of profile
\(j\), and also for their interaction with user \(i\)'s observable characteristics.
Specifically, we control for the attractiveness score, age, height and education level (measured in a scale from 1 to 3) of the profile evaluated. In addition, for each of these variables we control for the square of the positive and negative difference between the value for the user evaluating and that of the profile evaluated. Finally, we also control for whether the users share the same race and religion.
Finally, \(\epsilon_{ijt}\) is an idiosyncratic shock that follows an extreme value distribution.
In Table~\ref{tab: estimation results logit with fixed effects}, we report the
estimation results.
\begin{table}[htpb]
\caption{Estimation Results}\label{tab: estimation results logit with fixed effects}
\centerline{\scalebox{0.75}{\begin{tabular}{lcc}
\toprule
& (1)            & (2) \\
\midrule
Batch size                  & -0.004$^{***}$ & -0.004$^{***}$\\
& (0.0004)       & (0.0004)\\
Score                      & 0.832$^{***}$  & 0.832$^{***}$\\
& (0.014)        & (0.014)\\
Score - Positive difference                       & 0.012$^{***}$  & 0.012$^{***}$\\
& (0.003)        & (0.003)\\
Score - Negative difference                       & -0.011$^{***}$ & -0.011$^{***}$\\
& (0.003)        & (0.003)\\
Age                        & -0.026$^{***}$ & -0.026$^{***}$\\
& (0.004)        & (0.004)\\
Age - Positive difference                         & -0.002$^{***}$ & -0.002$^{***}$\\
& (0.0004)       & (0.0004)\\
Age - Negative difference                         & -0.0009$^{**}$ & -0.0009$^{**}$\\
& (0.0005)       & (0.0005)\\
Height                     & 0.055$^{***}$  & 0.055$^{***}$\\
& (0.008)        & (0.008)\\
Height - Positive difference                      & 0.002$^{***}$  & 0.002$^{***}$\\
& (0.0006)       & (0.0006)\\
Height - Negative difference                      & -0.004$^{***}$ & -0.004$^{***}$\\
& (0.0006)       & (0.0006)\\
Education level                        & 0.060$^{**}$   & 0.061$^{**}$\\
& (0.024)        & (0.024)\\
Education level - Positive difference                         & -0.001         & -0.0008\\
& (0.014)        & (0.014)\\
Education level - Negative difference                         & -0.077$^{***}$ & -0.077$^{***}$\\
& (0.016)        & (0.016)\\
Share religion                   & 0.078$^{***}$  & 0.078$^{***}$\\
& (0.013)        & (0.013)\\
Share race                  & 0.457$^{***}$  & 0.458$^{***}$\\
& (0.030)        & (0.030)\\
&                &  \\
\midrule
User  & $\checkmark$  & $\checkmark$\\
Date        &                & $\checkmark$\\
\midrule
Observations                      & 396,226        & 396,226\\
Pseudo R$^2$                      & 0.386        & 0.386\\
\bottomrule
\end{tabular}}}
\end{table}

\begin{figure}[htp!]
\centering
\caption{Distribution of Like Probabilities}\label{fig: distribution like probs}
\includegraphics[scale=0.5]{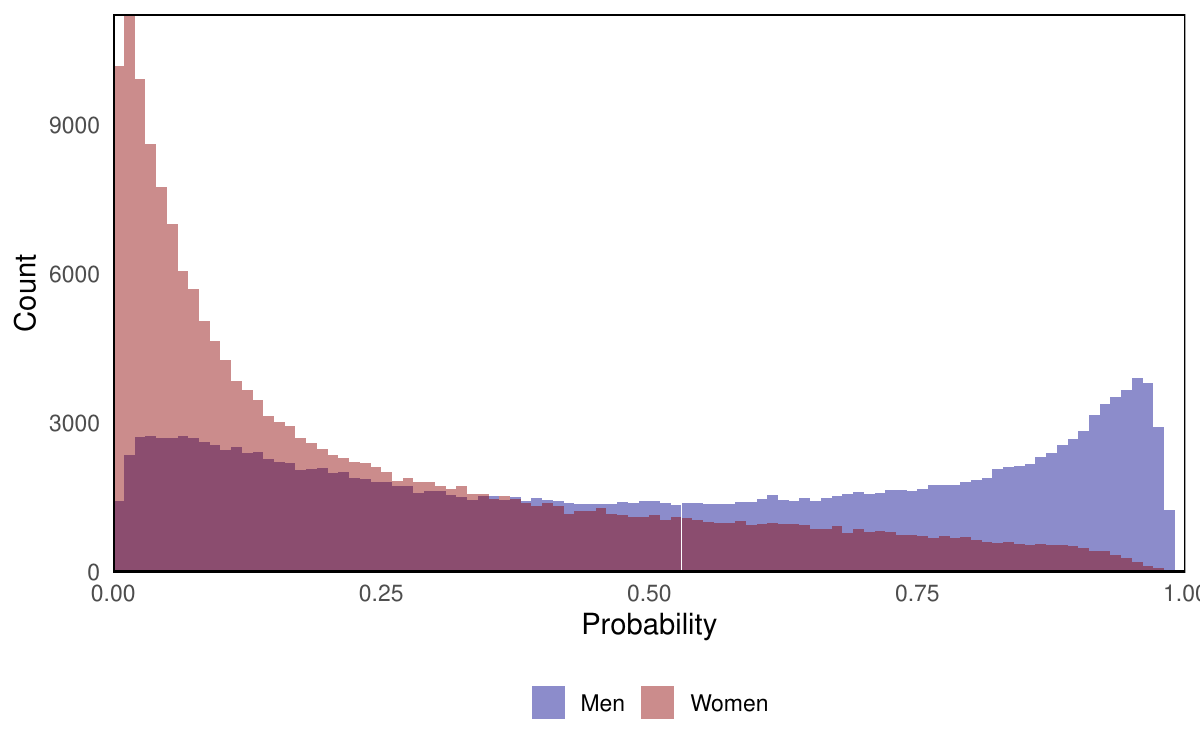}
\end{figure}

Using these coefficients, for each user \(i\) and profile \(j\in \potentials_i^1\) we compute the probability
\(\plike_{i,j}\). In Figure~\ref{fig: distribution like probs} we plot the distribution of like probabilities separating by gender, estimated using the parameters from column (2) in Table~\ref{tab: estimation results logit with fixed effects}. These are the probabilities we use on our simulation study.

In Table~\ref{tab: descriptives instance}, we report summary statistics for the sample used in the numerical experiments (standard deviations in parenthesis).
\begin{table}[htp!]
\caption{Descriptives of Instance}\label{tab: descriptives instance}
\centerline{\scalebox{0.85}{\begin{tabular}{lccccc}
\toprule
& \(N\) & Score & Potentials & Backlog & Like Prob. \\
\midrule
Women & 1682 & 4.462 & 109.895 & 0.120 & 0.295\\
&  & (0.056) & (2.860) & (1.040) & (0.005)\\
Men & 1193 & 2.563 & 133.477 & 0.029 & 0.527\\
&  & (0.044) & (3.591) & (0.385) & (0.005)\\    
\bottomrule
\end{tabular}}}
\end{table}

\begin{table}[htp!] \caption{Overall Results - Linear}\label{tab:overall_results}
\centerline{\scalebox{0.85}{
\begin{tabular}{llcccccc}
\toprule
 & & DH-int & DH & DHT & Greedy & PM & Partner \\
\midrule
\multirow{2}{*}{Two-Directional} & \multirow{2}{*}{Overall}  & 4691.14 & 4634.00 & 4441.99 & 3947.85 & 3675.61 & 3597.18 \\
& & (58.81) & (48.62) & (56.83) & (57.20) & (52.22) & (55.73) \\
\midrule
\multirow{4}{*}{One-Directional} 
& \multirow{2}{*}{Men}  & 4246.17 & 4161.99 & 3749.95 & 3421.24 & 3671.00 & 3085.63 \\
& & (53.99) & (58.79) & (49.88) & (45.46) & (53.33) & (49.34) \\
& \multirow{2}{*}{Women}  & 4250.23 & 4198.66 & 3951.47 & 3559.55 & 3674.48 & 3258.73 \\
& & (58.05) & (47.81) & (52.37) & (47.37) & (53.87) & (49.80) \\
\bottomrule
\end{tabular}}}
\vspace{0.5em}
\begin{minipage}{\linewidth}
\footnotesize \emph{Note:} We report the average number of matches produced assuming linear history effect. Standard deviations are reported in parentheses.
\end{minipage}
\end{table}

%% file: 0_main.bbl
\begin{thebibliography}{39}
\providecommand{\natexlab}[1]{#1}
\providecommand{\url}[1]{\texttt{#1}}
\expandafter\ifx\csname urlstyle\endcsname\relax
  \providecommand{\doi}[1]{doi: #1}\else
  \providecommand{\doi}{doi: \begingroup \urlstyle{rm}\Url}\fi

\bibitem[Adamczyk(2011)]{adamczyk2011improved}
Marek Adamczyk.
\newblock Improved analysis of the greedy algorithm for stochastic matching.
\newblock \emph{Information Processing Letters}, 111\penalty0 (15):\penalty0
  731--737, 2011.

\bibitem[Agrawal et~al.(2010)Agrawal, Ding, Saberi, and
  Ye]{agrawal2010correlation}
Shipra Agrawal, Yichuan Ding, Amin Saberi, and Yinyu Ye.
\newblock Correlation robust stochastic optimization.
\newblock In \emph{Proceedings of the twenty-first annual ACM-SIAM symposium on
  Discrete Algorithms}, pages 1087--1096. SIAM, 2010.

\bibitem[Aouad and Saban(2022)]{aouad21}
Ali Aouad and Daniela Saban.
\newblock Online assortment optimization for two-sided matching platforms.
\newblock \emph{Management Science}, 0\penalty0 (0), 2022.

\bibitem[Arnosti and Shi(2020)]{Arnosti20}
Nick Arnosti and Peng Shi.
\newblock Design of lotteries and wait-lists for affordable housing allocation.
\newblock \emph{Management Science}, 66\penalty0 (6):\penalty0 2291--2307,
  2020.

\bibitem[Ashlagi et~al.(2022)Ashlagi, Krishnaswamy, Makhijani, Saban, and
  Shiragur]{ashlagi19}
Itai Ashlagi, Anilesh Krishnaswamy, Rahul Makhijani, Daniela Saban, and
  Kirankumar Shiragur.
\newblock Technical note - assortment planning for two-sided sequential
  matching markets.
\newblock \emph{Operations Research}, 70\penalty0 (5):\penalty0 2784--2803,
  2022.

\bibitem[Bansal et~al.(2012)Bansal, Gupta, Li, Mestre, Nagarajan, and
  Rudra]{bansal2012lp}
Nikhil Bansal, Anupam Gupta, Jian Li, Juli{\'a}n Mestre, Viswanath Nagarajan,
  and Atri Rudra.
\newblock When lp is the cure for your matching woes: Improved bounds for
  stochastic matchings.
\newblock \emph{Algorithmica}, 63:\penalty0 733--762, 2012.

\bibitem[Besbes et~al.(2021)Besbes, Castro, and Lobel]{besbes21}
Omar Besbes, Francisco Castro, and Ilan Lobel.
\newblock Surge pricing and its spatial supply response.
\newblock \emph{Management Science}, 67\penalty0 (3):\penalty0 1350--1367,
  2021.

\bibitem[Besbes et~al.(2023)Besbes, Fonseca, Lobel, and Zheng]{fonseca23}
Omar Besbes, Yuri Fonseca, Ilan Lobel, and Fanyin Zheng.
\newblock Signaling competition in two-sided markets.
\newblock 2023.

\bibitem[Celdir et~al.(2024)Celdir, Cho, and Hwang]{celdir24}
Musa~Eren Celdir, Soo-Haeng Cho, and Elina~H. Hwang.
\newblock Popularity bias in online dating platforms: Theory and empirical
  evidence.
\newblock \emph{Manufacturing \& Service Operations Management}, 26\penalty0
  (2):\penalty0 537--553, 2024.

\bibitem[Chen et~al.(2009)Chen, Immorlica, Karlin, Mahdian, and Rudra]{Chen09}
Ning Chen, Nicole Immorlica, Anna~R. Karlin, Mohammad Mahdian, and Atri" Rudra.
\newblock Approximating matches made in heaven.
\newblock In \emph{Automata, Languages and Programming}, pages 266--278.
  Springer Berlin Heidelberg, 2009.

\bibitem[Chen and Nasiry(2020)]{Chen20}
Ningyuan Chen and Javad Nasiry.
\newblock Does loss aversion preclude price variation?
\newblock \emph{Manufacturing \& Service Operations Management}, 22\penalty0
  (2):\penalty0 383--395, 2020.

\bibitem[Cui and Hamilton(2024)]{cui24}
Titing Cui and Michael Hamilton.
\newblock Pricing strategies for online dating platforms, 2024.

\bibitem[Fisher et~al.(1978)Fisher, Nemhauser, and Wolsey]{fisher78}
Marshall~L Fisher, George~L Nemhauser, and Laurence Wolsey.
\newblock An analysis of approximations for maximizing submodular set
  functions--ii.
\newblock \emph{Polyhedral combinatorics}, 1978.

\bibitem[Gamlath et~al.(2019)Gamlath, Kale, and Svensson]{gamlath2019beating}
Buddhima Gamlath, Sagar Kale, and Ola Svensson.
\newblock Beating greedy for stochastic bipartite matching.
\newblock In \emph{Proceedings of the Thirtieth Annual ACM-SIAM Symposium on
  Discrete Algorithms}, pages 2841--2854. SIAM, 2019.

\bibitem[Gandhi et~al.(2006)Gandhi, Khuller, Parthasarathy, and
  Srinivasan]{gandhi06}
Rajiv Gandhi, Samir Khuller, Srinivasan Parthasarathy, and Aravind Srinivasan.
\newblock Dependent rounding and its applications to approximation algorithms.
\newblock \emph{Journal of the ACM (JACM)}, 53\penalty0 (3):\penalty0 324--360,
  2006.

\bibitem[Goyal and Udwani(2023)]{goyal2023online}
Vineet Goyal and Rajan Udwani.
\newblock Online matching with stochastic rewards: Optimal competitive ratio
  via path-based formulation.
\newblock \emph{Operations Research}, 71\penalty0 (2):\penalty0 563--580, 2023.

\bibitem[Guo et~al.(2025)Guo, Jiang, and Shen]{guo25}
Mengzi~Amy Guo, Hansheng Jiang, and Zuo-Jun~Max Shen.
\newblock Multiproduct dynamic pricing with reference effects under logit
  demand.
\newblock \emph{Manufacturing \& Service Operations Management}, 27\penalty0
  (5):\penalty0 1645--1663, 2025.

\bibitem[Halaburda et~al.(2018)Halaburda, Piskorski, and
  Yildirim]{Halaburda2018}
Hanna Halaburda, Mikolaj~Jan Piskorski, and Pinar Yildirim.
\newblock {Competing by Restricting Choice: The Case of Search Platforms}.
\newblock \emph{Management Science}, 64\penalty0 (8):\penalty0 3574--3594,
  2018.

\bibitem[He et~al.(2025)He, Zhang, and Zheng]{He25}
Taotao He, Yating Zhang, and Huan Zheng.
\newblock Assortment optimization under history-dependent effects, 2025.

\bibitem[Hu and Nasiry(2018)]{Hu18}
Zhenyu Hu and Javad Nasiry.
\newblock Are markets with loss-averse consumers more sensitive to losses?
\newblock \emph{Management Science}, 64\penalty0 (3):\penalty0 1384--1395,
  2018.

\bibitem[Immorlica et~al.(2022)Immorlica, Lucier, Manshadi, and
  Wei]{Immorlica22}
Nicole Immorlica, Brendan Lucier, Vahideh Manshadi, and Alexander Wei.
\newblock Designing approximately optimal search on matching platforms.
\newblock \emph{Management Science}, Forthcoming, 2022.

\bibitem[Jeloudar et~al.(2021)Jeloudar, Lo, Pollner, and Saberi]{Jeloudar21}
Mobin Jeloudar, Irene Lo, Tristan Pollner, and Amin Saberi.
\newblock Decentralized matching in a probabilistic environment.
\newblock In \emph{Proceedings of the 22nd ACM Conference on Economics and
  Computation}, pages 635--653, 2021.

\bibitem[Kagan et~al.(2025)Kagan, Leider, and Sahin]{kagan25}
Evgeny Kagan, Stephen Leider, and Ozge Sahin.
\newblock Sequential decision making: From decision elicitation to strategy
  identification.
\newblock \emph{Management Science}, (forthcoming), 2025.

\bibitem[Kanoria and Saban(2021)]{Kanoria2017}
Yash Kanoria and Daniela Saban.
\newblock Facilitating the search for partners on matching platforms.
\newblock \emph{Management Science}, 67\penalty0 (10):\penalty0 5990--6029,
  2021.

\bibitem[Karp et~al.(1990)Karp, Vazirani, and Vazirani]{kvv}
R.M. Karp, U.V. Vazirani, and V.V Vazirani.
\newblock An optimal algorithm for online bipartite matching.
\newblock In \emph{Proceedings of the 22nd Annual ACM Symposium on Theory of
  Computinh}, 1990.

\bibitem[K{\"o}k et~al.(2015)K{\"o}k, Fisher, and Vaidyanathan]{kok2015}
A.~K{\"o}k, M.~Fisher, and R.~Vaidyanathan.
\newblock \emph{Retail Supply Chain Management: Quantitative Models and
  Empirical Studies}, pages 175--236.
\newblock Springer US, 2015.

\bibitem[Lee et~al.(2009)Lee, Sviridenko, and Vondr{\'a}k]{lee2009submodular}
Jon Lee, Maxim Sviridenko, and Jan Vondr{\'a}k.
\newblock Submodular maximization over multiple matroids via generalized
  exchange properties.
\newblock In \emph{International Workshop on Approximation Algorithms for
  Combinatorial Optimization}, pages 244--257. Springer, 2009.

\bibitem[Lee and Niederle(2014)]{Lee2014}
Soohyung Lee and Muriel Niederle.
\newblock {Propose with a rose? Signaling in internet dating markets}.
\newblock \emph{Experimental Economics}, 18\penalty0 (4):\penalty0 731--755,
  2014.

\bibitem[Long et~al.(2020)Long, Nasiry, and Wu]{long20}
Xiaoyang Long, Javad Nasiry, and Yaozhong Wu.
\newblock A behavioral study on abandonment decisions in multistage projects.
\newblock \emph{Management Science}, 66\penalty0 (5):\penalty0 1999--2016,
  2020.

\bibitem[Manshadi and Rodilitz(2022)]{manshadi22}
Vahideh Manshadi and Scott Rodilitz.
\newblock Online policies for efficient volunteer crowdsourcing.
\newblock \emph{Management Science}, 68\penalty0 (9):\penalty0 6572--6590,
  2022.

\bibitem[Manshadi et~al.(2022)Manshadi, Rodilitz, Saban, and
  Suresh]{rodilitz22}
Vahideh Manshadi, Scott Rodilitz, Daniela Saban, and Akshaya Suresh.
\newblock Online algorithms for matching platforms with multi-channel traffic.
\newblock 2022.

\bibitem[Mehta and Panigrahi(2012)]{mehta2012online}
Aranyak Mehta and Debmalya Panigrahi.
\newblock Online matching with stochastic rewards.
\newblock In \emph{2012 IEEE 53rd Annual Symposium on Foundations of Computer
  Science}, pages 728--737. IEEE, 2012.

\bibitem[\"{O}zer and Zheng(2016)]{Ozer15}
\"{O}zalp \"{O}zer and Yanchong Zheng.
\newblock Markdown or everyday low price? the role of behavioral motives.
\newblock \emph{Management Science}, 62\penalty0 (2):\penalty0 326--346, 2016.

\bibitem[Rios and Ghosh(2024)]{rios24}
Ignacio Rios and Pramit Ghosh.
\newblock Competition in optimal stopping: Behavioral insights.
\newblock \emph{Manufacturing \& Service Operations Management}, 26\penalty0
  (6):\penalty0 2256--2273, 2024.

\bibitem[Rios et~al.(2023)Rios, Saban, and Zheng]{rios2021}
Ignacio Rios, Daniela Saban, and Fanyin Zheng.
\newblock Improving match rates in dating markets through assortment
  optimization.
\newblock \emph{Manufacturing {\&} Service Operations Management}, 25\penalty0
  (4):\penalty0 1304--1323, 2023.

\bibitem[Rochet and Tirole(2003)]{Rochet2003}
Jean-Charles Rochet and Jean Tirole.
\newblock {Two-Sided Markets}.
\newblock \emph{Journal of the European Economic Association}, pages 990--1029,
  2003.

\bibitem[Torrico et~al.(2021)Torrico, Carvalho, and Lodi]{torrico21}
Alfredo Torrico, Margarida Carvalho, and Andrea Lodi.
\newblock Multi-agent assortment optimization in sequential matching markets.
\newblock 2021.

\bibitem[Vondr\'ak(2008)]{vondrak08}
Jan Vondr\'ak.
\newblock Optimal approximation for the submodular welfare problem in the value
  oracle model.
\newblock In \emph{Proceedings of the 40th Annual ACM Symposium on the Theory
  of Computing (STOC)}, pages 67--74, 2008.

\bibitem[Wang(2018)]{Wang18}
Ruxian Wang.
\newblock When prospect theory meets consumer choice models: Assortment and
  pricing management with reference prices.
\newblock \emph{Manufacturing \& Service Operations Management}, 20\penalty0
  (3):\penalty0 583--600, 2018.

\end{thebibliography}
